\documentclass[11pt,twoside]{article}

\usepackage[ansinew]{inputenc}
\usepackage{amsfonts}
\usepackage{amssymb,amsmath}
\usepackage{latexsym}

\newcommand {\debeq}	{\begin{eqnarray*}}
\newcommand {\fineq}	{\end{eqnarray*}}
\newcommand {\lbd}	{\lambda}
\newcommand     {\eps}  {\epsilon}
\newcommand     {\vareps}       {\varepsilon}

\newcommand	{\expf}	[1]
{\mbox{e}^{#1}}

\newcommand	{\tendinfty}	
{\rightarrow\infty}

\newcommand	{\intgen}	
{\int_0^\infty}

\newcommand	{\ZZ}{\mathbb{Z}}

\newcommand	{\NN}{\mathbb{N}}
\newcommand	{\DoubleS}{\mathbb{S}}

\newcommand	{\TT}{\mathbb{T}}

\newcommand	{\UU}{\mathbb{U}}
\newcommand	{\CC}{\mathbb{C}}

\newcommand	{\RR}{\mathbb{R}}
\newcommand	{\PP}{\mathbb{P}}
\newcommand	{\EE}{\mathbb{E}}

\newtheorem	{thm}		{Theorem}[section]

\newtheorem	{dfn}		[thm]{Definition}
\newtheorem	{lem} 	[thm]	{Lemma}

\newtheorem     {rem}           {Remark}
\newtheorem	{prop}	[thm]{Proposition}
\newtheorem	{cor}		[thm]{Corollary}

\newcommand	{\indic}	[1]
{{\bf{1}}_{\{#1\}}}
\newcommand	{\indicbis}	[1]
{{\bf{1}}_{#1}}

\newcommand	{\Xtau}			{X^{(\tau)}}
\newcommand	{\Xtaut}			{X_t^{(\tau)}}
\newcommand	{\Xtauf}	[1]	{X_{#1}^{(\tau)}}

\linethickness{0.4pt}

\begin{document}

\title{The contour of splitting trees is a Lévy process}
\author{\textsc{By Amaury Lambert%\footnote{Université Pierre et Marie Curie and \'{E}cole Normale Supérieure, Paris}
}
}
\date{}\maketitle
\noindent\textsc{Laboratoire de Probabilités et Modèles Aléatoires\\
UMR 7599 CNRS and UPMC Univ Paris 06\\
Case courrier 188\\
4, Place Jussieu\\
F-75252 Paris Cedex 05, France}\\
\textsc{E-mail: }amaury.lambert@upmc.fr\\
\textsc{URL: }http://ecologie.snv.jussieu.fr/amaury/

\begin{abstract}
\noindent
Splitting trees are those random trees where individuals give birth at constant rate during a lifetime with general distribution, to i.i.d. copies of themselves.
The width process of a splitting tree is then a binary, homogeneous  Crump--Mode--Jagers (CMJ) process, and is not Markovian unless the lifetime distribution is exponential (or a Dirac mass at $\{\infty\}$). 
Here, we allow the birth rate to be infinite, that is, pairs of birth times and lifespans of newborns form a Poisson point process along the lifetime of their mother, with possibly infinite intensity measure.% $\Lambda$ on $(0,\infty]$ called the \emph{lifespan measure}, which has $\intgen (1\wedge r)\Lambda (dr)<\infty$.

A splitting tree is a random (so-called) chronological tree. Each element of a chronological tree is a (so-called) existence point $(v,\tau)$ of some individual $v$ (vertex) in a discrete tree, where $\tau$ is a nonnegative real number called chronological level (time). We introduce a total order on existence points, called linear order, and a mapping $\varphi$ from the tree into the real line which preserves this order. The inverse of $\varphi$ is called the exploration process, and the projection of this inverse on chronological levels the contour process.
  
For splitting trees truncated up to level $\tau$, we prove that thus defined contour process is a Lévy process reflected below $\tau$ and killed upon hitting 0. This allows to derive  properties of (i) splitting trees: conceptual proof of Le Gall--Le Jan's theorem in the finite variation case, exceptional points, coalescent point process, age distribution; (ii) CMJ processes: one-dimensional marginals, conditionings, limit theorems, asymptotic numbers of individuals with infinite vs finite descendances.
\end{abstract}  	
\medskip
%Received December 2006.\\
\textit{Running head.} The contour of splitting trees.

\section{Introduction}
A splitting tree \cite{GK}  is a tree formed by individuals with i.i.d. lifespans, who give birth at the same constant rate, while they are alive, to copies of themselves. The number of individuals alive at time $t$ evolves with $t$  according to non-Markovian dynamics, known as the Crump--Mode--Jagers (CMJ) process. Actually, general CMJ processes feature birth intensities that may vary through the lifetime of an individual as well as clutches of random sizes. Thus, the CMJ processes considered here are actually \emph{binary homogeneous CMJ processes}. On the other hand, we consider splitting trees that are more general than those considered in \cite{GK}, since here individuals may have infinitely many offspring. This is done by assuming that for each individual, the birth times and lifespans of her offspring form a Poisson point process with intensity $dt\,\Lambda(dr)$, where $\Lambda$ is a Lévy measure on $(0,\infty]$ called the \emph{lifespan measure}, which has $\intgen (1\wedge r)\Lambda (dr)<\infty$.\\
\\  
Contour processes of splitting trees as defined here are also different from those in \cite{GK}.
Our contour process  $(X_s, s\ge 0)$ visits once and once only all the instants at which any individual of the tree is alive, that we call  \emph{existence points} (the set of existence points of a given individual is merely its lifetime, which is in bijection with an interval of the real line). To avoid confusion with the usual denomination of `time' given to the index of stochastic processes (viz. the contour process in the present setting), we will call \emph{level}, or \emph{chronological level}, the real, physical time in which live the individuals of the tree. With the same goal, we will also try and give a Greek letter (such as $\tau$ or $\sigma$) to these levels, and a Latin letter (such as $t$ or $s$) to time as index of stochastic processes. In the case when each individual has finitely many offspring, one can set a rule for the contour as follows: when the contour process is about  to visit a birth level, it jumps to the death level of the newborn and then decreases linearly (at speed $-1$) along the lifetime of this individual until it encounters a new birth event. When the contour process ends up its visit of an individual, its value is thus the birth level of this individual. It then continues the visit of its mother's lifetime at the level where (when) it had been interrupted. Hereafter, this process  will be called \textit{jumping chronological contour process}, abbreviated as JCCP.\\
\\
The key result of the present work is that thus defined contour process $X$ for splitting trees is a Lévy process. Our inspiration comes from a previous study \cite{P} in the critical case with exponential lifespans (the only case when the CMJ process is Markovian, except the Yule case, where lifespans are a.s. infinite).

Our result yields a new interpretation of a now famous connection between Lévy processes and branching processes \cite{LGLJ, DLG, L2} that can be worded as follows: the genealogy of (continuous-state) Markov branching processes can be coded by a scalar Markov process (called genealogy-coding process in \cite{L2}), which is a Lévy process in the subcritical and critical cases. This process is usually an abstract object that has to be considered \emph{ex nihilo}, whereas here, it is defined as the contour process of a predefined tree (another interpretation was given in terms of queues in \cite{LGLJ}). 

This interpretation of a spectrally positive Lévy process as a contour process living in the space of chronological levels of a splitting tree yields an elegant way of considering and inferring properties related to the individuals alive at a fixed level $\tau$ (number,  coalescence levels, ages,...) in the tree.\\
In addition, it relates for the first time the genealogies of continuous-state branching processes defined in an apparently unrelated way in \cite{LGLJ} and \cite{BLG}, in the finite variation case. In the former work, the starting object is a Lévy process $X$ with no negative jumps starting from $\chi$ and killed when it hits 0, and to each timepoint $t$ is associated a height $H_t$ (generation, integer number), \label{3} which is given by the following functional of the path of $X$ 
\begin{equation}
\label{eqn : hauteur}
H_t:=\mathrm{Card}\{0\le s\le t : X_{s-}<\inf_{s\le r\le t} X_r^{}\}.
\end{equation}
The Lebesgue measure of the set of timepoints with height $n$ is a nonnegative real number $Z_n$, and $(Z_n;n\ge 0)$ is shown to satisfy both the Markov property and the branching property. \label{2} Here and elsewhere, we will say that a stochastic process $Z$ with nonnegative values satisfies the \emph{branching property}, if for any two independent copies $Z'$ and $Z''$ of $Z$ started respectively at $x$ and $y$, $Z'+Z''$ has the law of $Z$ started at $x+y$.\\
In the latter work, the starting object is a sequence of i.i.d. subordinators with zero drift, $S_1, S_2,\ldots$, and the population size at generation $n$ is $Z'_n:=S_n\circ \cdots\circ S_1 (\chi)$, with $Z'_0=\chi$. The genealogy of this continuous population is defined as follows: a point $c$ in generation $n$ ($c\in[0,Z'_{n}]$) is the daughter of a point $b$ in generation $n-1$ ($b\in[0,Z'_{n-1}]$), if $S_n(b-)<c<S_n(b)$. In particular, all points of generation $n-1$ that are no jump times of $S_n$ have no descendance.

Actually, these two genealogies can be coupled simultaneously starting from a single random object, namely, a splitting tree $\TT$: the aforementioned Lévy process $X$ is actually the \emph{JCCP} (jumping chronological contour process, as defined previously) of $\TT$, and $Z_n=Z'_n$ is the \emph{sum of all lifespans} of individuals belonging to generation $n$. Let us try and explain this briefly. \label{explication} First, at time $t$, the JCCP $X$ visits (the existence point at real time $X_t$ of) an individual whose generation in the discrete tree is the integer $H_t$, which can be seen to be exactly as in \eqref{eqn : hauteur} (see forthcoming Corollary \ref{cor : height process}). Second, the time $Z_n$ spent by the contour $X$ at height $n$, which is the total time spent by the height process $H$ at $n$, is exactly the Lebesgue measure of the set of existence points, in the splitting tree, of individuals of generation $n$, which is exactly the \emph{sum of their lifespans}. Third, one can embed the existence points of individuals of generation $n$ into the real half-line by arranging lifetimes `end to end' on the interval $[0,Z_n]$. Then by construction of the splitting tree, $Z_{n+1}$ is the sum of all atoms of i.i.d. Poisson measures with intensity $\Lambda$ defined on lifetime intervals whose lengths sum up to $Z_n$. As a consequence, conditional on $Z_n=z$, we have obtained that $Z_{n+1}= S_{n+1}(z)$, where the $(S_i;i\ge 1)$ are independent subordinators with zero drift and Lévy measure $\Lambda$ (see also the proof of Theorem \ref{thm : Jirina}).  This explains why $Z=Z'$, and that the genealogy defined thanks to these subordinators is exactly that associated with the topology of the splitting tree.\\
In his seminal article \cite{Ji}, \label{5} M. Jirina introduced for the first time (multidimensional) Markov processes with continuous-state space satisfying the branching property. He studied the whole class of such processes in discrete time, but only a subclass of them in continuous time (\emph{pure-jump} processes). For this reason, and because the term `Jirina process' for continuous-state branching processes in continuous time, has progressively disappeared in the last fifteen years (in favor of `CSBP', or `CB-process'), I propose to call \emph{Jirina processes} those branching processes in discrete time and continuous state-space like $(Z_n;n\ge 0)$.

Last, to clear up the difference between CMJ processes and Markov branching processes such as Jirina or Bienaymé--Galton--Watson (BGW) processes, note that both count the `population size' as time runs, but the former count the number of individuals alive at the \emph{same chronological level}, whereas the latter count those indexed by the successive \emph{generations} of the tree. Replacing time by generations guarantees their Markov property to branching processes, but makes their genealogy harder to infer \cite{O', L3}.\\
\\
In the next section, we set up classical notation on discrete trees and define plane chronological trees, that we endow with a genealogical structure, a distance, a linear order, a Lebesgue measure, and a closure. We also recall well-known facts about spectrally positive Lévy processes as well as Jirina processes. 

In Section 3, we introduce and study an order-preserving bijection $\varphi$ between the closure of a chronological tree with finite Lebesgue measure, and a compact interval of the real line. The inverse $\varphi^{-1}$ of this bijection is called the \emph{exploration process}. It is its projection on chronological levels $p_2 \circ \varphi^{-1}$, which is called \emph{jumping chronological contour process}, or JCCP.

It is only in Section 4 that we consider random chronological trees, called \emph{splitting trees} (as defined earlier). We study the properties of a splitting tree, and prove that the JCCP of its  truncation up to level $\tau$ is a Lévy process reflected below $\tau$ and killed upon hitting 0.

This last result allows to derive a number of properties of splitting trees and CMJ processes, which is done in the last section. For splitting trees, an intuitive proof of Le Gall--Le Jan's theorem is given in the finite variation case (cf. above);
the set $\Gamma$ of levels where the population size is infinite is considered: when $\Lambda$ is finite, $\Gamma$ is empty, and when $\Lambda$ is infinite, $\Gamma$ has zero Lebesgue measure but is everywhere dense a.s.; the coalescence levels of individuals alive at the same level are shown to be independent and with the same distribution, that we specify; the law of ages and residual lifetimes of individuals from a same level is computed. For CMJ processes, the one-dimensional marginal is shown to be \emph{modified geometric}; \label{9} the supercritical CMJ process conditioned on extinction is characterized; various limit theorems are given, among which the convergence in distribution, conditional on non-extinction, of the numbers of individuals with infinite vs finite descendances to $(p\xi, (1-p) \xi)$, where $\xi$ is an exponential random variable with parameter $p=1-\intgen r e^{-\eta r} \Lambda(dr)$, $\eta$ being the Malthusian parameter.\\
Actually, the \label{63}  set of points with infinite descendance, or skeleton, has a discrete branching structure (that of a Yule tree). In a forthcoming work \cite{Lprep2}, we extend the study of splitting trees with finite variation to splitting trees with infinite variation, and prove in particular that the branching structure of the skeleton is again discrete (see also \cite{BFM}). Note that in the infinite variation case, lifespans are not even locally summable, and that in the presence of a Brownian component, we cannot stick to the tree space we deal with here. We also point out that an application of the present work to allelic partitions is available \cite{Lprep1}.\\
\\
The reader who might like to avoid technicalities can proceed as follows: in Section 2, just focus on the definitions and the two statements; restrict the reading of Section 3 to Theorem \ref{thm : dfn exploration}, Definition \ref{dfn : jccp} and Theorem \ref{thm : allure jccp}; skip proofs in the last two sections.

\section{Preliminaries on trees and stochastic processes}

\subsection{Discrete trees}
Let $\NN$ denote the set of positive integers.
Locally finite rooted trees \cite{N, DLG} can be coded thanks to the so-called Ulam--Harris--Neveu labelling. Each vertex of the tree is represented by a finite sequence of integers as follows.  The root of the tree is $\emptyset$, the $j$-th child of $u=(u_1,\ldots, u_n)\in\NN^n$, is $uj$, where $vw$ stands for the concatenation of the sequences $v$ and $w$, here $uj=(u_1,\ldots, u_n, j)$. Set $\vert u \vert=n$ the \emph{generation}, or \emph{genealogical height}, of $u$. More rigorously, let
$$
{\cal U} = \bigcup_{n=0}^{\infty} \NN^n
$$
where $\NN^0 = \{\emptyset\}$. A \emph{discrete tree} $\cal T$ is a subset of $\cal U$ such that\\
\indent
(i) $\emptyset\in \cal T$\\
\indent
(ii)
if $v=uj\in{\cal T}$, where $j\in \NN$, then $u\in\cal T$\\
\indent
(iii)
for every $u\in\cal T$, there is a nonnegative integer $K_u\le \infty$ (the \emph{offspring number}) such that $uj\in\cal T$ if and only if $j\in\{1,\ldots, K_u\}$.\\
\\
Note that individuals can have infinitely (but countably) many offspring.
We write $u\prec v$ if $u$ is an \emph{ancestor} of $v$, that is, there is a sequence $w$ such that $v=uw$.
For any $u=(u_1,\ldots, u_n)$, $u\vert k$ denotes the ancestor $(u_1,\ldots, u_k)$ of $u$ at generation $k$.
We denote by $u\wedge v$ the \emph{most recent common ancestor}, in short \emph{mrca}, of $u$ and $v$, that is, the sequence $w$ with highest generation such that $w\prec u$ and $w\prec v$.

\subsection{Chronological trees}
Chronological trees are particular instances of $\RR$-trees. For further reading on $\RR$-trees, see e.g. \cite{E, DT, LGSurvey} and the references therein.
The $\RR$-trees we consider here can roughly be seen as the set of edges of some discrete tree embedded in the plane, where each edge length is a \emph{lifespan}.

\begin{figure}[ht]

\unitlength 1.8mm % = 4.55pt
\linethickness{0.4pt}
\ifx\plotpoint\undefined\newsavebox{\plotpoint}\fi % GNUPLOT compatibility
\begin{picture}(67.13,34)(-8,0)
\thicklines
\put(7,4){\line(0,1){12}}
\thinlines
%\vector{dash}(7,13)(10,13)
\put(10,13){\vector(1,0){.04}}\put(6.96,12.96){\line(1,0){.75}}
\put(8.46,12.96){\line(1,0){.75}}
%\end
\thicklines
\put(10,13){\line(0,1){8}}
\thinlines
%\vector{dash}(10,19)(13,19)
\put(13,19){\vector(1,0){.04}}\put(9.96,18.96){\line(1,0){.75}}
\put(11.46,18.96){\line(1,0){.75}}
%\end
\thicklines
\put(13,19){\line(0,1){5}}
\thinlines
%\vector{dash}(10,17)(15,17)
\put(15,17){\vector(1,0){.04}}\put(9.96,16.96){\line(1,0){.833}}
\put(11.62,16.96){\line(1,0){.833}}
\put(13.29,16.96){\line(1,0){.833}}
%\end
\thicklines
\put(15,17){\line(0,1){4}}
\put(18,18){\line(0,1){5}}
\put(24,8){\line(0,1){6}}
\thinlines
%\vector{dash}(24,12)(27,12)
\put(27,12){\vector(1,0){.04}}\put(23.96,11.96){\line(1,0){.75}}
\put(25.46,11.96){\line(1,0){.75}}
%\end
\thicklines
\put(27,12){\line(0,1){15}}
\thinlines
%\vector{dash}(27,25)(30,25)
\put(30,25){\vector(1,0){.04}}\put(26.96,24.96){\line(1,0){.75}}
\put(28.46,24.96){\line(1,0){.75}}
%\end
\thicklines
\put(30,25){\line(0,1){3}}
\thinlines
%\vector{dash}(15,18)(18,18)
\put(18,18){\vector(1,0){.04}}\put(14.96,17.96){\line(1,0){.75}}
\put(16.46,17.96){\line(1,0){.75}}
%\end
%\vector{dash}(30,26)(33,26)
\put(33,26){\vector(1,0){.04}}\put(29.96,25.96){\line(1,0){.75}}
\put(31.46,25.96){\line(1,0){.75}}
%\end
\thicklines
\put(33,26){\line(0,1){4}}
\thinlines
%\vector{dash}(33,29)(36,29)
\put(36,29){\vector(1,0){.04}}\put(32.96,28.96){\line(1,0){.75}}
\put(34.46,28.96){\line(1,0){.75}}
%\end
\thicklines
\put(36,29){\line(0,1){3}}
\thinlines
%\vector{dash}(33,28)(38,28)
\put(38,28){\vector(1,0){.04}}\put(32.96,27.96){\line(1,0){.833}}
\put(34.62,27.96){\line(1,0){.833}}
\put(36.29,27.96){\line(1,0){.833}}
%\end
\thicklines
\put(38,28){\line(0,1){2}}
\thinlines
%\vector{dash}(27,21)(35,21)
\put(35,21){\vector(1,0){.04}}\put(26.96,20.96){\line(1,0){.889}}
\put(28.73,20.96){\line(1,0){.889}}
\put(30.51,20.96){\line(1,0){.889}}
\put(32.29,20.96){\line(1,0){.889}}
\put(34.07,20.96){\line(1,0){.889}}
%\end
\thicklines
\put(35,21){\line(0,1){4}}
\thinlines
%\vector{dash}(27,15)(40,15)
\put(40,15){\vector(1,0){.04}}\put(26.96,14.96){\line(1,0){.929}}
\put(28.81,14.96){\line(1,0){.929}}
\put(30.67,14.96){\line(1,0){.929}}
\put(32.53,14.96){\line(1,0){.929}}
\put(34.38,14.96){\line(1,0){.929}}
\put(36.24,14.96){\line(1,0){.929}}
\put(38.1,14.96){\line(1,0){.929}}
%\end
\thicklines
\put(40,15){\line(0,1){8}}
\thinlines
%\vector{dash}(40,21)(43,21)
\put(43,21){\vector(1,0){.04}}\put(39.96,20.96){\line(1,0){.75}}
\put(41.46,20.96){\line(1,0){.75}}
%\end
\thicklines
\put(43,21){\line(0,1){4}}
\thinlines
%\vector{dash}(40,19)(45,19)
\put(45,19){\vector(1,0){.04}}\put(39.96,18.96){\line(1,0){.833}}
\put(41.62,18.96){\line(1,0){.833}}
\put(43.29,18.96){\line(1,0){.833}}
%\end
\thicklines
\put(45,19){\line(0,1){4}}
\thinlines
%\vector{dash}(45,22)(48,22)
\put(48,22){\vector(1,0){.04}}\put(44.96,21.96){\line(1,0){.75}}
\put(46.46,21.96){\line(1,0){.75}}
%\end
\thicklines
\put(48,22){\line(0,1){7}}
\thinlines
%\vector{dash}(48,25)(53,25)
\put(53,25){\vector(1,0){.04}}\put(47.96,24.96){\line(1,0){.833}}
\put(49.62,24.96){\line(1,0){.833}}
\put(51.29,24.96){\line(1,0){.833}}
%\end
%\vector{dash}(48,23)(55,23)
\put(55,23){\vector(1,0){.04}}\put(47.96,22.96){\line(1,0){.875}}
\put(49.71,22.96){\line(1,0){.875}}
\put(51.46,22.96){\line(1,0){.875}}
\put(53.21,22.96){\line(1,0){.875}}
%\end
\thicklines
\put(55,23){\line(0,1){3}}
\thinlines
%\vector{dash}(40,16)(57,16)
\put(57,16){\vector(1,0){.04}}\put(39.96,15.96){\line(1,0){.944}}
\put(41.84,15.96){\line(1,0){.944}}
\put(43.73,15.96){\line(1,0){.944}}
\put(45.62,15.96){\line(1,0){.944}}
\put(47.51,15.96){\line(1,0){.944}}
\put(49.4,15.96){\line(1,0){.944}}
\put(51.29,15.96){\line(1,0){.944}}
\put(53.18,15.96){\line(1,0){.944}}
\put(55.07,15.96){\line(1,0){.944}}
%\end
\thicklines
\put(57,16){\line(0,1){5}}
\thinlines
%\vector{dash}(57,19)(60,19)
\put(60,19){\vector(1,0){.04}}\put(56.96,18.96){\line(1,0){.75}}
\put(58.46,18.96){\line(1,0){.75}}
%\end
\thicklines
\put(60,19){\line(0,1){3}}
\thinlines
%\vector{dash}(24,10)(62,10)
\put(62,10){\vector(1,0){.04}}\put(23.96,9.96){\line(1,0){.974}}
\put(25.9,9.96){\line(1,0){.974}}
\put(27.85,9.96){\line(1,0){.974}}
\put(29.8,9.96){\line(1,0){.974}}
\put(31.75,9.96){\line(1,0){.974}}
\put(33.7,9.96){\line(1,0){.974}}
\put(35.65,9.96){\line(1,0){.974}}
\put(37.6,9.96){\line(1,0){.974}}
\put(39.55,9.96){\line(1,0){.974}}
\put(41.49,9.96){\line(1,0){.974}}
\put(43.44,9.96){\line(1,0){.974}}
\put(45.39,9.96){\line(1,0){.974}}
\put(47.34,9.96){\line(1,0){.974}}
\put(49.29,9.96){\line(1,0){.974}}
\put(51.24,9.96){\line(1,0){.974}}
\put(53.19,9.96){\line(1,0){.974}}
\put(55.14,9.96){\line(1,0){.974}}
\put(57.08,9.96){\line(1,0){.974}}
\put(59.03,9.96){\line(1,0){.974}}
\put(60.98,9.96){\line(1,0){.974}}
%\end
\thicklines
\put(62,10){\line(0,1){3}}
\thinlines
%\vector{dash}(7,8)(24,8)
\put(24,8){\vector(1,0){.04}}\put(6.96,7.96){\line(1,0){.944}}
\put(8.84,7.96){\line(1,0){.944}}
\put(10.73,7.96){\line(1,0){.944}}
\put(12.62,7.96){\line(1,0){.944}}
\put(14.51,7.96){\line(1,0){.944}}
\put(16.4,7.96){\line(1,0){.944}}
\put(18.29,7.96){\line(1,0){.944}}
\put(20.18,7.96){\line(1,0){.944}}
\put(22.07,7.96){\line(1,0){.944}}
%\end
%\vector{dash}(10,14)(20,14)
\put(20,14){\vector(1,0){.04}}\put(9.96,13.96){\line(1,0){.909}}
\put(11.77,13.96){\line(1,0){.909}}
\put(13.59,13.96){\line(1,0){.909}}
\put(15.41,13.96){\line(1,0){.909}}
\put(17.23,13.96){\line(1,0){.909}}
\put(19.05,13.96){\line(1,0){.909}}
%\end
\thicklines
\put(20,14){\line(0,1){3}}
\thinlines
%\vector{dash}(20,16)(22,16)
\put(22,16){\vector(1,0){.04}}\put(19.96,15.96){\line(1,0){.667}}
\put(21.29,15.96){\line(1,0){.667}}
%\end
\thicklines
\put(22,16){\line(0,1){2}}
\thinlines
\put(4,4){\vector(0,1){30}}
%\dashline(3,4)(67.13,4)
\put(2.96,3.96){\line(1,0){.987}}
\put(4.93,3.96){\line(1,0){.987}}
\put(6.9,3.96){\line(1,0){.987}}
\put(8.88,3.96){\line(1,0){.987}}
\put(10.85,3.96){\line(1,0){.987}}
\put(12.82,3.96){\line(1,0){.987}}
\put(14.8,3.96){\line(1,0){.987}}
\put(16.77,3.96){\line(1,0){.987}}
\put(18.74,3.96){\line(1,0){.987}}
\put(20.72,3.96){\line(1,0){.987}}
\put(22.69,3.96){\line(1,0){.987}}
\put(24.66,3.96){\line(1,0){.987}}
\put(26.63,3.96){\line(1,0){.987}}
\put(28.61,3.96){\line(1,0){.987}}
\put(30.58,3.96){\line(1,0){.987}}
\put(32.55,3.96){\line(1,0){.987}}
\put(34.53,3.96){\line(1,0){.987}}
\put(36.5,3.96){\line(1,0){.987}}
\put(38.47,3.96){\line(1,0){.987}}
\put(40.45,3.96){\line(1,0){.987}}
\put(42.42,3.96){\line(1,0){.987}}
\put(44.39,3.96){\line(1,0){.987}}
\put(46.37,3.96){\line(1,0){.987}}
\put(48.34,3.96){\line(1,0){.987}}
\put(50.31,3.96){\line(1,0){.987}}
\put(52.29,3.96){\line(1,0){.987}}
\put(54.26,3.96){\line(1,0){.987}}
\put(56.23,3.96){\line(1,0){.987}}
\put(58.21,3.96){\line(1,0){.987}}
\put(60.18,3.96){\line(1,0){.987}}
\put(62.15,3.96){\line(1,0){.987}}
\put(64.13,3.96){\line(1,0){.987}}
\put(66.1,3.96){\line(1,0){.987}}
%\end
\thicklines
\put(53,25){\line(0,1){2}}
\thinlines
%\vector{dash}(48,27)(50,27)
\put(50,27){\vector(1,0){.04}}\put(47.96,26.96){\line(1,0){.667}}
\put(49.29,26.96){\line(1,0){.667}}
%\end
\thicklines
\put(50,27){\line(0,1){3}}
\thinlines
%\vector{dash}(50,28)(52,28)
\put(52,28){\vector(1,0){.04}}\put(49.96,27.96){\line(1,0){.667}}
\put(51.29,27.96){\line(1,0){.667}}
%\end
\thicklines
\put(52,28){\line(0,1){3}}
\put(28.5,11.62){\makebox(0,0)[lc]{$\alpha(u)$}}
\put(25.75,27.25){\makebox(0,0)[rc]{$\omega(u)$}}
\put(28.62,19.37){\makebox(0,0)[lc]{$\alpha(v)$}}
\put(25.37,17.37){\makebox(0,0)[cc]{$u$}}
\put(35.87,22.87){\makebox(0,0)[cc]{$v$}}
\put(5.5,9.75){\makebox(0,0)[cc]{$\emptyset$}}
\put(7.62,2.5){\makebox(0,0)[cc]{$\rho$}}
\end{picture}

\caption{ A representation of a \emph{chronological tree} $\TT$ with finite \emph{discrete part} $\cal T$. Horizontal axis has no interpretation, but horizontal arrows indicate filiation; vertical axis indicates chronological levels. Three elements of $\cal T$ are shown, its root $\emptyset$, a typical individual $v$ and her mother $u$. Various elements of $\TT$ are shown, its root $\rho$, the death level $\omega(u)$ of $u$, as well as the birth levels $\alpha(u)$ and $\alpha(v)$.}
\label{fig : arbre-CMJ}
\end{figure}
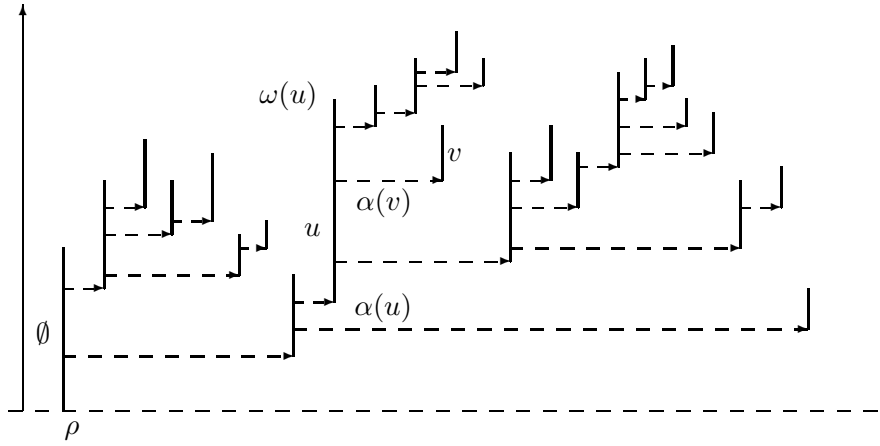

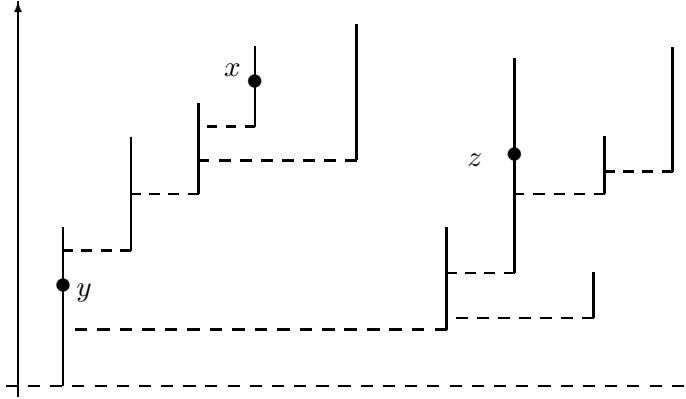
\begin{figure}[ht]

\unitlength 1.5mm % = 4.55pt
\linethickness{0.6pt}
\ifx\plotpoint\undefined\newsavebox{\plotpoint}\fi % GNUPLOT compatibility
\begin{picture}(65,41)(-10,0)
\put(6,1){\vector(0,1){35}}
%\dashline{1}(5,2)(65,2)
\put(4.96,1.96){\line(1,0){.984}}
\put(6.92,1.96){\line(1,0){.984}}
\put(8.89,1.96){\line(1,0){.984}}
\put(10.86,1.96){\line(1,0){.984}}
\put(12.82,1.96){\line(1,0){.984}}
\put(14.79,1.96){\line(1,0){.984}}
\put(16.76,1.96){\line(1,0){.984}}
\put(18.73,1.96){\line(1,0){.984}}
\put(20.69,1.96){\line(1,0){.984}}
\put(22.66,1.96){\line(1,0){.984}}
\put(24.63,1.96){\line(1,0){.984}}
\put(26.6,1.96){\line(1,0){.984}}
\put(28.56,1.96){\line(1,0){.984}}
\put(30.53,1.96){\line(1,0){.984}}
\put(32.5,1.96){\line(1,0){.984}}
\put(34.46,1.96){\line(1,0){.984}}
\put(36.43,1.96){\line(1,0){.984}}
\put(38.4,1.96){\line(1,0){.984}}
\put(40.37,1.96){\line(1,0){.984}}
\put(42.33,1.96){\line(1,0){.984}}
\put(44.3,1.96){\line(1,0){.984}}
\put(46.27,1.96){\line(1,0){.984}}
\put(48.23,1.96){\line(1,0){.984}}
\put(50.2,1.96){\line(1,0){.984}}
\put(52.17,1.96){\line(1,0){.984}}
\put(54.14,1.96){\line(1,0){.984}}
\put(56.1,1.96){\line(1,0){.984}}
\put(58.07,1.96){\line(1,0){.984}}
\put(60.04,1.96){\line(1,0){.984}}
\put(62.01,1.96){\line(1,0){.984}}
\put(63.97,1.96){\line(1,0){.984}}
%\end
\put(10,16){\line(0,-1){14}}
\put(16,24){\line(0,-1){10}}
\put(22,19){\line(0,1){8}}
\put(27,25){\line(0,1){7}}
\put(36,22){\line(0,1){12}}
\put(44,7){\line(0,1){9}}
\put(50,12){\line(0,1){19}}
%\dashline{1}(16,14)(10,14)
\put(15.96,13.96){\line(-1,0){.857}}
\put(14.24,13.96){\line(-1,0){.857}}
\put(12.53,13.96){\line(-1,0){.857}}
\put(10.81,13.96){\line(-1,0){.857}}
%\end
%\dashline{1}(22,19)(16,19)
\put(21.96,18.96){\line(-1,0){.857}}
\put(20.24,18.96){\line(-1,0){.857}}
\put(18.53,18.96){\line(-1,0){.857}}
\put(16.81,18.96){\line(-1,0){.857}}
%\end
%\dashline{1}(27,25)(22,25)
\put(26.96,24.96){\line(-1,0){.833}}
\put(25.29,24.96){\line(-1,0){.833}}
\put(23.62,24.96){\line(-1,0){.833}}
%\end
%\dashline{1}(36,22)(22,22)
\put(35.96,21.96){\line(-1,0){.933}}
\put(34.09,21.96){\line(-1,0){.933}}
\put(32.22,21.96){\line(-1,0){.933}}
\put(30.36,21.96){\line(-1,0){.933}}
\put(28.49,21.96){\line(-1,0){.933}}
\put(26.62,21.96){\line(-1,0){.933}}
\put(24.76,21.96){\line(-1,0){.933}}
\put(22.89,21.96){\line(-1,0){.933}}
%\end
%\dashline{1}(50,12)(44,12)
\put(49.96,11.96){\line(-1,0){.857}}
\put(48.24,11.96){\line(-1,0){.857}}
\put(46.53,11.96){\line(-1,0){.857}}
\put(44.81,11.96){\line(-1,0){.857}}
%\end
%\dashline{1}(44,7)(10.13,7)
\put(43.96,6.96){\line(-1,0){.996}}
\put(41.96,6.96){\line(-1,0){.996}}
\put(39.97,6.96){\line(-1,0){.996}}
\put(37.98,6.96){\line(-1,0){.996}}
\put(35.99,6.96){\line(-1,0){.996}}
\put(33.99,6.96){\line(-1,0){.996}}
\put(32,6.96){\line(-1,0){.996}}
\put(30.01,6.96){\line(-1,0){.996}}
\put(28.02,6.96){\line(-1,0){.996}}
\put(26.02,6.96){\line(-1,0){.996}}
\put(24.03,6.96){\line(-1,0){.996}}
\put(22.04,6.96){\line(-1,0){.996}}
\put(20.05,6.96){\line(-1,0){.996}}
\put(18.06,6.96){\line(-1,0){.996}}
\put(16.06,6.96){\line(-1,0){.996}}
\put(14.07,6.96){\line(-1,0){.996}}
\put(12.08,6.96){\line(-1,0){.996}}
%\end
\put(58,19){\line(0,1){5}}
\put(64,21){\line(0,1){11}}
\put(57,8){\line(0,1){4}}
%\dashline{1}(57,8)(44,8)
\put(56.96,7.96){\line(-1,0){.929}}
\put(55.1,7.96){\line(-1,0){.929}}
\put(53.24,7.96){\line(-1,0){.929}}
\put(51.38,7.96){\line(-1,0){.929}}
\put(49.53,7.96){\line(-1,0){.929}}
\put(47.67,7.96){\line(-1,0){.929}}
\put(45.81,7.96){\line(-1,0){.929}}
%\end
%\dashline{1}(58,19)(50,19)
\put(57.96,18.96){\line(-1,0){.889}}
\put(56.18,18.96){\line(-1,0){.889}}
\put(54.4,18.96){\line(-1,0){.889}}
\put(52.62,18.96){\line(-1,0){.889}}
\put(50.84,18.96){\line(-1,0){.889}}
%\end
%\dashline{1}(64,21)(58,21)
\put(63.96,20.96){\line(-1,0){.857}}
\put(62.24,20.96){\line(-1,0){.857}}
\put(60.53,20.96){\line(-1,0){.857}}
\put(58.81,20.96){\line(-1,0){.857}}
%\end
\put(25,30){\makebox(0,0)[cc]{$x$}}
\put(27,29){\circle*{1.12}}
\put(10,10.88){\circle*{1.12}}
\put(50,22.5){\circle*{1.12}}
\put(11.88,10.25){\makebox(0,0)[cc]{$y$}}
\put(46.5,22){\makebox(0,0)[cc]{$z$}}
\end{picture}

\caption{Three points $x, y, z$ in a finite chronological tree, satisfying  $y\prec x$ and $x\le y \le z$.
%The \emph{heights} (generations in the discrete tree) of points $x, y, z$ are respectively 3, 0, 2.
}
\label{fig : smalltree}
\end{figure}

Specifically, each individual of the underlying discrete tree possesses a \emph{birth level} $\alpha$ and a \emph{death level} $\omega$, both nonnegative real numbers such that $\alpha<\omega$, and (possibly zero) offspring whose birth times are distinct from one another and belong to the interval $(\alpha, \omega)$. We think of a \emph{chronological tree} as the set of all so-called \emph{existence points} of individuals (vertices) of the discrete tree. See Fig. \ref{fig : arbre-CMJ} and \ref{fig : smalltree} for graphical representations of a chronological tree.

\paragraph{Definition.}
More rigorously, let
$$
{\UU} = {\cal U}\times [0,+\infty),
$$
and set $\rho:=(\emptyset,0)$.\\
\\
We let $p_1$ and $p_2$ stand respectively for the canonical projections on $\cal U$ and $[0,+\infty)$.\\
\\
The first projection of any subset $\TT$ of $\UU$ will be denoted by ${\cal T}$
$$
{\cal T}:=p_1(\TT)=\{u\in \TT : \exists \sigma\ge 0, (u,\sigma)\in\TT\}.
$$
A \emph{chronological tree} $\TT$ is a subset of $\UU$ such that\\
\indent
(i) $\rho\in\TT$ (the \emph{root})\\
\indent
(ii) ${\cal T}$ is a discrete tree (as defined in the previous subsection)\\
\indent
(iii)
for any $u\in{\cal T}$, there are $0\le \alpha(u) <\omega (u)\le\infty$ such that $(u,\sigma)\in\TT$ if and only if $\sigma \in (\alpha(u), \omega(u)]$\\
\indent
(iv)
for any $u\in{\cal T}$ and $j\in\NN$ such that $uj\in\cal T$, $\alpha(uj) \in (\alpha(u), \omega(u))$.\\
\indent
(v)
for any $u\in\cal T$ and $i,j\in\NN$ such that $ui, uj\in\cal T$, 
$$
i\not=j \Rightarrow \alpha(ui) \not= \alpha (uj).
$$
For any $u\in{\cal T}$, $\alpha (u)$ is the birth level of $u$, $\omega(u)$ its death level, and we denote by $\zeta(u)$ its \emph{lifespan} $\zeta(u):=\omega(u)-\alpha(u)$.

Observe that (iii) implies that if $\TT$ is not reduced to $\rho$ then $\emptyset$ has a positive lifespan $\zeta(\emptyset)$. We will always assume that $\alpha(\emptyset)=0$.

%Recall that for any $u\in\cal T$, $K_u$ is the number of offspring of $u$. When $K_u<\infty$, it will be useful to use the following convention for any integer $j$
%$$
%j>K_u \Rightarrow\alpha(uj):= \omega(u).
%$$
The (possibly infinite) number of individuals alive at chronological level $\tau$ is denoted by $\Xi_\tau$
$$
\Xi_\tau = \mbox{Card}\{ v\in {\cal T} : \alpha(v) <\tau \le \omega(v)\}=\mbox{Card}\{ x\in {\TT} : p_2(x)=\tau\}\le\infty,
$$
and $(\Xi_\tau; \tau\ge 0)$ is usually called the \emph{width process}.

%Also for any $u\in{\cal T}$ and $j\in\NN$ such that $uj\in\cal T$, the birth time $\alpha(uj)$ of $u$'s daughter $uj$ is in $(\alpha(u), \omega(u))$. Points of the type $(u,\alpha(uj))$ are called \emph{branching points}. 

%\paragraph{Labellings.}
%We say that $\TT$ is numbered in the \emph{primogeniture labelling} if for any $u\in\cal T$ and $i,j\le K_u$,
%$$
%i<j\Leftrightarrow \alpha(ui) < \alpha (uj).
%$$
%If for any $u\in\cal T$ and $i,j\le K_u$,
%$$
%i<j\Leftrightarrow \alpha(ui) > \alpha (uj),
%$$
%we say that $\TT$ is numbered in the \emph{Gospel labelling}  (because the last daughter is ranked first).

\paragraph{Genealogical and metric structures.}
A chronological tree can naturally be equipped with the following genealogical structure and metric. 
For any $x,y\in\TT$ such that $x=(u,\sigma)$ and $y=(v,\tau)$, we will say that $x$ is an \emph{ancestor} of $y$, and write $x\prec y$ as for discrete trees, if $u\prec v$ and: 
\begin{itemize}
\item
if $u=v$, then $\sigma\le \tau$
\item 
if $u\not= v$, then $\sigma \le \alpha(uj)$, where $j$ is the unique integer such that $uj\prec v$. 
\end{itemize}
For $y=(v,\tau)$, the \emph{segment} $[\rho, y]$ is the set of ancestors of $y$, that is
\debeq
[\rho, y]	&:=& \{x\in\TT : x\prec y\}\\
					&=&\{(v,\sigma) : \alpha(v)<\sigma\le \tau \}\cup\{ (u,\sigma) : \exists k, u = v\vert k , \alpha(v\vert k)<\sigma\le \alpha(v\vert k+1)\}.
\fineq
For any $x,y\in\TT$, it is not difficult to see that there is a unique existence point $z\in\TT$ such that 
$[\rho, x]\cap [\rho, y]=[\rho, z]$. This point is the point of highest level in $\TT$ such that $z\prec x$ and $z\prec y$. In particular, notice that $p_1(z)=p_1(x)\wedge p_1(y)$ (i.e. $p_1(z)$ is the mrca of $p_1(x)$ and $p_1(y)$). The level $p_2(z)$ is called the \emph{coalescence level} of $x$ and $y$, and $z$ the \emph{coalescence point} (or most recent common ancestor) of $x$ and $y$, denoted as for discrete trees by $z=x\wedge y$. The \emph{segment} $[x,y]$ is then defined as
$$
[x,y]:=[\rho, x]\cup[\rho,y]\backslash [\rho,x\wedge y[,
$$
where a reversed bracket means that the corresponding extremity is excluded.
A natural distance $d$ on $\TT$ can readily be defined as
$$
d(x,y):= p_2(x)+p_2(y)-2p_2(x\wedge y).
$$
Note that $p_2$ is also the distance to the root.
%
%
%
%
%
%
%
%
%Observe that the tree structure and metric of a chronological tree do not depend on the numbering of daughters in $\cal T$.

\paragraph{Degree.}
The \textit{degree} of a point $x\in\TT$, i.e. the number of connected components of $\TT\backslash\{x\}$, can be equal to 1, 2 or 3. Apart from the root $\rho$, points of degree 1 are called \emph{death points} or \emph{leaves} and are those $x=(u,\sigma)$ such that $\sigma = \omega(u)$. Points of degree 2 are called \emph{simple points}. Points of degree 3 are called \emph{birth points} or \emph{branching points} and are those $x=(u,\sigma)$ such that $\sigma=\alpha(uj)$ for some integer $j\le K_u$. For example, if $x$ is not an ancestor of $y$ and $y$ is not an ancestor of $x$, then $x\wedge y$ is a branching point.

\paragraph{Grafting.} \label{11}
Let $\TT, \TT'$ be two chronological trees and $x=(u,\sigma)\in \TT$ a point of degree 2. For any positive integer $i$, we denote by $g(\TT', \TT, x, i)$ the tree obtained by \emph{grafting} $\TT'$ on $\TT$ at $x$, as descending from $ui$. More precisely, denote by $\tilde{\TT}$ the tree  obtained from $\TT$ by renaming points $(ukw,\tau)$ as $(u(k+1)w,\tau)$, for all $k\ge i$ and finite integer words $w$ such that $(ukw,\tau)\in\TT$. Then $g(\TT', \TT, x,i)$ is given by
$$
g(\TT', \TT, x, i):= \tilde{\TT}\cup\{ (uiw, \sigma+\tau) : (w,\tau)\in\TT'\}.  
$$
Observe that $g(\TT', \TT, x,i)$ indeed is a chronological tree. %, and that if $\TT$ was in the Gospel or primogeniture labelling, this is in general no longer the case after grafting. 

\paragraph{Planar embedding.}
The trees  we consider can be regarded as \emph{plane} trees satisfying the rule: `edges always grow to the right'.
For any $x\in\TT$, we denote by $\theta(x)$ the \emph{descendance} of $x$, that is, the subset of $\TT$ containing all $z\in\TT$ such that $x\prec z$.   The descendance of $x$ can be split out into its $l$(eft)-descendance $\theta_l (x)$ and $r$(ight)-descendance $\theta_r (x)$. Their definitions are as follows: if $x$ is not a branching point, $\theta_l (x)=\theta (x)$ and $\theta_r (x)=\emptyset$; if $x=(u,\sigma)$ is a branching point, then $\sigma=\alpha(uj)$ for some integer $j\le K_u$ and
$$
\theta_l (x) := \bigcup_{\vareps>0}\ \theta(u,\sigma+\vareps)
\quad\mbox{ and }\quad
\theta_r(x) := \{x\}\cup\bigcup_{\vareps>0}\ \theta(uj,\sigma+\vareps) .
$$
Actually, since by definition $x$ belongs to its descendance $\theta(x)$, it has to belong to either $\theta_l(x)$ or $\theta_r(x)$. In the case when $x$ is a branching point, the most convenient convention (which we adopt) is that $x\in\theta_r(x)$.

 Then for any $x \in\TT$, the complement in $\TT$ of $[\rho,x]\cup\theta(x)$ can be partitioned into two forests that we call its \emph{left-hand component} $L(x)$ and its \emph{right-hand component} $R(x)$ as follows.
$$
L(x):=\bigcup_{y:\ x\in\theta_r(y)} \theta_l(y)
\quad\mbox{ and }\quad
R(x):= \bigcup_{y:\ x\in\theta_l(y)} \theta_r(y) .
$$
For $x\in\TT$, the branching points that are ancestor points of $x$ belong to one of the following sets
$$
G(x):=\{y\prec x : x\in\theta_r(y) \}
\quad\mbox{ or }\quad
D(x) := \{y\prec x : x\not\in\theta_r(y)\not=\emptyset\} .
$$
Elements of $G(x)$ (resp.) are called \emph{left (resp. right) branching ancestor points} of $x$, in short \emph{lbap} (resp. \emph{rbap}).

\paragraph{Linear order `$\le$'.}
Let $x, y\in\TT$. Assume that $x\wedge y\not\in\{x,y\}$. Then either $y\in\theta_r(x\wedge y)$ (and  then $x\in\theta_l(x\wedge y)$) or  $x\in\theta_r(x\wedge y)$ (and  then $y\in\theta_l(x\wedge y)$). As a consequence, either $x\in L(y)$ and $y\in R(x)$, or $y\in L(x)$ and $x\in R(y)$, so that the relation `$\leq$' defined on $\TT$ as follows is a \emph{total order}, or \emph{linear order} on $\TT$ (whereas `$\prec$' only defines a partial order):
$$
x\leq y \Leftrightarrow \left[ y\prec x \mbox{ or } x\in L(y)\right] \Leftrightarrow \left[ y\prec x \mbox{ or } y\in R(x) \right].
$$
Note that the genealogical order and the linear order are opposite on a segment $[\rho,x]$. This will have important consequences when our trees are random.

It is also important to notice that if $\TT$ is not reduced to the root, then for any $x\in\TT$,
$$
(\emptyset,\omega(\emptyset))\leq x\leq (\emptyset,\alpha(\emptyset))=\rho.
$$

\paragraph{Lebesgue measure.}

The Borel $\sigma$-field of a chronological tree $\TT$ can be defined as the $\sigma$-field generated by segments. Defining the \emph{Lebesgue measure} of a segment as the distance between its extremities, Caratheodory's theorem ensures the existence of Lebesgue measure, say $\lbd$, on the Borel sets of $\TT$.

We will most of the time call \emph{length} the measure $\lbd(\DoubleS)$ of a Borel subset $\DoubleS$ of $\TT$. Note that the \emph{total length} of the tree $\lbd(\TT)$ is the sum of all lifespans
$$
\lbd(\TT) = \sum_{u\in{\cal T}} \zeta(u)\le \infty.
$$
We define the \emph{truncation} of $\DoubleS$ at level $\tau$ as the  subset of points of $\DoubleS$ whose existence level is lower than $\tau$
$$
C_\tau(\DoubleS) := \{x\in \DoubleS : p_2(x)\le \tau\} .
$$
We will say that $\DoubleS$
\begin{itemize}
\item
 has \emph{locally finite length} if for any finite level $\tau$, $C_\tau(\DoubleS)$ has finite length
\item
is \emph{finite} if it has  finite discrete part $p_1(\DoubleS)$
\item
 is \emph{locally finite} if for any finite level $\tau$, $C_\tau(\DoubleS)$ is finite.
\end{itemize}
Recall that $\Xi_\tau$ is the number of individuals alive at level $\tau$. 
\begin{prop}
\label{prop : finite length}
For any $\tau> 0$, 
\begin{equation}
\label{eqn : longueur troncation}
\int_{0}^\tau \Xi_\sigma \, d\sigma= \lbd(C_\tau(\TT)) .
\end{equation}
If $\TT$ has locally finite length, then  
\begin{itemize}
\item
$\Xi_\tau<\infty$ for Lebesgue-a.e. $\tau$
\item
 $\TT$ has finite length iff $\Xi_\tau=0$ for all sufficiently large $\tau$.
\end{itemize}
\end{prop}
\paragraph{Proof.}
First check that
$$
\lbd(C_\tau(\TT)) = \sum_{u\in{\cal T}} \left( (\omega(u)\wedge \tau)-(\alpha(u)\wedge \tau)\right) .
$$
Next write $\Xi_\sigma$ as
$$
\Xi_\sigma= \sum_{u\in{\cal T}} \indic{\alpha(u)<\sigma \le \omega(u)},
$$
and use Fubini's theorem to get
$$
\int_{0}^\tau \Xi_\sigma d\sigma = \sum_{u\in{\cal T}}\int_{0}^\tau d\sigma\, \indic{\alpha(u)<\sigma \le \omega(u)}=\sum_{u\in{\cal T}} \left( (\omega(u)\wedge \tau)-(\alpha(u)\wedge \tau)\right) ,
$$
which ends the proof of \eqref{eqn : longueur troncation}.

The first item of the list is a mere re-statement of \eqref{eqn : longueur troncation}. To see the direct sense of the equivalence stated in the second item, pick any individual $u$ in the discrete part $\cal T$ of $\TT$. If $n=\vert u\vert$, and $u_k= u\vert k$ for $k=0,1,\ldots,n$, then since $\alpha(u_k)<\omega(u_{k-1})$, we get by induction
$$
\omega(u) = \zeta(u_n) + \alpha(u_n) < \zeta(u_n)+\omega(u_{n-1})
= \zeta(u_n)+\zeta(u_{n-1}) + \alpha(u_{n-1})<\cdots
<\sum_{k=0}^n \zeta(u_k).
$$
But when $\TT$ has finite length $\ell$, $\sum_{k=0}^n \zeta(u_k)\le \sum_{v\in{\cal T}} \zeta(v)=\ell$, so that $\omega(u) <\ell$. Then all individuals are dead at time $\ell$, so that $\Xi_\tau=0$ for any $\tau \ge \ell$. 

The converse is elementary. Indeed, if there is $\tau_0$ such that $\Xi_\tau=0$ for all $\tau\ge \tau_0$, then $C_\tau(\TT)=\TT$ for all $\tau\ge \tau_0$. So if $\TT$ has locally finite length, it has finite length. \hfill$\Box$\\
\\
Finally, we indicate how to locally close splitting trees that have locally finite length.
\paragraph{Local closure.} Let $\partial {\cal T}$ denote the (possibly empty) \emph{local boundary} of the discrete tree $\cal T$
$$
\partial{\cal T}:=\{\mbox{infinite sequences }u : u\vert n \in{\cal T} \mbox{ for all } n \mbox{ and }(\alpha(u\vert n))_n \mbox{ is bounded}\}.
$$
We will sometimes write $\overline{\cal T}:={\cal T}\cup \partial{\cal T}$ for the \emph{local closure} of ${\cal T}$.

Now assume that $\TT$ has locally finite length. Then for any $u\in\partial {\cal T}$, one has $\lim_{n\tendinfty} \zeta(u\vert n)=0$, so we can define 
$$
\nu(u):=\lim_{n\uparrow \infty}\uparrow \alpha(u\vert n)= \lim_{n\uparrow \infty} \omega(u\vert n).
$$
This allows to define the \emph{local boundary} of $\TT$ as
$$
\partial\TT := \{(u,\nu(u)) : u\in\partial {\cal T}\}.
$$
We point out that this closure  \label{17} is local in the sense that all points in $\partial\TT$ are at finite distance from the root. Also note that the set $\partial\TT$ is not necessarily countable (examples of uncountable boundaries will be seen later when trees are random). Points in $\partial\TT$ can be thought of as points with zero lifespan and infinite height in the discrete genealogy. They have no descendance and degree 1, so to distinguish them from other leaves, we will call them \emph{limit leaves}. Note that our terminology is a little bit abusive, since the local boundary of $\TT$, taken in the usual sense (boundary associated to the metric $d$), should also comprise actual leaves, as well as the root. Also note that limit leaves are their own mother in the discrete genealogy.
 
All other properties of $\TT$ (genealogy, metric, order,...) trivially extend to its \emph{local closure} $\overline{\TT}:=\TT\cup \partial \TT$, and details need not be written down. In particular, Lebesgue measure extends to $\overline{\TT}$: its Borel $\sigma$-field is the Borel $\sigma$-field of $\TT$ completed with all sets $A\cup B$, where $A$ is a Borel set of $\TT$ and $B$ is any subset of $\partial \TT$, with $\lbd (A\cup B):=\lbd(A)$. Note that $\lbd(\partial \TT)=0$.

\subsection{Spectrally positive Lévy processes}
All results stated in this subsection are well-known and can be found in \cite{B}.

We denote by $(Y_t;t\ge0)$ a real-valued Lévy process (i.e., a process with independent and homogeneous increments, and a.s. càdlàg paths) with \emph{no negative jumps}, and by $P_x$ its distribution conditional on $Y_0=x$. 
Its Laplace exponent $\psi$ is defined by 
$$
E_0(\exp(-\lbd Y_t))= \exp(t\psi(\lbd)),
$$
and is specified by the Lévy--Khinchin formula
\begin{equation}
\label{eqn : exp Lap X}
\psi(\lbd) = \alpha\lbd + \beta \lbd^2
+\int_0^\infty (\expf{-\lbd r} -1 +\lbd r\indicbis{r<1})\Lambda(dr)\qquad\lbd\ge 0,
\end{equation}
where $\alpha\in\RR$, $\beta\ge 0$\ denotes the \emph{Gaussian coefficient}, and the \emph{Lévy measure} $\Lambda$\ is a $\sigma$-finite measure on $(0,\infty]$\ such that $\int_0^\infty (r^2\wedge 1)\Lambda (dr)<\infty$. We will sometimes assume that $q:=\Lambda(\{+\infty\})=-\psi(0)$ is positive. This amounts to killing the process at rate $q$.

The paths of $Y$ have finite variation a.s. if and only if $\beta=0$ and $\int_0^1 r \Lambda (dr)<\infty$. Otherwise the paths of $Y$ have infinite variation   a.s.

When $Y$ has increasing paths a.s., it is called a \emph{subordinator}. In that case, $\psi(\lbd)<0$ for any positive $\lbd$, and we will prefer
to define its Laplace exponent  as $-\psi$.
Since a subordinator has finite variation, its Laplace exponent  can be written as
$$
-\psi(\lbd) = \mbox{d}\lbd
+\int_0^\infty (1-e^{-\lbd r})\,\Lambda(dr)\qquad\lbd\ge 0,
$$
where $\mbox{d}\ge 0$ is called the \emph{drift coefficient}.

Next assume that $Y$ is not a subordinator. Then, since $\psi$ is convex and ultimately positive,
$$
\lim_{\lbd\tendinfty}\psi(\lbd) = +\infty.
$$
Denote by $\eta$ the largest root of $\psi$. If $\psi(0)<0$ (case when $\Lambda(\{+\infty\})\not=0$), then $\eta>0$ is the unique root of $\psi$. If $\psi(0)=0$, then either $\psi'(0^+)< 0$ and
$\eta>0$, so that $\psi$\ has exactly two roots (0 and $\eta$), or $\psi'(0^+)\ge 0$ and $\psi$\ has a unique root $\eta=0$. The inverse of $\psi$\ on $[\eta, \infty)$ is denoted by $\phi : [0, \infty)\rightarrow [\eta, \infty)$ and has in particular $\phi(0)=\eta$.

We write $T_A = \inf\{t\geq 0 : Y_t \in A\}$ for the first entrance time of $Y$ in a Borel set $A$ of $\RR$, and $T_y$ for $T_{\{y\}}$. It is known that
$$
E_0 (e^{-q T_{-x}}) = e^{-\phi(q)x}\quad \quad \: q\ge 0 , x \ge 0.
$$
In particular, $P_0 (T_{-x}<\infty)=e^{-\eta x}$.

Finally, there exists a unique continuous function $W : [0,+\infty)\rightarrow [0,+\infty)$, with Laplace transform
$$
\int_0^{\infty} e^{-\lambda x} W(x) dx = \frac{1}{\psi (\lambda)} \qquad \lambda > \eta,
$$
such that for any $ 0<x<a$,
\begin{equation}
\label{F1}
P_{x} (T_0 <T_{(a, +\infty)} ) = W(a-x)/W(a).
\end{equation}
The function $W$\ is strictly increasing and called the \emph{scale function}.

\subsection{Jirina processes}

We call \emph{Jirina process} a branching process in discrete time and continuous state-space. Specifically, a Jirina process is a  time-homogeneous Markov chain $(Z_n;n\ge 0)$ with values in $[0,+\infty]$ satisfying the branching property (w.r.t. initial condition). Writing $Z_n(x)$ for the value at generation $n$ of the Jirina process starting from $Z_0=x\in[0,+\infty)$, the branching property implies that for each integer $n$, $(Z_n(x)\,; x\ge 0)$ has i.i.d. nonnegative increments. In particular, $(Z_1(x)\,; x\ge 0)$
is a \emph{subordinator}, that we prefer to denote $S$. Let $F$ be the Laplace exponent of $S$, d its drift coefficient and $\Lambda$ its Lévy measure. 
%$$
%\EE(\exp(-\lbd S(x)))=\exp(-xF(\lbd))\qquad\lbd, x\ge 0.
%$$

By the (homogeneous) Markov property, there are i.i.d. subordinators $(S_n)_{n\ge 1}$ distributed as $S$, such that, conditional on $Z_0, Z_1,\ldots, Z_n$, 
$$
Z_{n+1} = S_{n+1}\circ Z_n.
$$
In particular, by Bochner's subordination,  the process $x\mapsto Z_n(x)$ is a subordinator with Laplace exponent $F_n$ the $n$-th iterate of $F$, so that
$$
\EE_x(\exp(-\lbd Z_n))=\exp(-xF_n(\lbd))\qquad\lbd,x\ge 0.
$$
We say that $Z$ is a Jirina process with \emph{branching mechanism} $F$. We will sometimes write $m:=F'(0^+)=\mbox{d}+\intgen r\Lambda(dr)\le \infty$. If $F(0)=0$, the Jirina process \label{18a}  is said subcritical, critical or supercritical according to whether $m<1$, $=1$ or $>1$. 
If $F(0)>0$ (i.e. $\Lambda(\{+\infty\})>0$), it is said supercritical.

We keep the same symbol $\eta$ as that used in the previous subsection, for the following quantity
$$
\eta:=\sup\{\lbd\ge 0 : F(\lbd)=\lbd\}<\infty,
$$
because we have in mind that the Lévy processes we will consider later on will have Laplace exponent $\lbd\mapsto \lbd -F(\lbd)$. Last, define the event of \emph{extinction}, denoted $\mathrm{Ext}$, as
$$
\mbox{Ext}:=\{\lim_{n\tendinfty}Z_n =0\} .
$$
\begin{prop}
\label{prop : properties Jirina}
Let $Z$ be a Jirina process with branching mechanism $F$. Then
$$
\PP_x(\mathrm{Ext}) = e^{-\eta x}\qquad x\ge 0. 
$$
Furthermore,
$$
\big\{\sum_{n\ge 0} Z_n <\infty\big\}=\mathrm{Ext} \qquad\qquad \mbox{a.s.}
$$
\end{prop}
Notice that if $\Lambda$ is infinite or d $>0$, then $Z_n >0$  a.s. at each generation $n$, even on extinction. On the other hand, if $S$ is a compound Poisson process, then the Borel--Cantelli lemma ensures that on the event of extinction, $Z_n=0$ for all sufficiently large $n$.

\paragraph{Proof.}
Recall that $F$ is concave increasing, so that $F(x)\ge x$ for $x\in[0,\eta]$ and  $F(x) \le x$ for $x\ge\eta$. This entails the convergence to $\eta$ of the sequence $(F_n(\lbd))_n$ for any initial value $\lbd> 0$. As a consequence, \label{21}
$$
\lim_{n\tendinfty}\EE_x(\exp(-\lbd Z_n))=e^{-\eta x}\qquad x,\lbd >0. 
$$
This last convergence implies that $\PP_x(\mathrm{Ext})\le \exp(-\eta x)$, because
$$
e^{-\eta x}\ge \lim_{n\tendinfty}\EE_x(\exp(-\lbd Z_n),\mathrm{Ext})= \PP_x(\mathrm{Ext})\qquad x,\lbd >0.
$$
%Taking extreme values for $\lambda$, this last convergence implies that, as $n\tendinfty$, $\PP_x(Z_n\in [0,\varepsilon))$ converges  to $\exp(-\eta x)$ and $\PP_x(Z_n\in (M, \infty))$ converges to $1-\exp(-\eta x)$. Then there is an event $A$ with probability $\exp(-\eta x)$ on which $Z_n$ converges in probability to 0; conversely, it converges to $+\infty$ on its complement, say $B$. The Borel--Cantelli lemma along with the Markov property then ensure that $\{\limsup Z_n>0\}\in B$ and $\{\liminf Z_n<\infty\}\in A$, which yields the pathwise convergence to 0 on $A$, and to $+\infty$ on $B$.\\
%\\
Now for any real number $\lbd\ge 0$, define by induction the sequence $(v_n(\lbd))_n$ as $v_0(\lbd) = \lbd$ and
$$
v_{n+1}(\lbd) = \lbd + F(v_n(\lbd))\qquad n\ge 0.
$$
Then an elementary induction argument shows that
$$
\EE_x\left(\exp\left(-\lbd \sum_{i=0}^n Z_i\right)\right) = \exp (-xv_n (\lbd))\qquad n \ge 0.
$$
On the one hand, by definition, $v_{n+1}=G_\lbd(v_n)$, where $G_\lbd:y\mapsto\lbd +F(y)$ is concave increasing and has $G_\lbd(0)\ge 0$. This implies that $G_\lbd$ has at least one fix point, and that its largest fix point, say $\phi(\lbd)$,
satisfying 
\begin{equation}
\label{eqn : inversion}
\phi(\lbd)-F\circ \phi(\lbd)=\lbd\qquad \lbd \ge 0,
\end{equation}
is also the limit of the sequence $(v_n(\lbd))$. On the other hand, Lebesgue convergence theorem gives the Laplace transform of
$$
T:=\sum_{i=0}^\infty Z_i \in(0,\infty]
$$
as
\begin{equation}
\label{eqn : TL somme Jirina}
\EE_x(\exp(-\lbd T))=\exp(-x\phi(\lbd)).
\end{equation}
As a consequence, $T<\infty$ with probability $\exp(-x\phi(0))$. But thanks to (\ref{eqn : inversion}), $\phi(0)$ is the largest fix point of $F$, so $\phi(0)=\eta$. Now since $\{T <\infty\}\subseteq \mbox{Ext}$, we get
$$
e^{-\eta x}\ge \PP_x(\mathrm{Ext}) \ge \PP_x(T <\infty)=e^{-\eta x},
$$
which implies that $\{T <\infty\}$ and $\mbox{Ext}$ have the same probability and hence coincide a.s.\hfill$\Box$\\
\\
Here the proof is over, but we want to point out that according to (\ref{eqn : inversion}), $\phi$ is the inverse of $\lbd\mapsto \lbd-F(\lbd)$, which itself is the Laplace exponent of the spectrally positive Lévy process $t\mapsto Y_t=S_t-t$. Referring to the last subsection on Lévy processes, equation \eqref{eqn : TL somme Jirina} thus shows that $T$ has the same law under $\PP_x$ as the first hitting time of 0 by the Lévy process $Y$ started at $x$. One of the goals of this paper is to shed some light on this relationship.

\section{The exploration process}
\label{sec : exploration process}

\subsection{Definition}

In most examples of $\RR$-trees, the latter are defined from their contour \cite{A1, A2, LG93, LGLJ}, which is a real function coding the genealogy, whereas here, we do the opposite (but see also \cite{Duq}).

Hereafter, $\TT$ denotes a chronological tree with \emph{finite total length} $\ell:=\lbd(\TT)$. The real interval $[0,\ell]$ is equipped with its Borel $\sigma$-field and Lebesgue measure, which we denote by `$\mbox{Leb}$'.\\
\\
For any $x\in \overline{\TT}$, set
$$
\DoubleS(x):=\{y\in\TT:y\le x\}.
$$
Since $\DoubleS(x)\setminus\{x\}$ is the union of segments of the form $]z,y]$ where $z=y\wedge x$ and $y$ ranges over the leaves $\le x$, and since the leaves of $\TT$ are in one-to-one correspondence with $\cal T$, which is at most countable, then $\DoubleS(x)$ is a Borel subset of $\TT$. It is then standard measure theory to prove that the real-valued mapping $\varphi :\overline{\TT}\mapsto [0,\ell]$ defined by
$$
\varphi(x):=\lbd(\DoubleS(x))\qquad x\in \overline{\TT}, 
$$
is measurable and pushes \label{26} forward Lebesgue measure on $\overline{\TT}$ to Lebesgue measure on $[0,\ell]$. For the sake of conciseness, we will say that $\varphi$ preserves the (Lebesgue) measure. By construction, it also preserves the order. Note that $\varphi(\emptyset,\omega(\emptyset))=0$ and $\varphi(\rho)=\ell$. 
 
Let us show that $\varphi$ is one-to-one and onto. First, because $\varphi$ preserves the order,  $\varphi(x)=\varphi(y)$ implies that $x\le y$ and $y\le x$, so that $x=y$, which proves that $\varphi$ is one-to-one. Second, let $t\in[0,\ell]$, and assume that $t\in D:= [0,\ell]\setminus\varphi(\TT)$. We want to show that then $t$ is the image by $\varphi$ of a limit leaf. As a first step, notice that because $\varphi$ preserves the measure, $D$ has zero Lebesgue measure. As a second step, we define for each $v\in\cal T$
$$
g(v):=\varphi(v,\omega(v))\quad\mbox{ and }\quad d(v):=\varphi(u,\alpha(v)),
$$
where $u$ denotes $v$'s mother. Then observe that for each $x\in\TT$ such that $v\prec p_1(x)$, one has
$$
g(v)\le \varphi(x) < d(v),
$$
so that in particular, either $p_1(x)=\emptyset$ or there is $v$ with $\vert v\vert = 1$ such that $\varphi(x)\in [g(v), d(v))$.
Assume that there is no $v$ with height 1 such that $t\in [g(v), d(v))$, and let 
$$
{\cal U}_1:= \{v \in{\cal T}: \vert v\vert =1,  g(v)<d(v)< t\}
\quad\mbox{ and }\quad
{\cal U}_2:= \{v \in{\cal T}: \vert v\vert =1, t<g(v)<d(v)\},
$$
so that ${\cal U}_1$ and ${\cal U}_2$ form a partition of the first generation in ${\cal T}$.
Then we set
$$
\sigma_1:=\inf_{v\in {\cal U}_1} \alpha(v)
\quad\mbox{ and }\quad
\sigma_2:=\sup_{v\in {\cal U}_2} \alpha(v).
$$
For each $i=1,2$, if $(\emptyset,\sigma_i)$ is a branching point, then there is an integer $k$ such that  $\sigma_i=\alpha(k)$ and we put $x_i:=(k,\omega(k))$, otherwise we put $x_i:=(\emptyset,\sigma_i)$.
Next observe that $\sup_{v\in{\cal U}_1} d(v)=\varphi(x_1)$ and $\inf_{v\in{\cal U}_2} g(v)=\varphi(x_2)$, so that
$$
\varphi(x_1)\le t\le \varphi(x_2).
$$
Actually, since by assumption $t$ is not in $\varphi(\TT)$, we get $\varphi(x_1)< t < \varphi(x_2)$, so in particular $x_1\not= x_2$ and $\sigma_1> \sigma_2$.
Now for any $x\in\TT$ such that $x_1< x < x_2$, and any pair $(v_1, v_2)\in {\cal U}_1\times {\cal U}_2$, we have
$$
g(v_1)<d(v_1)<\varphi(x)<g(v_2)<d(v_2),
$$
so that $p_1(x)=\emptyset$, and $x=(\emptyset,\sigma)$ for some $\sigma\in(\sigma_2, \sigma_1)$. Then it is easily seen that $\varphi(x)= \varphi(x_1) + \lbd([x,x_1])=\varphi(x_1) + \sigma_1 -\sigma$, and that this equality still holds for $x=x_2$. Therefore we get a contradiction, since we can always find $\sigma$ such that $t=\varphi(\emptyset, \sigma)$. In conclusion, there must be $v$ with height 1 such that $t\in (g(v), d(v))$ (by assumption, $t\not= g(v)$).

Repeating this argument to the subtree descending from $v$, and iterating, we get the existence of a unique infinite sequence $u\in\partial {\cal T}$ such that $t\in (g(u\vert n), d(u\vert n))$ for all $n$. These intervals are nested, and because $\TT$ has finite length, they have vanishing diameter ($d(u\vert n)-g(u\vert n)$ is the length of the subtree descending from $u\vert n$), so we get, as announced, $t=\varphi (u,\nu(u))$ for some $u\in \partial{\cal T}$.

As a conclusion, $\varphi$ is a bijection from $\overline{\TT}$ onto $[0,\ell]$ which preserves the order and the measure.\\
\\
Conversely, let $\psi$ be any order-preserving and measure-preserving bijection from $\overline{\TT}$  onto $[0, \ell]$. First, since  $\psi$ is order-preserving, $\psi^{-1}([0, \psi(x)])=\overline{\DoubleS(x)}$. Second, since $\psi$ is measure-preserving, $\lbd(\DoubleS(x))=\lbd(\overline{\DoubleS(x)})=\mathrm{Leb}([0, \psi(x)])=\psi(x)$. This can be recorded in the following statement.

\begin{thm}
\label{thm : dfn exploration}
The mapping $\varphi$ is the unique order-preserving and measure-preserving bijection from $\overline{\TT}$  onto the real interval $[0, \ell]$.  
\end{thm}

\begin{dfn}
\label{dfn : jccp}
The process $(\varphi^{-1}(t);t\in[0,\ell])$ is called the \emph{exploration process}.
Its second projection will be denoted by $(X_t; t\in [0, \ell])$ and called JCCP, standing for \emph{jumping chronological contour process}.
\end{dfn}

\begin{thm} 
\label{thm : allure jccp}
The exploration process is càdlàg (w.r.t. the distance $d$ on $\overline{\TT}$), and for any $t\in[0,\ell]$, $t$ is a jump time iff $\varphi^{-1}(t)$ is a leaf $(v, \omega(v))$ of $\TT$. In that case, $\varphi^{-1}(t-)=(u,\alpha(v))$, where $u$ is $v$'s mother in the discrete genealogy. 

As a consequence, the JCCP $(X_t; t\in [0, \ell])$ is also càdlàg and the size of each of its jumps is the lifespan of one individual. In addition,
\begin{equation}
\label{eqn : allure jccp}
X_t = -t + \sum_{ \varphi(v,\omega(v))\le t}\zeta(v)\qquad 0\le t\le\ell.
\end{equation}
%In particular, $X_0=\zeta(1)$ and $X_\ell = 0$.
\end{thm}%
The JCCP is a càdlàg function taking the values of all levels of all points in $\overline{\TT}$, once and once only, starting at the death level of the ancestor. In the \emph{finite} case, it follows this rule: 
when the visit of an individual $v$ with lifespan $(\alpha(v), \omega(v)]$ begins, the value of the JCCP is $\omega(v)$. The JCCP then visits lower chronological levels of $v$'s lifespan at constant speed $-1$. If $v$ has no child, then this visit lasts exactly the lifespan $\zeta(v)$ of $v$; if $v$ has at least one child, then the visit is interrupted each time a birth level of one of $v$'s daughters, say $w$, is encountered (youngest child first since the visit started at the death level). At this point, the JCCP jumps from $\alpha(w)$ to $\omega(w)$ and starts the visit of the existence levels of $w$. Since the tree has finite length, the visit of $v$ has to terminate: it does so at the chronological level $\alpha(v)$ and continues the exploration of the existence levels of $v$'s mother, at the level where it had been interrupted.
This procedure then goes on recursively until level $0$ is encountered ($0=\alpha(\emptyset)=$ birth level of the root).\\
\\
A chronological tree and its associated JCCP are represented on Fig. \ref{fig : JCCP}, which, because of its size, was moved to page \pageref{fig : JCCP}.

\begin{rem} 
In the case when the tree has finite discrete part, the JCCP has another interpretation \cite[Fig.1 p.230]{LGLJ} in terms of queues. Each jump $\Delta_t$ is interpreted as a customer of a one-server queue arrived at time $t$ with a load $\Delta_t$. This server treats the customers' loads at constant speed 1 and has priority LIFO (last in -- first out). \label{29} The tree structure is derived from the following rule: each customer is the mother of all customers who interrupted her while she was being served. A natural ranking of siblings (customers who interrupt the same service) is the order of their arrivals. The value $X_t$ of the JCCP  is the remaining load in the server at time $t$.
\end{rem}
Now we make two statements that the proof of Theorem \ref{thm : allure jccp} will require.

\paragraph{Claim 1.}
For any $x,y \in\TT$ such that $x\prec y$, we have
$$
d(x,y)\le \varphi(x) -\varphi(y).
$$

\paragraph{Claim 2.}
Let $x\in \TT$ and $(x_n)_{n\ge 0}$ a sequence of points of $\TT$ converging to $x$, such that one of the following conditions holds
\begin{itemize}
\item[(a)] $x_0\prec x_1 \prec x_2\prec\cdots \prec x$
\item[(b)] $x\prec\cdots\prec x_2\prec x_1\prec x_0$ and either $x$ is not a branching point, or it is one but $x_n\in\theta_r(x)$ for all $n$.
\end{itemize}
Then $\lim_{n\tendinfty} \varphi(x_n) =\varphi(x)$.

\paragraph{Proof of Claim 1.}
Recall that $d(x,y)=\lbd([x,y])$ and notice that $\DoubleS(x)\supseteq\DoubleS(y)\cup[x,y]$, so that, taking the measure of each side,
$$
\lbd (\DoubleS(x))\ge \lbd (\DoubleS(y))+\lbd([x,y])-\lbd(\DoubleS(y)\cap[x,y]).
$$ 
But $\DoubleS(y)\cap[x,y]=\{y\}$ and $\lbd (\{y\})=0$, so we get
$d(x,y)=\lbd([x,y])\le \lbd (\DoubleS(x))-\lbd (\DoubleS(y))=\varphi(x)-\varphi(y)$.
\hfill$\Box$

\paragraph{Proof of Claim 2.}
Since the proofs in both cases (a) and (b) are very similar, we only write it in case (a). Set $t:=\varphi(x)$ and $s_n:=\varphi(x_n)$. The genealogical ordering of the sequence induces the following ranking $x_0\ge x_1 \ge x_2\ge \cdots \ge x$, so that $(s_n)$ is a nonincreasing sequence whose limit $s:=\lim_{n\tendinfty} s_n$ has $s\ge t$. Let us suppose that $s>t$, and define $y:=\varphi^{-1}((s+t)/2)$. In particular, $x\le y\le x_n$ and $y\not=x, x_n$ for all $n$. Now let $z:=x\wedge y$. First, because $x\le y$ and $x\not=y$, we have $z\not=x$. Second, observe that for any $n$, $z$ and $x_n$ are on the segment $[\rho, x]$ so either $x_n\prec z$ or $z\prec x_n$; but since $x\le y\le x_n$ and $x_n\prec x$, we deduce that $z\le x_n$, so that $x_n\prec z$. In conclusion, $x_n\prec z\prec x$, so that $d(x_n,x)=d(x_n,z)+d(z,x)$. This brings about the contradiction, since $d(z,x)>0$, whereas $d(x_n,x)$ vanishes.
\hfill$\Box$

\paragraph{Proof of Theorem \ref{thm : allure jccp}.} \label{30}
Let $x\in\partial \TT$ and $t=\varphi(x)$. Recall from the beginning of this section that there is an infinite sequence $u\in \partial {\cal T}$ such that for any $n$, $t\in (g(u\vert n), d(u\vert n))$. Also notice that the intervals $(g(u\vert n), d(u\vert n))$ form a sequence of nested intervals decreasing to $t$. But for any $r,s$ in this interval, $\varphi^{-1}(r)$ and $\varphi^{-1} (s)$ are in the chronological subtree ${\cal T}_n$ descending from $u\vert n$, so that the distance between those two points in $\overline{\TT}$ is at most $\lbd ({\cal T}_n)=d(u\vert n)-g(u\vert n)$, which vanishes as $n\tendinfty$. This shows that $\varphi^{-1}$ is continuous at $t$, and allows us to discard $\partial\TT$ in the remainder of the proof. 

Let us show that $\varphi^{-1}$ is right-continuous. Let $t\in[0,\ell)$, and write $y=\varphi^{-1}(t)$. Next let $(t_n)_{n\ge 0}$ be a decreasing sequence converging to $t$, and write $y_n=\varphi^{-1}(t_n)$. For any $z\in]\rho, y[$, if $y_n\wedge y\not\in ]z,y[$, then $y_n$ is explored after $z$, that is, $\varphi(y_n)\ge \varphi(z) > \varphi(y)$. Since by assumption $\varphi(y_n)$ converges to $\varphi(y)$, we conclude that $y_n\wedge y\in ]z,y[$ for all sufficiently large $n$. This yields the convergence of $y_n\wedge y$ to $y$, and since $y_0\wedge y\prec y_1\wedge y \prec y_2\wedge y\prec\cdots \prec y$, we can apply Claim 2 (a) to the sequence $(y_n\wedge y)_n$, which gives 
$$
\lim_{n\tendinfty}\varphi(y_n\wedge y) = t.
$$
Next, since $y_n\wedge y\prec y_n$, we can apply Claim 1, which gives $d(y_n\wedge y, y_n)\le \varphi(y_n\wedge y)-\varphi(y_n)$. But both terms in the r.h.s. of the foregoing inequality converge to $t$, so that $d(y_n\wedge y, y_n)$ vanishes. We conclude with the triangular inequality $d(y_n,y)\le d(y_n,y_n\wedge y)+d(y_n\wedge y, y)$ which implies that $d(y_n,y)$ vanishes. In other words, for any decreasing sequence $(t_n)$ converging to $t$, $\varphi^{-1}(t_n)$ converges to $\varphi^{-1}(t)$, that is, $\varphi^{-1}$ is right-continuous.\\
\\
Next, we prove that $\varphi^{-1}$ has left-limits that can be characterized as in the theorem. Similarly as previously, let $t\in(0,\ell]$, and write $(v,\tau)=y=\varphi^{-1}(t)$. This time, let $(t_n)_{n\ge 0}$ be an increasing sequence converging to $t$, and write $y_n=\varphi^{-1}(t_n)$.

First, assume that $y$ is not a leaf. If $y$ is not a branching point either, then the proof that $(y_n)$ converges to $y$ can easily be adapted from that for the right-continuity, but appealing to Claim 2 (b) rather than (a). The same argument still applies if $y$ is a branching point, but to be allowed to appeal to Claim 2 (b), one first has to prove that $y_n\in\theta_r(y)$ for all sufficiently large $n$. Let us check that. For any $n$ and any $z\in\theta_r(y)$, if $y_n\not\in\theta_r(y)$ then $y_n\le z\le y$, so that $\varphi(y_n)<\varphi(z)<\varphi(y)$. But by assumption $(\varphi(y_n))$ converges to $\varphi(y)$, which proves that $y_n\in\theta_r(y)$ for all sufficiently large $n$.

Second, more interestingly, assume that $y$ is a leaf and $y\in\TT$, so that $\tau=\omega(v)$. Write $u$ for the mother of $v$, and set $z:=(u, \alpha(v))$, and $z':=(u,\omega(u))$. Since $y_n\le y$, $y_n\not\in\theta_r (z)$. Moreover, applying the same argument as in the end of the last paragraph, we find that for all sufficiently large $n$, $y_n\in\theta_l (z)$, or otherwise said, $y_n\wedge z'\in[z,z'[$. Actually, this same reasoning can be applied to any $z''\in[z,z'[$, namely $y_n\wedge z'\in[z,z''[$ for all sufficiently large $n$, which proves that $(y_n\wedge z')$ converges to $z$. Next, let $z_n$ be the midpoint of the segment $[z, y_n]$. Since $y_n\in\theta_l(z)$ for $n$ large enough, we have $z\prec z_n \prec y_n$, so that $\varphi(y_n)\le\varphi(z_n)$ and because $y\in R(z_n)$, $\varphi(z_n)\le t$. The immediate consequence is that
$$
\lim_{n\tendinfty} \varphi(z_n)=t.
$$
Now applying Claim 1 to $z_n\prec y_n$, we get $d(z_n,y_n)\le \varphi(z_n) - \varphi(y_n)$, and since both terms in the r.h.s. converge to $t$, we deduce that $d(z_n,y_n)$ vanishes. But by definition of $z_n$, $d(z,y_n)=2d(z_n,y_n)$, which shows that $(y_n)$ converges to $z$. In other words $\varphi^{-1}(t-)=(u,\alpha(v))$.\\
\\
The proof now focuses on the JCCP $X=(X_t;t\in[0,\ell])$. First, $X=p_2\circ \varphi^{-1}$ is càdlàg as a mere consequence of the fact that $\varphi^{-1}$ is càdlàg and $p_2$ is continuous (whereas $p_1$ is not, though). Now for any $t\in(0,\ell)$, $t$ is a jump time of $X$ only if \label{32} it is a jump time of $\varphi^{-1}$. When it is so, we know that $\varphi^{-1}(t)=(v,\omega(v))$ for some individual $v$, and that
$\varphi^{-1}(t-)=(u,\alpha(v))$, where $u$ is $v$'s mother. This shows that $X_t= \omega(v)$ and $X_{t-}=\alpha(v)$, so that the jump size at time $t$, $\Delta X_t=X_t-X_{t-}$, is equal to the lifespan $\zeta(v)$ of $v$. 

Now we prove \eqref{eqn : allure jccp}.
As said in the beginning of the proof, we can assume that $x=\varphi^{-1}(t)\in\TT$. We will need the following notation: for each point $y=(v,\tau)\in\TT$, we will set $\bar{y}:=(u,\alpha(v))$, where $u$ is $v$'s mother. Notice that since $\bar{y}\prec y$, we always have $y\le \bar{y}$. Then $\DoubleS(x)$ can be written as the union of $\CC(x)$ and $\UU(x)$, where
$$
\CC(x):=\{y\in\TT : y\le\bar{y}\le x\}
\quad\mbox{ and }\quad
\UU(x):=\{y\in\TT : y\le x\le\bar{y}\}.
$$
The only intersection between these subsets is $x$, and only if $x$ is a branching point.
We will use twice the following observation: if $(v,\tau)\in \DoubleS(x)$, then for any $\sigma\in[\tau,\omega(v)]$, since $(v,\tau)\prec (v,\sigma)$, we have $(v,\sigma)\le (v,\tau)\le x$, so that $(v,\sigma)\in\DoubleS(x)$.
The first consequence of this observation is that for any $y\in\CC(x)$, with $v=p_1(y)$, $(v,\sigma)\in \CC(x)$ for any $\sigma\in(\alpha(v),\omega(v)]$. We write ${\cal C}(x)$ the set of such vertices of the discrete tree. The second consequence is that for any $y\in\UU(x)$, with $y\not=x$ and $v=p_1(y)$, there is $\tau\in(\alpha(v),\omega(v))$ such that 
$$
(v,\sigma) \in\UU(x)\Leftrightarrow \sigma\in(\tau,\omega(v)],
$$
which implies that $(v,\sigma)\le x \le (v,\tau)$. Now thanks to Claim 1, unless $(v,\tau)=x$, $z=(v,\tau)$  is a branching point, $(v,\sigma)\in\theta_l(z)$ and $x\in \theta_r(z)$. In particular, if $w$ denotes $p_1(x)$, we have $v\prec w$ in the discrete tree, so that $p_1(\UU(x))$ is the set of ancestors of $w$, which has cardinality $n+1$, where $n=\vert w\vert$, and, thanks to the last display,
$$
\UU(x)=\big(w\times [p_2(x), \omega(w)]\big)\cup\bigcup_{k=0}^{n-1}\big((w\vert k)\times (\alpha(w\vert k+1), \omega(w\vert k)] \big).
$$
Taking the Lebesgue measure, we get
$$
\lbd(\UU(x))=\omega(w)-p_2(x) + \sum_{k=0}^{n-1} (\omega(w\vert k)-\alpha(w\vert k+1))
$$
On the other hand,
$$
\lbd(\CC(x))= \sum_{v\in {\cal C}(x)} \zeta(v), 
$$
so that, recalling that the intersection of $\CC(x)$ and $\UU(x)$ has zero Lebesgue measure, and that their union equals $\DoubleS(x)$, we deduce
$$
\lbd(\DoubleS(x))=\left(\sum_{v : (v,\omega(v))\le x}\zeta(v)\right) - \left(p_2(x)-\alpha(w) +\sum_{k=0}^{n-1} (\alpha(w\vert k+1)-\alpha(w\vert k))\right).
$$
Now notice that all terms in the second expression in the r.h.s. cancel out to $p_2(x)$, which yields
$$
p_2(x)=-\varphi(x)+\sum_{v : (v,\omega(v))\le x}\zeta(v).
$$
Writing $t=\varphi(x)$ yields \eqref{eqn : allure jccp}.
\hfill $\Box$

\subsection{Properties of the JCCP}

Actually, the chronological tree itself can be recovered from its JCCP  (modulo labelling of siblings). In the next two  statements, we provide some useful applications of this correspondence.

For each $t\in[0,\ell]$, set
$$
\hat{t} := \sup\{s\le t : X_s < X_{t}\}\vee 0\qquad 0\le t\le\ell. 
$$
\begin{thm}
\label{thm : proprietes jccp}
 Let $x=(u,\sigma)$ and $y=(v,\tau)$ denote any two points in $\overline{\TT}$, and set $s=\varphi(x)$ and $t=\varphi(y)$. Then the following hold :\\
\\
\indent
(i) 
The first visit to $v$ is $\hat{t}$ 
$$
\varphi (v, \omega(v))=\hat{t}.
$$
In particular, if $y\in\partial \TT$ then $\hat{t}=t$. 
If $t$ is a jump time of $X$, then $t = \hat{t}$ as well, and the first visit to the mother $u$ of $v$ in ${\cal T}$ is given by
$$
 \varphi (u, \omega(u)) = \sup\{s\le t : X_s < X_{t-}\}.
$$

\indent
(ii)
Ancestry between $x$ and $y$ :
$$
y\prec x\;\Leftrightarrow \;\hat{t}\le s\le t
$$

\indent
(iii) 
Coalescence level between $x$ and $y$ (assume e.g. $s\le t$) :
$$
p_2(x\wedge y) = \inf_{s\le r \le t} X_r.
$$
\end{thm}

For any $t\in[0,\ell]$, we define the process $X^{(t)}$ on $[0,t]$ as
$$
X^{(t)}_r := X_{t-}-X_{(t-r)-}\qquad r\in[0,t],
$$
with the convention that $X_{0-}=0$. We also set \label{36}
$$
H_t:=\vert p_1\circ\varphi^{-1}(t)\vert
$$
the generation, or genealogical height in $\overline{\cal T}$, of the individual $v=p_1\circ\varphi^{-1}(t)$ visited at time $t$ (recall that $\vert v\vert$ denotes the length of the integer word $v$).

The following corollary states \label{35} that if $v$ has finite height $H_t$, then the record times of $X^{(t)}$ are exactly those times when each of $v$'s ancestors is visited for the first time by the exploration process (which actually holds also when $v \in\partial{\cal T}$). It  also characterizes the height process $(H_t;t\ge 0)$ of genealogical heights, or \emph{height process}, as a functional of the path of the JCCP.
\begin{cor}
\label{cor : height process}
Let $y=\varphi^{-1}(t)$, $v=p_1(y)$, and $t_k$ the first visit to $v_k=v\vert k$ (ancestor of $v$  belonging to generation $k$), that is,
$$
t_k:=\varphi(v_k, \omega(v_k)).
$$
\indent
(i) if $n:=\vert v \vert<\infty$, then $y\not\in\partial\TT$, and one can define recursively the record times of $X^{(t)}$ by $s_1=t-\hat{t}$ and 
$$
s_{k+1} = \inf
\{s\ge s_k :  X_s^{(t)}> X_{s_k}^{(t)}\}\qquad\; k \ge 0.
$$
Then
$$
t_{k} = t- s_{n-k+1}\qquad\; 0\le k \le n.
$$
\indent
(ii) in the general case, recall $H_t\le \infty$ is the genealogical height of $v\in \overline{\cal T}$. Then $H_t$ is given by
\begin{eqnarray}
\label{eqn : height}
H_t & = &\mathrm{Card}\{0\le s\le t :  X_s^{(t)}=\sup_{0\le r\le s} X_r^{(t)}\}\nonumber\\
& = &\mathrm{Card}\{0\le s\le t : X_{s-}<\inf_{s\le r\le t} X_r^{}\}. 
\end{eqnarray}
\end{cor}

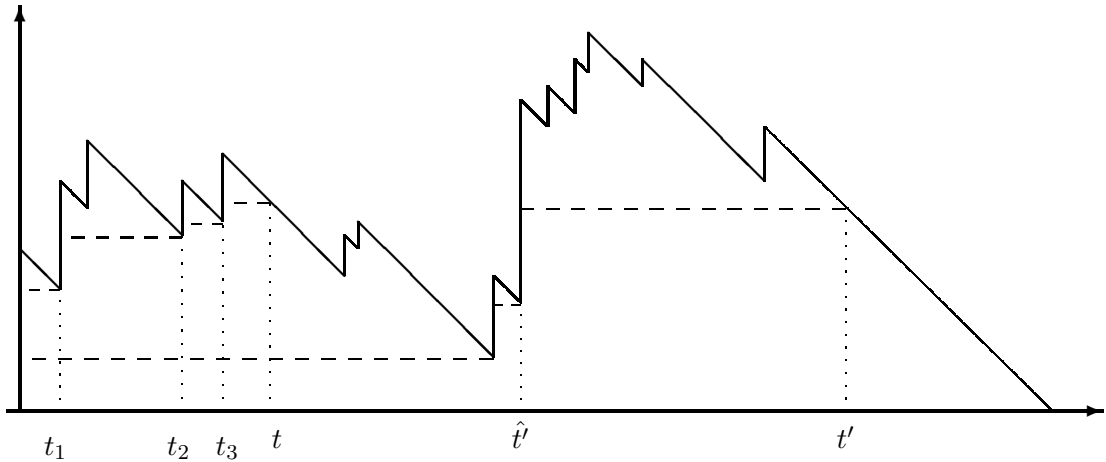
\begin{figure}[ht]
%\input{hauteur}

%TeXCAD Picture [hauteur.tex]. Options:
%\grade{\on}
%\emlines{\off}
%\epic{\off}
%\beziermacro{\on}
%\reduce{\on}
%\snapping{\off}
%\quality{8.00}
%\graddiff{0.01}
%\snapasp{1}
%\zoom{8.0000}
\unitlength 1.8mm % = 3.98pt
\linethickness{0.4pt}
\ifx\plotpoint\undefined\newsavebox{\plotpoint}\fi % GNUPLOT compatibility
\begin{picture}(82,35)(0,0)
%\line{dash}(32,49)(35,49)
\put(0,0){}
%\line{dash}(35,55)(38,55)
\put(0,0){}
%\line{dash}(35,53)(40,53)
\put(0,0){}
%\line{dash}(49,48)(52,48)
\put(0,0){}
%\line{dash}(52,61)(55,61)
\put(0,0){}
%\line{dash}(40,54)(43,54)
\put(0,0){}
%\line{dash}(55,62)(58,62)
\put(0,0){}
%\line{dash}(58,65)(61,65)
\put(0,0){}
%\line{dash}(58,64)(63,64)
\put(0,0){}
%\line{dash}(52,57)(60,57)
\put(0,0){}
%\line{dash}(52,51)(65,51)
\put(0,0){}
%\line{dash}(65,57)(68,57)
\put(0,0){}
%\line{dash}(65,55)(70,55)
\put(0,0){}
%\line{dash}(70,58)(73,58)
\put(0,0){}
%\line{dash}(73,61)(78,61)
\put(0,0){}
%\line{dash}(73,59)(80,59)
\put(0,0){}
%\line{dash}(65,52)(82,52)
\put(0,0){}
%\line{dash}(82,55)(85,55)
\put(0,0){}
%\line{dash}(49,46)(87,46)
\put(0,0){}
%\line{dash}(32,44)(49,44)
\put(0,0){}
%\line{dash}(35,50)(45,50)
\put(0,0){}
%\line{dash}(45,52)(47,52)
\put(0,0){}
%\line{dash}(73,63)(75,63)
\put(0,0){}
%\line{dash}(75,64)(77,64)
\put(0,0){}
%-
%-
%-
%-
%-
%-
%-
%-
%-
%-
%-
%-
%-
%-
%-
%-
%-
%-
%-
%-
%-
\thicklines
\put(4,5){\vector(0,1){30}}
\put(4,17){\line(1,-1){3}}
%\emline(7,22)(9,20)
\multiput(7,22)(.0240964,-.0240964){83}{\line(0,-1){.0240964}}
%\end
\put(9,25){\line(1,-1){7}}
\put(16,22){\line(1,-1){3}}
\put(19,24){\line(1,-1){9}}
%\emline(28,18)(29,17)
\multiput(28,18)(.0238095,-.0238095){42}{\line(0,-1){.0238095}}
%\end
\put(29,19){\line(1,-1){10}}
%\emline(39,15)(41,13)
\multiput(39,15)(.0240964,-.0240964){83}{\line(0,-1){.0240964}}
%\end
%\emline(41,28)(43,26)
\multiput(41,28)(.0240964,-.0240964){83}{\line(0,-1){.0240964}}
%\end
%\emline(43,29)(45,27)
\multiput(43,29)(.0240964,-.0240964){83}{\line(0,-1){.0240964}}
%\end
\put(45,27){\line(0,1){4}}
%\emline(45,31)(46,30)
\multiput(45,31)(.0238095,-.0238095){42}{\line(0,-1){.0238095}}
%\end
\put(46,30){\line(0,1){3}}
\put(46,33){\line(1,-1){4}}
\put(50,29){\line(0,1){2}}
\put(50,31){\line(1,-1){9}}
\put(59,22){\line(0,1){4}}
\put(3,5){\vector(1,0){81}}
\put(7,14){\line(0,1){8}}
\put(9,20){\line(0,1){5}}
\put(16,18){\line(0,1){4}}
\put(19,19){\line(0,1){5}}
\put(28,15){\line(0,1){3}}
\put(29,17){\line(0,1){2}}
\put(39,9){\line(0,1){6}}
\put(41,13){\line(0,1){15}}
\put(43,26){\line(0,1){3}}
\thinlines
%\dottedline(65,20)(65,5)
\multiput(64.95,19.95)(0,-.9375){17}{{\rule{.4pt}{.4pt}}}
%\end
%\dottedline(41,13)(41,5)
\multiput(40.95,12.95)(0,-.8889){10}{{\rule{.4pt}{.4pt}}}
%\end
%\dashline{1}(65,20)(41,20)
\put(64.95,19.95){\line(-1,0){.96}}
\put(63.03,19.95){\line(-1,0){.96}}
\put(61.11,19.95){\line(-1,0){.96}}
\put(59.19,19.95){\line(-1,0){.96}}
\put(57.27,19.95){\line(-1,0){.96}}
\put(55.35,19.95){\line(-1,0){.96}}
\put(53.43,19.95){\line(-1,0){.96}}
\put(51.51,19.95){\line(-1,0){.96}}
\put(49.59,19.95){\line(-1,0){.96}}
\put(47.67,19.95){\line(-1,0){.96}}
\put(45.75,19.95){\line(-1,0){.96}}
\put(43.83,19.95){\line(-1,0){.96}}
\put(41.91,19.95){\line(-1,0){.96}}
%\end
%\dashline{1}(22.5,20.38)(19,20.38)
\put(22.45,20.33){\line(-1,0){.875}}
\put(20.7,20.33){\line(-1,0){.875}}
%\end
%\dashline{1}(19,18.88)(16,18.88)
\put(18.95,18.83){\line(-1,0){.75}}
\put(17.45,18.83){\line(-1,0){.75}}
%\end
%\dashline{1}(16,17.88)(7,17.88)
\put(15.95,17.83){\line(-1,0){.9}}
\put(14.15,17.83){\line(-1,0){.9}}
\put(12.35,17.83){\line(-1,0){.9}}
\put(10.55,17.83){\line(-1,0){.9}}
\put(8.75,17.83){\line(-1,0){.9}}
%\end
%\dashline{1}(7,14)(4,14)
\put(6.95,13.95){\line(-1,0){.75}}
\put(5.45,13.95){\line(-1,0){.75}}
%\end
%\dottedline(22.5,20.25)(22.5,5)
\multiput(22.45,20.2)(0,-.9531){17}{{\rule{.4pt}{.4pt}}}
%\end
%\dottedline(19,18.63)(19,5)
\multiput(18.95,18.58)(0,-.9736){15}{{\rule{.4pt}{.4pt}}}
%\end
%\dashline{1}(41,12.88)(39,12.88)
\put(40.95,12.83){\line(-1,0){.667}}
\put(39.62,12.83){\line(-1,0){.667}}
%\end
%\dashline{1}(39,8.88)(4,8.88)
\put(38.95,8.83){\line(-1,0){.972}}
\put(37.01,8.83){\line(-1,0){.972}}
\put(35.06,8.83){\line(-1,0){.972}}
\put(33.12,8.83){\line(-1,0){.972}}
\put(31.17,8.83){\line(-1,0){.972}}
\put(29.23,8.83){\line(-1,0){.972}}
\put(27.28,8.83){\line(-1,0){.972}}
\put(25.34,8.83){\line(-1,0){.972}}
\put(23.39,8.83){\line(-1,0){.972}}
\put(21.45,8.83){\line(-1,0){.972}}
\put(19.51,8.83){\line(-1,0){.972}}
\put(17.56,8.83){\line(-1,0){.972}}
\put(15.62,8.83){\line(-1,0){.972}}
\put(13.67,8.83){\line(-1,0){.972}}
\put(11.73,8.83){\line(-1,0){.972}}
\put(9.78,8.83){\line(-1,0){.972}}
\put(7.84,8.83){\line(-1,0){.972}}
\put(5.89,8.83){\line(-1,0){.972}}
%\end
\put(18.5,2.25){\makebox(0,0)[lc]{$t_3$}}
\put(41,2.25){\makebox(0,0)[cb]{$\hat{t'}$}}
\put(65,2.25){\makebox(0,0)[cb]{$t'$}}
\put(23,2.25){\makebox(0,0)[cb]{$t$}}
%\emline(59,26)(80.25,5)
\multiput(59,26)(.0243692661,-.0240825688){872}{\line(1,0){.0243692661}}
%\end
%\dottedline(16,17.88)(16,5)
\multiput(15.95,17.83)(0,-.9908){14}{{\rule{.4pt}{.4pt}}}
%\end
%\dottedline(7,14)(7,5)
\multiput(6.95,13.95)(0,-.9){11}{{\rule{.4pt}{.4pt}}}
%\end
\put(15.75,2.25){\makebox(0,0)[cc]{$t_2$}}
\put(7.5,2.25){\makebox(0,0)[rc]{$t_1$}}
\end{picture}

\caption{The JCCP of some finite chronological tree with jumps in solid line. Set $u=p_1\circ \varphi^{-1}(t)$ (resp. $v=p_1\circ \varphi^{-1}(t')$) the individual visited at time $t$ (resp. $t'$). The first time when $v$ is visited by the exploration process is $\hat{t'}$. The first visits to the $H_t=3$ ancestors of $u$ are also shown. }
\label{fig : hauteur}
\end{figure}

Quantities defined in the  previous two statements are represented in Fig. \ref{fig : hauteur}.

\paragraph{Proof of Theorem \ref{thm : proprietes jccp}.}
(i) We first show that when $y\in\partial \TT$, then $\hat{t}=t$.
Indeed, recall that $v$ is an infinite sequence, and that $t=\lim_{n\uparrow\infty} \uparrow t_n$, where $t_n=\varphi(v\vert n,\omega(v\vert n))$. Also thanks to Theorem \ref{thm : allure jccp}, $X_{t_n-}=\alpha(v\vert n)$, so $X_{t_n-}$ increases to $X_{t-}=X_t$. This shows that $\hat{t}=t$.

As a by-product we also get the displayed equation in the case when $y\in\partial \TT$, since then there is only one pointwise visit $t$ of $v$, so the first visit is $t$ and we know that $t=\hat{t}$. From now on, we can discard the points of $\partial \TT$.

Set $z=(v,\omega(v))$, where $v=p_1\circ\varphi^{-1}(t)$ is the individual visited at time $t$ (and recall $y=\varphi^{-1}(t)$). Observe that $z\le y$ because $y\prec z$, so that $\varphi(z) \le t$,  and recall from Theorem \ref{thm : allure jccp} that $\varphi(z)$ is a jump time, with $X_{\varphi(z)-}=\alpha(v)<p_2(y)=X_t$. Now for any $s\in[\varphi(z), t]$, set $x=\varphi^{-1}(s)$ and $a=x\wedge z$. Since both $a$ and $y$ belong to $[\rho, z]$, either $a\prec y$ or $y\prec a$. We will suppose that $a\prec y$ and find a contradiction. Recall $z\le x$. Either $x\prec z$, so $a=x$ and $y\le x$, or $x\in R(z)$, so that $x\in R(y)$, and similarly $y\le x$. But unless $y=x$, this contradicts $s\le t$. As a consequence, $y\prec a$, so that $y\prec x$, and $X_t=p_2(y)\le p_2(x)=X_s$. Summing everything up, $X_{\varphi(z)-}<X_t$ and for any $s\in[\varphi(z), t]$, $X_s\ge X_t$. From these, we deduce that $\hat{t} = \varphi(z)$.

For the second statement, let $t$ be a jump time, and set $z=(u,\omega(u))$, where $u$ is the mother of $v$ in the discrete tree $\cal T$ (it exists since $\varphi^{-1}(t)$ is a leaf and so  is not in $\partial \TT$).  Then from Theorem \ref{thm : allure jccp} we know that $y=(v,\omega(v))$, that is, $t=\varphi(v,\omega(v))=\hat{t}$, that $X_{t-}=\alpha(v)$, and similarly $X_{\varphi(z)-}=\alpha(u)$. Because $y\in R(z)$, we get $\varphi(z)< t$, and because $v$ is a daughter of $u$, we get $\alpha(u)<\alpha(v)$, that is, $X_{\varphi(z)-}<X_{t-}$. Now for $s\in]\varphi(z), t[$, set $x=\varphi^{-1}(s)$ and $a=x\wedge z$. We want to show that $X_s>X_{t-}$, which also reads  $p_2(x)> \alpha(v)$.
Since $s>\varphi(z)$, either $x\prec z$ or $x\in R(z)$. If $x\prec z$, then $p_2(x)> \alpha(v)$ otherwise we would get $x\prec y$ and then $t\le s$. If $x\in R(z)$, then $a$ is a rbap of $z$ with $p_2(a)\ge \alpha(v)$ otherwise we would get $x\in R(y)$ and then $s\ge t$. But since $a\prec x$ and $a\not= x$, $p_2(x)> p_2(a)$, and we deduce again that $p_2(x)> \alpha(v)$. So we get $X_s>X_{t-}$ for all $s\in]\varphi(z), t[$, and the proof ends recalling that $X_{\varphi(z)-}<X_{t-}$.\\

(ii) 
Again set $z=(v,\omega(v))$, with $v=p_1(y)$ and $y=\varphi^{-1}(t)$. Recall that in the first paragraph above,
we have shown that for any $s\in[\varphi(z), t]$, $y\prec x(=\varphi^{-1}(s))$. Now thanks to (i), $\hat{t} = \varphi(z)$, so that for any $s\in[\hat{t}, t]$, $y\prec x$. Conversely, assume that $y\prec x$. First notice that $x\le y$, so that $s\le t$. Second, we want to show that $z\le x$, so that $\hat{t}\le s$. Set $a=x\wedge z$. Then either $x=a$, so that $x\prec z$ (and $z\le x$), or $a$ is a  branching ancestor point of $z$. Now since $y\prec x$ and $y\prec z$, we get $y\prec a$, and $a\in[y,z]$. Finally, because $p_1(y)=p_1(z)(=v)$, there is no lbap in the segment $[y,z]$, so that $a$ is a rbap of $z$ and $x\in R(z)$ (then $z\le x$ again). \\

(iii) 
Set $I_{s,t}=\inf_{s\le r \le t} X_r$. If $y\prec x$, then $p_2(x\wedge y)= p_2(y) =X_t$. Now, by (ii), $\hat{t}\le s \le t$, and by (i),
$I_{\hat{t},t}=X_t$, so that $ I_{s,t}=X_t$. Thus, we have proved (iii) when $x\wedge y =y$. Next set $z=x\wedge y$ and assume without loss of generality that $y\in R(x)$,  so $z$ is a branching point $(u,\alpha(v))=\varphi^{-1}(r)$. Since $z\prec y$, $ t\le r$. Also we have $z\prec x$, so $\hat{r}\le s$, which means that $X_h\ge X_r$ for any $h\in[s,r]\subset[s,t]$. In particular, $I_{s,t}\ge X_r=p_2(x\wedge y)$.
For the converse inequality, notice that $x\le a\le y$, where $a=(v,\omega(v))$ because $a\in R(x)$ and $y\in R(a)$. Also by definition, $p_1(x)\wedge p_1(y)=u$ and since $y\in R(x)$, $v$ is an ancestor individual of $p_1(y)$ but not of $p_1(x)$ (otherwise $p_1(x)\wedge p_1(y)=v$).  This ensures that $a\not=x$, and since $x\le a\le y$, $s< \varphi(a) \le t$. 
But $a$ is a leaf, so thanks to Theorem \ref{thm : allure jccp}, $\varphi^{-1}(\varphi(a)^-)=(u,\alpha(v))$, and $X_{\varphi(a)-}=\alpha(v)= p_2(x\wedge y)$. Because $s< \varphi(a) \le t$, we get $I_{s,t}\le  p_2(x\wedge y)$.\hfill$\Box$

\paragraph{Proof of Corollary \ref{cor : height process}.} 
We show (i) by descending induction on $k\in\{0,\ldots,n\}$. 
Since $\vert v\vert=n$, $v_n = v$, and  by the first statement in Theorem \ref{thm : proprietes jccp} (i), $t_n=\varphi(v, \omega (v)) = \hat{t}=t-s_1$, so the result holds for $k=n$. Next let $k\le n-1$, and assume that $t_{k+1}=t-s_{n-k}$. Thanks to the second statement in Theorem \ref{thm : proprietes jccp} (i), since $v_k$ is the mother of $v_{k+1}$,
$$
t_k= \varphi(v_k , \omega(v_k))= \sup\{s\le t_{k+1}: X_s < X_{t_{k+1}-}\}.
$$
Elementary calculations yield
$$
\sup\{s\le t_{k+1}: X_s < X_{t_{k+1}-}\} = t-\inf\{s\ge t-t_{k+1}: X_s^{(t)} > X_{t-t_{k+1}}^{(t)}\}.
$$
But since $t-t_{k+1}=s_{n-k}$ by assumption, we deduce that
$$
t_k=t-\inf\{s\ge s_{n-k}: X_s^{(t)} > X_{s_{n-k}}^{(t)}\}=t-s_{n-k+1}.
$$
As for (ii), if $v \not \in\partial{\cal T}$, just use (i) to see that the number of records of $X^{(t)}$ is the number of ancestors of $\varphi^{-1}(t)$ in $\cal T$, which is also its height. Then check that the records of the past supremum of $X^{(t)}$ are those of \label{38} the future infimum $J^{(t)}$ of $X$ on $[0,t]$, that is,
$$
J^{(t)}_s= \inf_{s\le r\le t} X_r\qquad 0\le s\le t.
$$
Now if $v\in\partial{\cal T}$, then on the one hand $v$ is infinite so by definition $H_t=\infty$; on the other hand, apply (ii) to $v\vert n+1$: write $t_n=\varphi (v\vert n, \alpha(v \vert n+1))$ and recall that $X_s\ge X_{t_{n}}$ for $s\in[t_{n},t]$, so that the future infimum $J^{(t)}$ of $X$ on $[0,t]$ has at least $n$ records, and let $n\tendinfty$.
\hfill $\Box$

\section{Splitting trees}

In this section, we consider random chronological trees, called \emph{splitting trees}, and whose width process, in the locally finite case, is a \emph{binary, homogeneous Crump--Mode--Jagers process}.

\subsection{Definition}

A splitting tree is a random chronological tree characterized by a $\sigma$-finite measure $\Lambda$ on $(0,\infty]$ called the \emph{lifespan measure}, satisfying
$$
\int_{(0,\infty]} (r\wedge 1) \Lambda(dr) <\infty.
$$
Let $\PP_\chi$ denote the law of a splitting tree starting with one ancestor individual $\emptyset$ having deterministic lifetime \label{40}
$(0, \chi]$, where $\chi$ is only allowed to equal $\infty$ when $\Lambda(\{+\infty\})>0$.

We give a recursive characterization of the family of probability measures $\PP=(\PP_\chi)_{\chi\ge 0}$ as follows. 
Recall from the Preliminaries on chronological trees that $g(\TT', \TT, x,i)$ is the tree obtained by grafting $\TT'$ on $\TT$ at $x$, as descending from $p_1(x)i$.

Let $(\alpha_i, \zeta_i)_{i\ge 1}$ be the atoms of a Poisson measure on $(0,\chi)\times (0, +\infty]$ with intensity measure $\mbox{Leb}\otimes\Lambda$ (where `Leb' stands for Lebesgue measure). Then $\PP$ is the unique family of probability measures on chronological trees $\TT$ satisfying
$$
\TT = \bigcup_{n\ge 1} g(\TT_n,\emptyset\times(0,\chi), (\emptyset,\alpha_n), n),
$$
where, conditionally on the Poisson measure, the $(\TT_n)$ are \emph{independent} splitting trees, and for each integer $n$, conditional on $\zeta_n=\zeta$, $\TT_n$ has law $\PP_{\zeta}$. 

In other words, for each individual $v$ of the tree, conditional on $\alpha(v)$ and $\omega(v)$, the pairs $(\alpha(vi), \zeta(vi))_{i\ge 1}$ made of the \emph{birth levels} and \emph{lifespans} of $v$'s offspring are the atoms of a Poisson measure on $(\alpha(v), \omega(v))\times(0, +\infty]$ with intensity measure $\mbox{Leb}\otimes\Lambda$. In addition, conditionally on this Poisson measure, descending subtrees issued from these offspring are independent.

\label{41}
Here, we have to say a word about the order of siblings. First, note that the right-hand side of the last display can also be understood as the sequential grafting of trees on $\emptyset\times(0,\chi)$ as $n$ increases. Second, as a referee pointed out, we need to specify how we label the atoms $(\alpha_i, \zeta_i)_{i\ge 1}$ to characterize $\PP$. In the case when $\Lambda(\{+\infty\})=0$, all lifespans are finite a.s., and we can assume that atoms can be ranked in such a way that for any two distinct integers $i,j$,
\begin{equation}
\label{eqn : order of siblings 1}
i<j\Longleftrightarrow \zeta_i>\zeta_j \mbox{ or } (\zeta_i=\zeta_j \mbox{ and } \alpha_i<\alpha_j).
\end{equation}
Note that this order carries over to the case when $\Lambda(\{+\infty\})>0$ but $\chi<\infty$.
In the case when $\Lambda(\{+\infty\})>0$ and $\chi=\infty$, we fix a bijection $h$ from $\NN$ onto $\NN\times\ZZ^+$, and we proceed into three steps. First, all atoms $(\alpha,\zeta)$ with infinite lifespan (second marginal) are ranked in increasing order of their birth date (first marginal). Then we call $i$-th cluster of atoms, the subset of those atoms with finite lifespan whose birth date is between the birth dates of the $(i-1)$-th and $i$-th atoms with infinite lifespan (birth date of $0$-th atom is $0$). Second, for each cluster separately, we can (and do) rank the atoms of the cluster in the order defined in \eqref{eqn : order of siblings 1}. Third and last, writing $h(n)=(h_1(n), h_2(n))$, we label atoms altogether in such a way that 
\begin{equation}
\label{eqn : order of siblings 2}
(\alpha_n,\zeta_n)\mbox{ has rank }h_2(n)\mbox{ in the }h_1(n)\mbox{-th cluster},
\end{equation}
where the atom with rank $0$ in the $i$-th cluster is the $i$-th atom with infinite lifespan itself.\\
\\
Check that when $\Lambda$ is \emph{finite}, the splitting tree is locally finite a.s.  The terminology `\emph{splitting trees}' is usually restricted to this case \cite{GK}. We will say a few more words about it in the next subsection.

\paragraph{Two branching processes.}
Recall that $\Xi_\tau$ denotes the number of individuals alive at $\tau$ (width of $\TT$ at level $\tau$) :
$$
\Xi_\tau = \mbox{Card}\{ v\in {\cal T} : \alpha(v) <\tau \le \omega(v)\}\le \infty.
$$
Under $\PP_\chi$, the width process $\Xi$ is a branching process: allowing several ancestor individuals with lifespans equal to $\chi$ and i.i.d. descendances, the branching property holds as a function of the number of ancestors. Unless $\Lambda$ is exponential or a Dirac mass at $\{\infty\}$, this branching process is \emph{not} Markovian.

Next set $Z_n$ the sum of all lifespans of individuals belonging to generation $n$ :
$$
Z_n:=\sum_{v\in{\cal T}:\vert v\vert =n}\zeta(v)=\lbd(\{x\in\TT:\vert p_1(x)\vert =n\}).
$$
Under $\PP_\chi$, $Z$ is another branching process, which itself is Markovian, as seen in the next statement.
\begin{thm}
\label{thm : Jirina} 
Under $\PP_\chi$,
$(Z_n;n\ge0)$ is a Jirina process starting from $\chi$, with branching mechanism $F$ given by
$$
F(\lbd):=\int_0^\infty (1-e^{-\lbd r})\,\Lambda(dr)\qquad\lbd\ge 0.
$$
In addition, $\TT$ has locally finite length a.s. under $\PP$, so $\Xi_\tau$ is a.s. finite for Lebesgue-a.e. $\tau$,  and the following events coincide a.s.
\begin{itemize}
\item[(i)] $\TT$ has finite length
\item[(ii)] $\lim_{n\tendinfty} Z_n =0$
\item[(iii)] $\Xi_\tau =0$ for all sufficiently large $\tau$. 
\end{itemize}
\end{thm}
\begin{rem}
The last statement allows to use the same terminology for splitting trees as for Jirina processes. For example, recall from the Preliminaries that the event defined a.s. equivalently by (i), (ii) and (iii) is called \emph{extinction} and denoted $\mathrm{Ext}$. So when $F(0)=0$, we set $m:=F'(0^+)=\intgen r\Lambda(dr)$, and we say that $\TT$ is \label{18b} subcritical, critical or supercritical according whether $m<1$, $=1$ or $>1$. When $F(0)\not=0$ (that is, $\Lambda(\{+\infty\})\not=0$), $\TT$ is said supercritical. 
\end{rem}

\paragraph{Proof.}
By construction, 
$$
Z_1=\sum_{s\le \chi}\Delta_s^1,
$$
where $(\Delta_s^1; s\ge 0)$ is a Poisson point process with intensity $\mbox{Leb}\otimes\Lambda$. Since $\int_0^1 r\Lambda (dr)<\infty$ by assumption, $Z_1$ has the value at time $\chi$ of a subordinator $S$ with Laplace exponent $F$. 
Now we reason by induction on the generation number $n$. By definition, conditionally on the knowledge of generation $n$, 
$$
Z_{n+1}=\sum_{v:\vert v\vert=n}\;\sum_{s\le \zeta(v)}\Delta_s^{(v)},
$$
where the point processes $(\Delta_s^{(v)}; s\ge 0)$ are i.i.d. Poisson point processes with common intensity $\mbox{Leb}\otimes\Lambda$. It is then a standard property of Poisson measures that, since $Z_n=\sum_{v:\vert v\vert=n}\zeta(v)$, we can write
$$
Z_{n+1}=\sum_{s\le Z_n}\Delta_s^n,
$$
where $(\Delta_s^n; s\ge 0)$ is a Poisson point process with intensity $\mbox{Leb}\otimes\Lambda$.
Then by induction on $n$, $Z_n$ has the law of the $n$-th composition (in Bochner's sense) of i.i.d. subordinators distributed as $S$, and evaluated at $\chi$. This is precisely saying that $Z$ is a Jirina process starting from $\chi$ with branching mechanism $F$ (see Preliminaries).\\
\\
To show that $\TT$ has locally finite length a.s. we are going to use Proposition \ref{prop : finite length}. Set $\Xi_\tau^{(n)}$ the number of individuals from generation $n$ living at time $\tau$
$$
\Xi_\tau^{(n)}:= \mbox{Card}\{ x\in {\TT} : \vert p_1(x)\vert=n , p_2(x)=\tau \}.
$$
Observe that $\Xi_\tau^{(0)}=\indic{\chi\ge \tau}<\infty$, and let $n\ge 0$.
By definition, conditional on $(\Xi_\sigma^{(n)};\sigma\le \tau)$, 
$$
\Xi_\tau^{(n+1)}=\sum_i \indic{\alpha_i<\tau\le \alpha_i+\zeta_i} ,
$$
where $(\alpha_i, \zeta_i)_{i\ge 1}$ are the atoms of a Poisson measure on $(0,\tau)\times (0,\infty]$ with intensity measure $\Xi_\sigma^{(n)}d\sigma \otimes \Lambda$. This entails that the conditional distribution of $\Xi_\tau^{(n+1)}$ is Poisson with parameter 
$$
\int_0^\tau d\sigma\,\Xi_\sigma^{(n)}\int_{[\tau-\sigma,\infty]} \Lambda (dr).
$$
With the notations $\bar{\Lambda}(r):=\Lambda([r,\infty])$ and $f_n(\tau):=\EE_\chi(\Xi_\tau^{(n)})$, we thus get
\begin{equation}
\label{eqn : recurrence generation}
f_{n+1}(\tau)= \int_0^\tau d\sigma\,f_n(\sigma)\bar{\Lambda}(\tau-\sigma).
\end{equation}
Next denote by $G_n$ the Laplace transform of $f_n$ \label{42}
$$
G_n(\mu) := \intgen  f_n(t)\, e^{-\mu t}\,dt\qquad \mu\ge 0,
$$
and assume that the series of generic term $G_n(\mu)$ converges for some $\mu$. Then Fubini's theorem ensures that $t\mapsto f(t):=\sum_n f_n(t)$ is integrable against $t\mapsto e^{-\mu t}$. In particular, for any $\tau>0$, $\int_0^\tau f(t)\, dt\le e^{\mu\tau} \int_0^\tau f(t)\, e^{-\mu t}\, dt \le e^{\mu\tau} \intgen f(t)\, e^{-\mu t} \,dt<\infty$.  
 Now applying Fubini's theorem to  $\Xi_\sigma = \sum_n \Xi_\sigma^{(n)}$, we get $\EE_\chi (\Xi_\sigma) =f(\sigma)$. A third and last application of Fubini's theorem, along with equation \eqref{eqn : longueur troncation} in Proposition \ref{prop : finite length}, yields
$$
\EE_\chi\lbd(C_\tau (\TT)) = \EE_\chi \int_0^\tau \Xi_\sigma\, d\sigma = \int_0^\tau f(\sigma)\, d\sigma <\infty\qquad\mbox{ for any }  \tau>0.
$$
In conclusion, it is sufficient that $\sum_n G_n (\mu)<\infty $ for some $\mu$ to ensure that $\TT$ has locally finite length a.s. (and, actually, has even integrable local length). Now we prove that it is indeed the case that $\sum_n G_n (\mu)<\infty $ for some $\mu$. Thanks to \eqref{eqn : recurrence generation}, we have
$$
G_{n+1}(\mu) = \intgen dt\,e^{-\mu t} \int_0^t ds\,f_n(s)\,\bar{\Lambda}(t-s) =  
%\intgen ds\,f_n(s)\int_s^\infty dt\, e^{-\mu t} \,\bar{\Lambda}(t-s)= 
G_n(\mu) \intgen dt\, e^{-\mu t} \,\bar{\Lambda}(t).
$$
Now a straightforward calculation shows that $\intgen dt\, e^{-\mu t} \,\bar{\Lambda}(t)= F(\mu)/\mu$, so that
$$
G_n(\mu) = \left(\frac{F(\mu)}{\mu}\right)^n \, G_0(\mu)\qquad \mu>0,
$$
and $G_0(\mu)=\int_0^\chi e^{-\mu t}\, dt<\infty$ for any $\mu\ge 0$. Now recall from the Preliminaries that there is some finite nonnegative real number $\eta$ such that $\mu >\eta\Leftrightarrow F(\mu) <\mu$, which entails the convergence of $\sum_n G_n(\mu)$ for all $\mu >\eta$.\\
\\
We finish the proof with the a.s. equality of the events defined in (i), (ii) and (iii).
Because $\lambda(\TT)=\sum_{n\ge 0} Z_n$ and $Z$ is a Jirina process, we can apply Proposition \ref{prop : properties Jirina} to get the a.s. equality of the events defined in (i) and (ii). The a.s. equality between those defined in (i) and (iii) comes from an appeal to Proposition \ref{prop : finite length}, recalling that $\TT$ has a.s. locally finite length.\hfill $\Box$

\subsection{The finite case}

\label{43}

The \emph{Crump--Mode--Jagers process}, or \emph{CMJ process}, is sometimes called the general branching process (see \cite{Taib} for a complete overview), because it counts the size of a branching population defined under the very general assumption that `reproduction schemes' of all individuals are i.i.d.. Specifically, the reproduction scheme of an individual with given birth date $\alpha$ is the joint knowledge of 
\begin{itemize}
\item her lifespan $\zeta\in(0,\infty]$ %(or her death date $\omega=\alpha+\zeta$)
\item the successive ages $(0<)\ \sigma_1 < \sigma_2<\cdots(<\zeta)$ at which she gives birth
\item
the (integer) size $\xi_i$ of the clutch she begot at the $i$-th birth time $\alpha+\sigma_i$, for all $i\ge 1$.
\end{itemize}
It is clear that a branching population with i.i.d. reproduction schemes can be constructed recursively starting from a finite number of individuals with given birth dates. 

To stick to the present framework, we further assume that
\begin{enumerate}
\item all clutch sizes (the $\xi_i$'s) are a.s. equal to $1$ (\emph{binary splitting})
	\item conditional on the lifespan $\zeta$, the point process $(\sigma_i)$ is a \emph{Poisson point process} on $(0,\zeta)$ with intensity $b$ (\emph{homogeneous} reproduction scheme)
	\item the common distribution of lifespans is $\Lambda(\cdot)/b$, where $\Lambda$ is some positive measure on $(0,\infty]$ with mass $b$ called the \emph{lifespan measure}.
\end{enumerate}
In other words, each individual gives birth at rate $b$ during her lifetime $(\alpha, \omega]$, to independent copies of herself whose lifespan common distribution is $\Lambda(\cdot)/b$. Check that this definition of a splitting tree is exactly the same as that given previously, in the special case when $\Lambda$ is finite.
In that case, the width process $(\Xi_\tau; \tau\ge 0)$, is a \emph{homogeneous binary Crump--Mode--Jagers process}.

Trees satisfying the foregoing assumptions could also be called general binary trees with constant birth rate, or homogeneous binary Crump--Mode--Jagers trees, but we will stick to the name of `\emph{splitting trees}'. On the other hand, this terminology is unfortunate because it evokes renewing binary fission (Yule tree).

\paragraph{A third branching process.} In the finite case, one can define ${\cal Z}_n$ the number of individuals belonging to generation $n$
$$
{\cal Z}_n := \mbox{Card}\{ v\in {\cal T} : \vert v\vert =n \}.
$$
From the definition, for any individual $v$, the total offspring number of $v$, conditional on $\zeta(v)=z$, is a Poisson random variable with parameter $bz$. Then it is easy to see that $({\cal Z}_n; n\ge 0)$ is a Bienaymé--Galton--Watson process started at 1, with offspring generating function $f$
$$
f(s) = \int_{(0,\infty)}b^{-1}\Lambda(dz)\ e^{-bz(1-s)}\qquad s\in[0,1],
$$
and the \emph{per capita} mean number of offspring can easily be computed to be equal to $m$ (recall $m=\intgen r\Lambda(dr)$).

\paragraph{Remark for modeling purpose.}
This birth--death scheme can be seen alternatively as a constant birth intensity measure $b\,\mathrm{Leb}$ combined with an \emph{age-dependent death intensity measure} $\mu$ given by
$$
\PP(\zeta\ge z) = \exp [-\mu((0,z))] \qquad z>0,
$$
which forces the equality
$$
\exp [-\mu((0,z))]= \bar{\Lambda}(z)/b\qquad z>0,
$$
where $\bar{\Lambda}(z)=\Lambda([z,\infty])$, and yields the following equation for $\mu$
$$
\mu(dz) = \frac{\Lambda(dz)}{\bar{\Lambda}(z)}\qquad z>0.
$$
It has to be emphasized that \label{44} in survival analysis (used in systems reliability, medical research, actuarial science, conservation biology,...), the (density of the) measure $\mu$ is called the \emph{hazard function} (probability of failing during time interval $dz$ conditional on survival up to time $z$), whereas the (density of the) measure $\Lambda/b$ is called the \emph{failure rate}.

If for some reason one has to proceed the other way round (the birth rate $b$ and some death intensity measure $\mu$ are given), notice that the lifespan measure $\Lambda$ is then
$$
\Lambda(dz) = b\mu(dz)\exp [-\mu((0,z))]\qquad z>0.
$$
The requirement that $\mu$ has to fulfill for lifespans to be a.s. finite is $\mu((0,\infty))=\infty$ (corresponding to $\Lambda(\{+\infty\})=0$). On the contrary, if $\mu$ is the null measure on $(0,\infty)$ (i.e. $\Lambda(\{+\infty\})=b$), then the CMJ process is called a \emph{Yule process} (pure-birth process with constant rate $b$).

\subsection{Law of the JCCP of a splitting tree}

From now on, we consider a splitting tree $\TT$ with lifespan measure $\Lambda$, and whose ancestor has lifespan $\zeta(\emptyset)=\chi$. Recall that its law is denoted $\PP_\chi$.

From Theorem \ref{thm : Jirina}, we know that $\TT$ has locally finite length a.s., and finite length iff extinction occurs. In particular, $\TT$ has a JCCP on $\mbox{Ext}$, denoted by $(X_t, t\in[0,\lbd(\TT)])$, and any of its finite truncations $C_\tau(\TT)$ has a JCCP as well, denoted  by $(\Xtaut, t\in[0,\lbd(C_\tau(\TT))])$. Although similar with the notation used in Corollary \ref{cor : height process}, this one has a totally different meaning.

Also set 
$$
\overline{\Xi}_\tau:=\mbox{Card}\{x\in\overline{\TT}:p_2(x)=\tau\}
\qquad\tau\ge 0.
$$
Observe that for any $\tau'\ge \tau$, the set $\{x\in\overline{\TT}:p_2(x)=\tau\}$ is equal to $\{x\in\overline{C_{\tau'}(\TT)}:p_2(x)=\tau\}$, so that
$$
\overline{\Xi}_\tau=\mbox{Card}\{t:X^{(\tau')}_t=\tau\}\qquad0\le\tau\le\tau',
$$
which does not depend on $\tau'$. \label{46}
\begin{lem}
\label{lem : Xi fini}
For any fixed $\chi,\tau$, $\PP_\chi(\overline{\Xi}_\tau<\infty)=1$.
\end{lem}
An important consequence of this lemma is that for any \emph{given} $\chi$ and $\tau\le \tau'$, both $X$ (on extinction) and $X^{(\tau')}$ hit $\tau$ a finite number of times a.s. under $\PP_\chi$.

\paragraph{Proof.} Assume that there is $\chi,\tau$, such that $b_{\chi,\tau}:=\PP_\chi(\overline{\Xi}_\tau=\infty)>0$. Then observe that, by the definition of splitting trees, $\chi\mapsto b_{\chi,\tau}$ is nondecreasing on $[0,\tau]$ and constant on $[\tau,+\infty)$. This ensures that $b_{\tau,\tau}>0$. But using the same kind of argument, as well as translation invariance, we get $b_{\sigma,\sigma}\ge b_{\tau,\tau}$ for any $\sigma\ge\tau$ and then $b_{\chi',\sigma}\ge b_{\tau,\tau}$ for any $\chi'\ge \sigma\ge\tau$. As a consequence, we get that for any $\chi'>\tau$,
$$
\EE_{\chi'}\int_0^\infty \indic{\overline{\Xi}_\sigma=\infty}\,d\sigma \ge  (\chi'-\tau)b_{\tau,\tau} >0.  
$$
Now observe that for any Borel set $A\in\overline{\TT}$, $\mbox{Leb} (p_2(A))\le \lbd(A)$. Moreover, $\overline{\Xi}_\sigma=\infty$ only if ${\Xi}_\sigma=\infty$ or $\sigma\in p_2(\partial\TT)$. But on the one hand, thanks to Theorem \ref{thm : Jirina}, $\{\sigma : \Xi_\sigma=\infty\}$ has a.s. zero Lebesgue measure, and on the other hand, $\partial \TT$, hence $p_2 (\partial\TT)$, also have zero Lebesgue measure. This is in contradiction with the last display.\hfill$\Box$\\
\\
We denote by $Y$ the spectrally positive Lévy process $t\mapsto Y_t:= -t+\sum_{s\le t}\Delta_s$, where $(\Delta_t, t\ge 0)$ is a Poisson point process with intensity measure $\mathrm{Leb}\otimes\Lambda$. In particular, $Y$ is a Lévy process with finite variation, whose Laplace exponent (see Preliminaries) will be denoted by $\psi$ 
$$
\psi(\lbd):=\lbd-F(\lbd)=\lbd-\intgen (1-\exp(-\lbd r))\ \Lambda (dr)\;\qquad\; \lbd\ge 0.
$$ 
The following statement is the fundamental result of this section. It is a little bit surprising at first sight, in the sense that, eventhough $(\varphi^{-1}(t);t\ge0)$ is not Markovian, its second projection is.
Recall that $T_A$ is the first hitting time of $A$.

\begin{thm}
\label{thm : jccp}
The law of $\Xtau$ is characterized by (i); conditional on $\mbox{\rm Ext}$, the law of $X$ is characterized by (ii).\\ 
\indent
(i)
Define recursively $t_0=0$, and $t_{i+1}=\inf\{t> t_i : \Xtaut \in\{0,\tau\}\}$. Then under $\PP_\chi$, the killed paths $e_i:= (\Xtauf{t_i+t}, 0\le t< t_{i+1}-t_i)$, $i\ge 0$, form a sequence of i.i.d. excursions, distributed as the Lévy process $Y$ killed at $T_0\wedge T_{(\tau,+\infty)}$, ending at the first excursion hitting 0 before $(\tau,+\infty)$. These excursions all start at $\tau$, but the first one, which  starts at $\min(\chi,\tau)$. In other words, $\Xtau$ has the law of $Y$ reflected below $\tau$ and killed upon hitting 0.\\
\\
\indent
(ii)
Under $\PP_\chi(\cdot \mid \mathrm{Ext})$, $X$ has the law of the Lévy process $Y$, started at $\chi$, conditioned on, and killed upon, hitting 0.  
\end{thm}

\paragraph{Proof.}
(i) \label{48}
In what follows, we will stick to the notation $\varphi^{-1}$ for the exploration process of $C_\tau(\TT)$ and we define $\ell:=\lbd(C_\tau(\TT))$. Then let $t\in[0,\ell)$ and $x=\varphi^{-1}(t)\in C_\tau(\overline{\TT})=\overline{C_\tau(\TT)}$. Set $u=p_1(x)\in\overline{\cal T}$ and $n=\vert u\vert \le\infty$.  Set also $u_k=u\vert k$. The segment $[\rho, x[$ of ancestor points of $x$ (except itself) is the union of segments $u_k\times [\alpha(u_k),\alpha(u_{k+1}))$, for $0\le k \le n$, with the convention $\alpha(u_{n+1})=p_2(x) = X_t^{(\tau)}$ when $n<\infty$. Next denote by $\tau_{k,i}\in(\alpha(u_{k}),\omega(u_{k}))$ the birth level of $u_ki$, and by $I_k$ the set of integers $i$ such that $\tau_{k,i}<\alpha(u_{k+1})$, that is, $u_k i$ is a child of $u_k$ born \emph{before} $u_{k+1}$. Finally, we denote by $\TT_{k,i}$ the subtree that can be seen as grafted at $\tau_{k,i}$, that is, 
$$
\TT_{k,i}:=\{(v,\sigma):(u_k iv,\sigma+\tau_{k,i})\in\theta_r ((u_k, \tau_{k,i}))\}\qquad 0\le k\le n, i\in I_k. 
$$
Now, we can write the subtree $\TT_{\mathrm{post-}t}^{(\tau)}$ of points of $C_\tau(\TT)$ still unvisited at time $t$ (which is the union of $[\rho,x[$ and $R(x)$ in $C_\tau(\TT)$) as
$$
\TT_{\mathrm{post-}t}^{(\tau)}:= C_\tau\left(\bigcup_{0\le k\le n, i\in I_k} g\left(\TT_{k,i}, [\rho,x[,(u_k,\tau_{k,i}), i\right)\right) .
$$
We point out that because truncation preserves the order, the post-$t$ exploration process of $C_\tau(\TT)$ is exactly the exploration process of $\TT_{\mathrm{post-}t}^{(\tau)}$.\\ 

Now we work under $\PP_\chi$. 
We call \emph{heightwise label transposition} any mapping $\theta: \UU\rightarrow \UU$, for which there are integers $k,i,j$ such that for any $y\in\UU$, written as $y=(v,\sigma)$ with $v=(v_1,\ldots, v_l)$,
$$
\theta(y)=\left\{
\begin{array}{cll}
(v,\sigma) & \mbox{ if } & l <k \mbox{ or } v_k\not\in\{i,j\}\\
(v^\star,\sigma)& \mbox{ if } &l \ge k \mbox{ and } v_k = j\\
(v^{\star\star},\sigma) &\mbox{ if } & l\ge k \mbox{ and } v_k = i,
\end{array} 
\right.$$
where $v^{\star}$ (resp. $v^{\star\star}$) is obtained from $v$ by setting $v_k$ equal to $i$ (resp. to $j$), and leaving other labels unchanged. Observe that heightwise label transpositions map chronological trees into chronological trees.%, and that the genealogical structures, the Lebesgue measure and the linear order are invariant under their action.

For any $r\ge 0$, let ${\cal F}_r$ denote the $\sigma$-field generated by $\{F(\varphi^{-1}(s)); s\le r, F \in {\cal H} \}$, where $\cal H$ denotes the set of measurable (real) functions on $\UU$ that are invariant under the action of \emph{all} heightwise label transpositions. This filtration is chosen to contain all past events that do not depend upon the choice of sibling labelling. Notice that since $p_2\in{\cal H}$,  $(X_s; s\le r)$ is ${\cal F}_r$-measurable, so it is sufficient to prove that $X^{(\tau)}$ is Markovian w.r.t. the filtration $({\cal F}_r)_{r\ge 0}$.

Now for all $0\le k\le n$, let $\zeta_{k,i}$ be the lifespan of the ancestor of $\TT_{k,i}$, and let $\mu_k$ be the random point measure on $(\alpha(u_{k}),\alpha(u_{k+1}))\times (0,\infty]$ with atoms $((\tau_{k,i},\zeta_{k,i}),i\in I_k)$.\\ 

\textbf{Claim 1.} For any $k=0,...n$, conditional on $\{\alpha(u_k),\alpha(u_{k+1})\}$, $\mu_k$ is a Poisson point measure on $(\alpha(u_k),\alpha(u_{k+1}))\times (0,\infty]$ with intensity measure $\mbox{Leb}\otimes \Lambda$, independent of ${\cal F}_t$.\\

To see this, notice that because $i\in I_k\Leftrightarrow \tau_{k,i}<\alpha(u_{k+1})$, the point measure $\mu_k$ is the restriction to $(\alpha(u_k),\alpha(u_{k+1}))\times (0,\infty]$ of the Poisson point measure, say $\nu_k$, on $(\alpha(u_k),\omega(u_{k}))\times (0,\infty]$ with intensity measure $\mbox{Leb}\otimes \Lambda$. Now we are reasoning conditional on $\{\alpha(u_k),\alpha(u_{k+1})\}$. First, the point  measure $\mu_k$ (conditionally) is a Poisson point measure on $(\alpha(u_k),\alpha(u_{k+1}))\times (0,\infty]$ with intensity measure $\mbox{Leb}\otimes \Lambda$. Second, as a consequence of  the recursive construction of $\TT$, $\mu_k$ depends on ${\cal F}_t$ only through the knowledge of $((\tau_{k,i},\zeta_{k,i}),i\not\in I_k)$. But since ${\cal F}_t$ is independent of sibling labelling, $\mu_k$ depends on ${\cal F}_t$ only through the knowledge of the atoms of the restriction of $\nu_k$ to $[\alpha(u_{k+1}),\omega(u_{k}))\times (0,\infty]$. As a consequence, the point measure $\mu_k$ is (conditionally) independent of ${\cal F}_t$. %Since one can relabel the atoms by the integers, e.g. in the order specified by \eqref{eqn : order of siblings 1} or \eqref{eqn : order of siblings 2}, without changing the (conditional) law of $\mu_k$, the claim follows.
\\

\textbf{Claim 2.} Conditional on $(\alpha(u_k); 0\le k\le n)$, the point measures ($\mu_k; 0\le k\le n)$ are independent.\\

This last claim is a mere consequence of the recursive construction of $\TT$, and of the fact that all points $((u_k, \sigma);0\le k\le n, \sigma\in(\alpha(u_{k}),\alpha(u_{k+1})))$ have disjoint descendances. 

Now rank the atoms $((\tau_{k,i},\zeta_{k,i}); 0\le k\le n, i\in I_k)$ in the order specified by \eqref{eqn : order of siblings 1} or \eqref{eqn : order of siblings 2}, according to whether $\Lambda(\{+\infty\})$ is zero or not. We denote by $(\tau_j, \zeta_j; j\ge 1)$ this relabelled set of atoms, as well as $\TT_j$ the tree corresponding to $(\tau_j, \zeta_j)$ (read $j$ as $j(k,i)$). For precisely the same reason as for Claim 2, we have\\

\textbf{Claim 3.} Conditional on $(\zeta_{j}; j\ge 1)$, the trees $\TT_{j}$ are independent and independent of ${\cal F}_t$, and for all $j\ge 1$, conditional on $\zeta_{j}=z$, $\TT_{j}$ has law $\PP_z$.\\

Next define
$$
\TT':= \bigcup_{ j}\; g\left(\TT_{j}, \emptyset\times[0,X_t^{(\tau)}[,(\emptyset,\tau_j), j\right).
$$
Thanks to Claims 1 and 2, conditional on $X_t^{(\tau)}=\sigma$, and conditional on $(\alpha(u_k); 0\le k\le n)$, the point measure $\mu:=\sum_{0\le k\le n}\mu_k$ is the superposition of independent Poisson point measures, all independent of ${\cal F}_t$, all with the same translation invariant intensity $\mbox{Leb}\otimes \Lambda$, on disjoint sets whose union is equal to $(0,\sigma)\times (0,\infty]$, up to the discrete set $(\alpha(u_k); 0\le k\le n)$. As a consequence, $\mu$ is a Poisson point measure with intensity $\mbox{Leb}\otimes \Lambda$ on $(0, \sigma)\times (0,\infty]$, which does not depend on $\{\alpha(u_{k+1}): 0\le k\le n\}$, and is hence independent of ${\cal F}_t$ conditional on $X_t^{(\tau)}=\sigma$. Now thanks to Claim 3, conditional on $X_t^{(\tau)}=\sigma$, $\TT'$ is independent of ${\cal F}_t$ and has law $\PP_\sigma$.

Now observe that any $x\in\TT_{\mathrm{post-}t}^{(\tau)}$ has $p_1(x)=u_kiv$ for some $v$, and that the mapping from $\TT_{\mathrm{post-}t}^{(\tau)}$ to $C_\tau(\TT')$ that maps $u_kiv$ into $jv$ ($j=j(k,i)$), and leaves $p_2(x)$ unchanged is a bijection which preserves the linear order and the second projection. As a consequence, the JCCPs of both trees are equal. This shows that conditional on $X_t^{(\tau)}=\sigma$, $(X_{t+s}^{(\tau)};s\ge 0)$ is independent of ${\cal F}_t$ and has the distribution of $(X_{s}^{(\tau)};s\ge 0)$ under $\PP_\sigma$.

Thus, we have shown that $X^{(\tau)}$ is Markovian, and since $\ell=\inf\{s\ge 0: X_{s}^{(\tau)}=0\}$, we get that $(X^{(\tau)}; 0\le s <\ell)$ is a Markov process killed upon hitting $0$.\\
\\ 
We want to prove that $X^{(\tau)}$ is a Feller process. Let $f$ be a continuous, hence bounded, function on $[0,\tau]$. 

\label{52}

As a first step, because $X^{(\tau)}$ is càdlàg, the dominated convergence theorem ensures the right-continuity of $t\mapsto \EE_\chi \left(f\big(X^{(\tau)}_t\big)\right)$, which is the first required condition for $X^{(\tau)}$ to be Feller.

As a second step, let $\chi'<\chi$, and recall that $(\emptyset, \chi')$ is  $\PP_\chi$-a.s. a simple point, that is, a point with degree 2 in the chronological tree, at which the exploration process is continuous. Also notice that $\varphi((\emptyset, \chi'))$ is the first hitting time $T_{\chi'}$ of $\chi'$ by $X^{(\tau)}$. As a consequence, the contour process of the (truncation of the) chronological tree issued from $\emptyset\times [0,\chi']$ (resp. $\emptyset\times [\chi',\chi]$) is $(X^{(\tau)}_{T_{\chi'}+t};t\ge 0)$ (resp. $(\Xtau_t;0\le t\le T_{\chi'})$), and in addition both contour processes are independent, and the former has the law of $\Xtau$ under $\PP_{\chi'}$. In particular, $\Xtau$ satisfies 
the Markov property at $T_{\chi'}$ under $\PP_{\chi}$. Using these observations, we now show that $\chi\mapsto \EE_\chi\left(f\big(X^{(\tau)}_t\big)\right)$ is continuous on $[0,\tau]$, which is the other required condition for $X^{(\tau)}$ to be Feller.
First, since $T_{\chi'}$ is the length of the subtree issued from $\emptyset\times [\chi',\chi]$ and truncated at $\tau$, it is stochastically dominated by $\lbd(\TT)$ under $\PP_{\chi-\chi'}$. Putting together the first statement of Theorem \ref{thm : Jirina} with the arguments developed in the proof of Proposition \ref{prop : properties Jirina}, we see that under $\PP_x$, $\lbd(\TT)=\sum_n Z_n$ ($=:T$ in the proof of Proposition \ref{prop : properties Jirina}) converges to 0 in probability as $x\to 0$.
This ensures that $\PP_\chi(T_{\chi'}>\eps)$ vanishes as $\chi-\chi'\to 0$. Second, we have the following equality
\debeq
\EE_\chi\left(f\big(X^{(\tau)}_t\big)\right) -\EE_{\chi'}\left(f\big(X^{(\tau)}_t\big)\right)
&=&\int_0^\eps
\EE_\chi\left(f\big(X^{(\tau)}_t\big)-f\big(X^{(\tau)}_{t+s}\big),T_{\chi'}\in ds\right)\\	&+&\EE_\chi\left(f\big(X^{(\tau)}_t\big)-f\big(X^{(\tau)}_{t+T_{\chi'}}\big), T_{\chi'}>\eps\right).%\\
	%&=&\int_0^\eps\PP_\chi (T_{\chi'}\in ds)\EE_{\chi'}\left(f\big(X^{(\tau)}_{t-s}\big)-f\big(X^{(\tau)}_{t}\big)\right)\\
	%&+& \EE_\chi\left(f\big(X^{(\tau)}_t\big)-\EE_{\chi'}\left(f\big(X^{(\tau)}_{t}\big)\right),T_{\chi'}>\eps\right)
\fineq
Recalling that $f$ is bounded, using dominated convergence along with the right-continuity of $X^{(\tau)}$ for the first term, we get the continuity of $\chi\mapsto \EE_\chi\left(f\big(X^{(\tau)}_t\big)\right)$.
The conclusion (see e.g. \cite{RY}) is that $X^{(\tau)}$ is a Feller process, and as such, it satisfies the strong Markov property.\\
\\
Now fix $\tau'>\tau>0$ and define $A^{\tau,\tau'}$ as
$$
A_t^{\tau,\tau'} = \int_0^t ds\,\indicbis{X_s^{(\tau')}\le\tau}\qquad t\ge 0,
$$
and let $a^{\tau,\tau'}$ be its right inverse
$$
a_t^{\tau,\tau'} = \inf\{s : A_s^{\tau,\tau'}>t\}\qquad t\ge 0.
$$ 
Then it is clear from our construction that $X^{(\tau')}\circ a^{\tau, \tau'}$ is the JCCP of $C_\tau\circ C_{\tau'}(\TT)$. But because $\tau'>\tau$, $C_\tau\circ C_{\tau'}= C_\tau$, so that  
$$
X^{(\tau)} = X^{(\tau')}\circ a^{\tau, \tau'}.
$$
Then define recursively $r_0=0$, $s_{i}=\inf\{s \ge r_{i} :  X_s^{(\tau')}\le\tau\}$ and $r_{i+1}=\inf\{r\ge s_i :  X_r^{(\tau')}>\tau\}$. Consider the killed paths $\eps_i:=(X_{t+s_i}^{(\tau')};0\le t<r_{i+1}-s_i)$, $i\ge 0$, which are the excursions of $X^{(\tau')}$ away from $[\tau,\tau']$. Because $X^{(\tau')}$ enters continuously in $[0,\tau]$ (it has no negative jumps), the initial value of each of these excursions is $\tau$, except, if $\chi<\tau$, $\eps_0$ which starts from $\chi$. Because $X^{(\tau')}$ is a killed strong Markov process, the excursions $(\eps_i)_{i\ge 0}$ form a sequence of i.i.d. excursions (except that $\eps_0$ may have a different starting point), distributed as the process $X^{(\tau')}$ killed upon exiting $(0,\tau]$, ending at the first excursion exiting  $(0,\tau]$ from the bottom.

Then observe that $A^{\tau,\tau'}\equiv t_i$ on $[r_i, s_{i}]$ and that $a^{\tau,\tau'}({t_i-})=r_i$ and $a^{\tau,\tau'}(t_{i})=s_i$, where $t_0=0$ and $t_i$ is the $i$-th hitting time of $\tau$ by $X^{(\tau)}$ defined in the theorem, so that $\eps_i$ is also the $i$-th excursion $e_i$ of $X^{(\tau)}$ away from $\tau$, for all $i$.

As a consequence, it only remains to prove that the strong Markov process $X^{(\tau')}$  and the Lévy process $Y$ with Laplace exponent $\psi$, both started at $\chi\le \tau$ and killed upon exiting $(0,\tau]$, have the same law. From what precedes, we know that this law does not depend on $\tau'$.
Recall from Theorem \ref{thm : allure jccp} that the process $X^{(\tau')}$ started at $\chi$ and killed upon exiting $(0,\tau]$ can be written
$$
X_t^{(\tau')} = \chi -t + \sum_{i : 0<t_i\le t}\zeta(u_i),\qquad 0\le t<T_0\wedge T_{(\tau,\infty)},
$$
where $(x_i=(u_i, \omega(u_i)); i\ge 0)$ range over the leaves of $C_{\tau'}(\TT)$, and $t_i=\varphi(x_i)$. 
By a standard truncation argument, we can assume that the lifespan measure $\Lambda$ has compact support and choose $\tau'$ sufficiently large for this support to be contained in $[0,\tau'-\tau]$. 
This guarantees that on $[0,T_0\wedge T_{(\tau,\infty)}]$, the  jumps of $X^{(\tau')}$ are exactly the lifespans of individuals (as if no branches were cut down). Then the jumping rate  of $X^{(\tau')}$  from $x$ to $(x+z)dz$ on $[0,T_0\wedge T_{(\tau,\infty)}]$ does not depend on $\tau$ (recall it does not depend on $\tau'$ either), so we denote it by $M(x,dz)$. Next recall from the proof of the Feller property above, that the law of $X^{(\tau')}$ started at $\chi$ and killed upon exiting $(0,\tau]$ is the same as that of $X^{(\tau')}$ started at $\chi+h$ and killed upon exiting $(h,\tau+h]$. But the jumping rate at $0^+$ of the former is $M(\chi, dz)$, whereas that of the latter is $M(\chi+h,dz)$, which proves that $M(x,dz)$ does not depend on $x$ either. As a consequence, the law of $X^{(\tau')}$ killed upon exiting $(0,\tau]$ is that of a Lévy process (killed upon exiting $(0,\tau]$) with no negative jumps and finite variation, drift $-1$ and Lévy measure, say $\Pi$. It only remains to show that $\Pi=\Lambda$. \label{56} Thanks to Theorem \ref{thm : proprietes jccp} and Corollary \ref{cor : height process}, the dates at which the ancestor gives birth are those levels $\sigma\in(0,\chi)$ such that $D(\sigma -)<D(\sigma)$, where 
$$
D(\sigma):= \inf\{t\ge 0: X^{(\tau')}_t = \chi -\sigma\} \qquad\sigma \in (0,\chi),
$$
and furthermore, the lifespan of the individual born at level $\sigma$ is the jump size $\Delta_\sigma$ of $X^{(\tau')}$ at $D(\sigma -)$. 
%Now for any $\sigma \in (0,\chi)$ such that $D(\sigma -)<D(\sigma)$, we set $\Delta_\sigma:=\Delta X^{(\tau')}_t$, writing $t$ for $D(\sigma -)$. 
Then we deduce that $(\Delta_\sigma ; \sigma \in(0,\chi))$ is a Poisson point process with intensity measure $\Lambda$. But since $X^{(\tau')}$ is a Lévy process with no negative jumps, finite variation, drift coefficient $-1$ and Lévy measure $\Pi$, it is known that $(\Delta_\sigma ; \sigma \in(0,\chi))$ is a Poisson point process with intensity measure $\Pi$. This shows that $\Pi=\Lambda$.\\
%If $\Lambda$ is finite, then $\Pi$ also is, and since the probability that there is no jump $\exp(-\chi\Pi((0,\infty]))$ is also the probability that the ancestor has no offspring $\exp(-\chi b)$, we get $\Pi((0,\infty])= b= \Lambda((0,\infty])$. Then the first jump is the lifespan of the youngest daughter of the ancestor, so that $\Pi(dr)/b=\Lambda(dr)/b$, which ends the proof in the finite case. Now if $\Lambda$ is infinite, so is $\Pi$ and by restricting the contour to the individuals whose ancestors all have lifespan greater than $\eps$, since jumps are exactly lifespans, we deduce from the finite case that the restrictions of $\Lambda$ and $\Pi$ to $(\eps,+\infty]$ are equal. Since this holds for any positive $\eps$, $\Pi=\Lambda$.\\
\\
\indent
(ii)
Under $\PP_\chi(\,\cdot\, , \mbox{Ext})$, thanks to Theorem \ref{thm : Jirina}, we see that (i) holds  for $\tau=\infty$, and we get that under $\PP_\chi(\,\cdot\, , \mbox{Ext})$,
$X$ is a Lévy process with Laplace exponent $\psi$, started at $\chi$, \emph{hitting 0}, and killed upon hitting 0. Now by Proposition \ref{prop : properties Jirina} and Theorem \ref{thm : Jirina}, $\PP_\chi(\mbox{Ext}) =\exp(-\eta \chi)$, and from the Preliminaries on Lévy processes, we know that this is also the probability that the Lévy process $Y$ with Laplace exponent $\lbd\mapsto \lbd -F(\lbd)$ hits 0. As a consequence, the statement on \emph{conditioned} distributions also holds.\hfill$\Box$

%From  Theorem \ref{thm : Jirina}, observe that
%$$
%\{\mbox{Ext}\}\Leftrightarrow \mbox{ there is some positive $\sigma$ such that $X^{(\tau)}=X$ for all }\tau\ge\sigma.
%$$
%Then define $Y^{(\tau)}$ as the first excursion of $X^{(\tau)}$ away from $\tau$. By (i), $Y^{(\tau)}$ has the law of the Lévy process $Y$ killed upon exiting $(0,\tau]$, and for any $\tau'>\tau$, 
%$$
%Y_t^{(\tau')} = Y_t^{(\tau)} \qquad t\le\inf\{s:Y_s^{(\tau')}>\tau\}.
%$$
%As a consequence, $Y^{(\tau)}$ converges a.s. as $\tau\tendinfty$ to a path $Y^{(\infty)}$ distributed as the Lévy process $Y$ killed upon exiting $(0,\infty)$. Following the equivalence initiating this proof, $Y^{(\infty)}=X$ iff extinction occurs, which ends the proof. 

\section{New properties of splitting trees and the Crump--Mode--Jagers process}

In this section, we will constantly use the notation $F(\lbd)=\intgen (1-e^{-\lbd r})\Lambda(dr)$, and  $\psi(\lbd)=\lbd-F(\lbd)$. The scale function $W$ is the positive function with Laplace transform $1/\psi$ (see Preliminaries).

Also recall that $\PP_\chi$ denotes the law of the splitting tree with $\zeta(\emptyset)=\chi$, whereas $P_\chi$ is the law of the spectrally positive Lévy process $Y$ with Laplace exponent $\psi$ started at $\chi$.

\subsection{Properties of splitting trees}

\subsubsection{A new proof of Le Gall and Le Jan's theorem}

\begin{thm}
Let $Y$ be a spectrally positive Lévy process with Laplace exponent $\psi$ such that $\psi(0)=0$ and $\psi'(0+)\ge 0$, started at $\chi$ and killed when it hits 0. Then define
$$
H_t = \mathrm{Card}\{0\le s\le t : Y_{s-}<\inf_{s\le r\le t} Y_r^{}\}\qquad t< T_0,
$$
and $L$ the local time process of $H$
$$
L_n := \int_0^{T_0} dt\,\indic{H_t=n}  \qquad n\ge 0.
$$
Then $(L_n; n\ge 0)$ is a Jirina process with branching mechanism $F$ starting from $L_0=\chi$, that is, for any integer $n \ge 1$,
$$
L_{n} = S_n\circ \cdots \circ S_1(\chi),
$$
where the $S_i$ are i.i.d. subordinators with Laplace exponent $F$.
\end{thm}

\label{57}

\begin{rem}
In \cite{LGLJ}, the Lévy process $Y$ is the continuous analogue of a downwards-skip-free random walk that can be extracted from a Bienaymé--Galton--Watson tree as follows. The jumps of this random walk are the offspring sizes, shifted by $-1$, of successive individuals of the tree taken in the depth-first search order. The height functional $H$ is chosen by analogy with the BGW case, where it allows to recover the generation number of individuals visited in the depth-first search order. Then it is shown \cite[Proposition 3.2]{LGLJ} that the time $L_n$ spent at level $n$ indeed is Markovian (in $n$), satisfies the branching property (in its initial value), and its probability transitions are displayed. \\
The foregoing statement is a slight refinement of Proposition 3.2 in \cite{LGLJ}, because it provides the law of $(L_n;n\ge 0)$ as the composition of subordinators, which in passing sheds light on the genealogy defined in \cite{BLG} by flows of subordinators. This relationship between both types of genealogies can be seen directly thanks to splitting trees (see proof below, and Introduction p.\pageref{explication}). 
\end{rem}

\paragraph{Proof.} Thanks to the tools set up in the previous sections, the proof of this theorem is straightforward.
First, since $\psi(0)=0$, thanks to Theorem \ref{thm : Jirina} and Proposition \ref{prop : properties Jirina}, the  splitting tree $\TT$ with lifespan measure $\Lambda$ and law  $\PP_\chi$ is subcritical ($\psi'(0+)> 0$) or critical ($\psi'(0+)= 0$), and so has finite length $\ell$ a.s. Now  Theorem \ref{thm : jccp}(ii)  allows to state that $Y$ is the JCCP of $\TT$, so in particular, $T_0=\ell$. 

Next, thanks to Corollary \ref{cor : height process}, we know that $H_t$ is also the genealogical height of $\varphi^{-1}(t)$ in $\cal T$, so that
$$
L_n = \int_0^{\ell} dt\,\indic{\vert p_1\circ\varphi^{-1}(t)\vert =n}  = \lbd(\{x\in\TT : \vert p_1(x)\vert=n\}) .
$$
This proves that $L_n = \sum_{v:\vert v \vert =n}\zeta(v)=Z_n$, so the proof ends with an appeal to Theorem \ref{thm : Jirina}.
\hfill$\Box$

\subsubsection{Exceptional points}
Recall that for any chronological tree $\TT$,
$$
\Xi_\tau=\mbox{Card}\{x\in\TT : p_2(x)=\tau\}\le\mbox{Card}\{x\in\overline{\TT} : p_2(x)=\tau\}=\overline{\Xi}_\tau.
$$
Recall from Lemma \ref{lem : Xi fini} that for any \emph{fixed} $\tau$, $\Xi_\tau=\overline{\Xi}_\tau<\infty$ a.s.
\begin{lem}
\label{lem : Xi infinite}
$\PP(\forall x=(u,\tau)\in\partial \TT,\, \Xi_\tau=\infty )=1$.
\end{lem}
An important consequence of this lemma is that a.s. under $\PP$ for all $\tau'\ge \tau\ge 0$,
$$
\Xi_\tau = \overline{\Xi}_\tau = \mbox{Card}\{t : X^{(\tau')}_t=\tau\}.
$$
In other words $(\Xi_\tau ;0\le \tau\le \tau')$ is a.s. equal to the occupation process of $X^{(\tau')}$. Indeed, if we had $\Xi_\tau<\overline{\Xi}_\tau$ then there would be $x\in\partial\TT$ such that $p_2(x)=\tau$, but in that case, $\Xi_\tau=\infty$, which would imply $\overline{\Xi}_\tau=\infty=\Xi_\tau$.

\paragraph{Proof.}
Assume that there is $x=(u,\tau)\in\partial \TT$ such that $\Xi_\tau<\infty$. Since $\TT$ has a.s. locally finite length, and we are interested in a local property (property of truncations), we can also assume, without loss of generality, that it has finite length. Recall that $\tau=\lim_{n} \omega(u\vert n)$ and that $\alpha(u\vert n)<\tau$ for all $n$. Since $\Xi_\tau<\infty$, we must have $\omega(u\vert n)<\tau$ for all but a finite number of $n$'s. In particular, there are infinitely many times $t_n:=\varphi(u\vert n, \omega(u\vert n))$ increasing to $t:=\varphi(x)$ such that $X_{t_n}<\tau=X_t$ (and $\lim_n X_{t_n}=\tau$), where $X$ is the JCCP of $\TT$. Setting $U:=\sup \{s<t : X_s\ge X_t\}$, we must have $U<t$, because otherwise, since $X$ has no negative jumps, it would hit $\tau$ infinitely often before time $t$. The bottomline is that, with positive probability, $X$ has infinitely many ladder times (record times of the past supremum) on a finite interval. But since $X$ is a spectrally positive Lévy process with finite variation, this happens with zero 
probability \cite[Chapter VII]{B}.\hfill $\Box$\\
\\
Actually, we have not even proved yet that $\partial \TT$ was nonempty with positive probability under $\PP$. Of course, when $\Lambda$ is finite, $\PP(\partial \TT \not=\emptyset)=0$. The following statement treats the infinite case. Fix $\tau >0$, and set 
$$
\Gamma :=\{\sigma\in[0,\tau] : \Xi_\sigma = \infty\}\quad \mbox{ and }G:=\{s\in[0,t] : H_s = \infty\},
$$
where $H_s$ is the height of $\varphi^{-1}(s)$, $\varphi^{-1}$ being the exploration process of $C_\tau (\TT)$, and $t=\lbd(C_\tau(\TT))$.

Thanks to the paragraph preceding Theorem \ref{thm : dfn exploration}, we know that $G$ has zero Lebesgue measure. Thanks to Theorem \ref{thm : Jirina}, we know that $\Gamma$ has a.s. zero Lebesgue measure. 
\begin{thm}
\label{thm : exceptional points}
If $\Lambda$ is infinite, conditional on $\Xi_\tau\not=0$, $\Gamma$ and $G$ are a.s. everywhere dense.
\end{thm}

\paragraph{Proof.}
We start with $G$. Recall the notation $X_s^{(t)}:=X_{t-}-X_{(t-s)-}$, $0\le s\le t$. By duality, and thanks to Theorem \ref{thm : jccp}, the law of $X^{(t)}$ on finite intervals away from 0 and $-\tau$ is that of the Lévy process $Y$. Now thanks to Corollary \ref{cor : height process}, $H_s$ is the number of record times of the past supremum of $X^{(s)}$. As a consequence, to have $H_s=\infty$, it is sufficient that
$$
\inf\{r\ge 0: X^{(s)}_r>0\}=0.
$$
But $\inf\{r\ge 0: X^{(s)}_r>0\}=\inf\{r\ge s: X^{(t)}_r>X^{(t)}_s\}$, so it suffices to show that for any finite interval $I$
$$ 
P(\{s\in I: \inf\{r\ge s: Y_r>Y_s\}=0\}\mbox{ is dense})=1.
$$
On the one hand, when $\Lambda$ is infinite, the Hausdorff dimension of the set of times when the path of $Y$ has Hölder exponent $h(s)<1$, is strictly positive  a.s. \cite{J}. On the other hand, it is easily seen that
$$
P(\forall s\in I, \inf\{r\ge s: Y_r>Y_s\}>0 \Rightarrow h(s)=1)=1,
$$
which entails
$$
P(\dim_H (\{s\in I: \inf\{r\ge s: Y_r>Y_s\}=0\}) >0)=1,
$$
which in turn yields the desired result (the interval $I$ being arbitrary).

Now let us prove that $\Gamma$ is everywhere dense a.s. Conditionally on $\Xi_\tau\not=0$, because $G$ is everywhere dense and $X$ is right-continuous, the closure of the range of $G$ by $X$ is $[0,\tau]$. Since $s\in G\Leftrightarrow\varphi(s)\in\partial\TT$, we can conclude with Lemma \ref{lem : Xi infinite}.\hfill$\Box$

\subsubsection{The coalescent point process}

Fix $\tau>0$.
For any chronological tree $\TT$, we let $(x_i(\tau);1\le i \le \Xi_\tau)$ denote the ranked points $x_1\le x_2 \le \cdots$ of $\TT$ such that $p_2(x_i)=\tau$. In particular, the vertices $p_1(x_i)$ of $\cal T$ are exactly the individuals alive at level $\tau$.

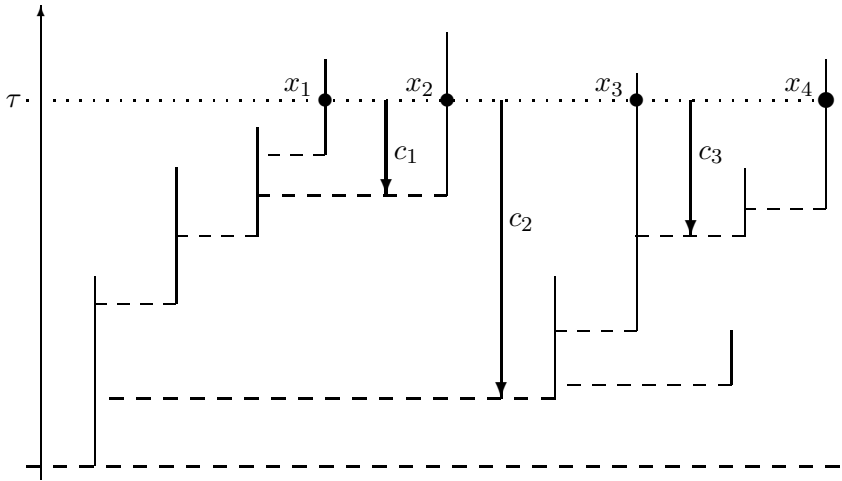
\begin{figure}[ht]

\unitlength 1.8mm % = 5.69pt
\linethickness{0.4pt}
\ifx\plotpoint\undefined\newsavebox{\plotpoint}\fi % GNUPLOT compatibility
\begin{picture}(65,45)(-4,3)
\put(6,5){\vector(0,1){35}}
%\dashline{1}(5,6)(65,6)
\put(4.96,5.96){\line(1,0){.984}}
\put(6.92,5.96){\line(1,0){.984}}
\put(8.89,5.96){\line(1,0){.984}}
\put(10.86,5.96){\line(1,0){.984}}
\put(12.82,5.96){\line(1,0){.984}}
\put(14.79,5.96){\line(1,0){.984}}
\put(16.76,5.96){\line(1,0){.984}}
\put(18.73,5.96){\line(1,0){.984}}
\put(20.69,5.96){\line(1,0){.984}}
\put(22.66,5.96){\line(1,0){.984}}
\put(24.63,5.96){\line(1,0){.984}}
\put(26.6,5.96){\line(1,0){.984}}
\put(28.56,5.96){\line(1,0){.984}}
\put(30.53,5.96){\line(1,0){.984}}
\put(32.5,5.96){\line(1,0){.984}}
\put(34.46,5.96){\line(1,0){.984}}
\put(36.43,5.96){\line(1,0){.984}}
\put(38.4,5.96){\line(1,0){.984}}
\put(40.37,5.96){\line(1,0){.984}}
\put(42.33,5.96){\line(1,0){.984}}
\put(44.3,5.96){\line(1,0){.984}}
\put(46.27,5.96){\line(1,0){.984}}
\put(48.23,5.96){\line(1,0){.984}}
\put(50.2,5.96){\line(1,0){.984}}
\put(52.17,5.96){\line(1,0){.984}}
\put(54.14,5.96){\line(1,0){.984}}
\put(56.1,5.96){\line(1,0){.984}}
\put(58.07,5.96){\line(1,0){.984}}
\put(60.04,5.96){\line(1,0){.984}}
\put(62.01,5.96){\line(1,0){.984}}
\put(63.97,5.96){\line(1,0){.984}}
%\end
\put(10,20){\line(0,-1){14}}
\put(16,28){\line(0,-1){10}}
\put(22,23){\line(0,1){8}}
\put(27,29){\line(0,1){7}}
\put(36,26){\line(0,1){12}}
\put(44,11){\line(0,1){9}}
\put(50,16){\line(0,1){19}}
%\dashline{1}(16,18)(10,18)
\put(15.96,17.96){\line(-1,0){.857}}
\put(14.24,17.96){\line(-1,0){.857}}
\put(12.53,17.96){\line(-1,0){.857}}
\put(10.81,17.96){\line(-1,0){.857}}
%\end
%\dashline{1}(22,23)(16,23)
\put(21.96,22.96){\line(-1,0){.857}}
\put(20.24,22.96){\line(-1,0){.857}}
\put(18.53,22.96){\line(-1,0){.857}}
\put(16.81,22.96){\line(-1,0){.857}}
%\end
%\dashline{1}(27,29)(22,29)
\put(26.96,28.96){\line(-1,0){.833}}
\put(25.29,28.96){\line(-1,0){.833}}
\put(23.62,28.96){\line(-1,0){.833}}
%\end
%\dashline{1}(36,26)(22,26)
\put(35.96,25.96){\line(-1,0){.933}}
\put(34.09,25.96){\line(-1,0){.933}}
\put(32.22,25.96){\line(-1,0){.933}}
\put(30.36,25.96){\line(-1,0){.933}}
\put(28.49,25.96){\line(-1,0){.933}}
\put(26.62,25.96){\line(-1,0){.933}}
\put(24.76,25.96){\line(-1,0){.933}}
\put(22.89,25.96){\line(-1,0){.933}}
%\end
%\dashline{1}(50,16)(44,16)
\put(49.96,15.96){\line(-1,0){.857}}
\put(48.24,15.96){\line(-1,0){.857}}
\put(46.53,15.96){\line(-1,0){.857}}
\put(44.81,15.96){\line(-1,0){.857}}
%\end
%\dashline{1}(44,11)(10.13,11)
\put(43.96,10.96){\line(-1,0){.996}}
\put(41.96,10.96){\line(-1,0){.996}}
\put(39.97,10.96){\line(-1,0){.996}}
\put(37.98,10.96){\line(-1,0){.996}}
\put(35.99,10.96){\line(-1,0){.996}}
\put(33.99,10.96){\line(-1,0){.996}}
\put(32,10.96){\line(-1,0){.996}}
\put(30.01,10.96){\line(-1,0){.996}}
\put(28.02,10.96){\line(-1,0){.996}}
\put(26.02,10.96){\line(-1,0){.996}}
\put(24.03,10.96){\line(-1,0){.996}}
\put(22.04,10.96){\line(-1,0){.996}}
\put(20.05,10.96){\line(-1,0){.996}}
\put(18.06,10.96){\line(-1,0){.996}}
\put(16.06,10.96){\line(-1,0){.996}}
\put(14.07,10.96){\line(-1,0){.996}}
\put(12.08,10.96){\line(-1,0){.996}}
%\end
%\dottedline(5,33)(64,33)
\multiput(4.96,32.96)(.98333,0){61}{{\rule{.5pt}{.5pt}}}
%\end
\put(36,33){\circle*{1.06}}
\put(58,23){\line(0,1){5}}
\put(64,25){\line(0,1){11}}
\put(64,33){\circle*{1.12}}
\put(57,12){\line(0,1){4}}
%\dashline{1}(57,12)(44,12)
\put(56.96,11.96){\line(-1,0){.929}}
\put(55.1,11.96){\line(-1,0){.929}}
\put(53.24,11.96){\line(-1,0){.929}}
\put(51.38,11.96){\line(-1,0){.929}}
\put(49.53,11.96){\line(-1,0){.929}}
\put(47.67,11.96){\line(-1,0){.929}}
\put(45.81,11.96){\line(-1,0){.929}}
%\end
%\dashline{1}(58,23)(50,23)
\put(57.96,22.96){\line(-1,0){.889}}
\put(56.18,22.96){\line(-1,0){.889}}
\put(54.4,22.96){\line(-1,0){.889}}
\put(52.62,22.96){\line(-1,0){.889}}
\put(50.84,22.96){\line(-1,0){.889}}
%\end
%\dashline{1}(64,25)(58,25)
\put(63.96,24.96){\line(-1,0){.857}}
\put(62.24,24.96){\line(-1,0){.857}}
\put(60.53,24.96){\line(-1,0){.857}}
\put(58.81,24.96){\line(-1,0){.857}}
%\end
\thicklines
\put(54,33){\vector(0,-1){10}}
\put(31.5,33){\vector(0,-1){7}}
\put(40,33){\vector(0,-1){22}}
\put(4,33){\makebox(0,0)[cc]{$\tau$}}
\put(25,34){\makebox(0,0)[cc]{$x_1$}}
\put(34,34){\makebox(0,0)[cc]{$x_2$}}
\put(48,34){\makebox(0,0)[cc]{$x_3$}}
\put(62,34){\makebox(0,0)[cc]{$x_4$}}
\put(33,29){\makebox(0,0)[cc]{$c_1$}}
\put(41.5,24){\makebox(0,0)[cc]{$c_2$}}
\put(55.5,29){\makebox(0,0)[cc]{$c_3$}}
\put(27,33){\circle*{1.03}}
\put(50,33){\circle*{1.06}}
\end{picture}

\caption{Illustration of a chronological tree showing the durations $c_1, c_2, c_3$ elapsed since \emph{coalescence} for each of the three consecutive pairs $(x_1, x_2), (x_2, x_3)$ and $(x_3, x_4)$ of the $\Xi_\tau=4$ individuals alive at level $\tau$.
%The \emph{heights} (generations in the discrete tree) of points $x_1, x_2, x_3, x_4$ are respectively 3, 3, 2, 4.
}
\label{fig : smalltreecoal}
\end{figure}

\begin{thm}
\label{thm : coalescent PP}
Conditional on $\Xi_\tau=n\ge 1$, let
$$
a_i:=\left\{
\begin{array}{ccl}
p_2(x_{i}\wedge x_{i+1}) & \mbox{ if } & i\in\{1,\ldots, n -1\}\\
0 & \mbox{ if } & i=n
\end{array}
\right.
$$
Then under $\PP_\chi (\cdot \mid \Xi_\tau\not=0)$, $(a_i;1\le i \le \Xi_\tau)$ is a sequence of i.i.d. r.v. stopped at its  first 0, whose common distribution is that of $A:=\inf Y_t$, where $Y$ is the Lévy process with Laplace exponent $\psi$ started at $\tau$ and killed upon exiting $(0,\tau]$. 

In particular, the duration $C$ elapsed since coalescence between two consecutive individuals (if any) has
$$
\PP(C\le \sigma)=\PP(A>\tau - \sigma \mid A \not=0) =\frac{1-1/W(\sigma)}{1-1/W(\tau)}\qquad 0\le\sigma\le\tau.
$$
Furthermore, the coalescence level between $x_j(\tau)$ and $x_k(\tau)$ ($j\le k$) is given by
$$
p_2(x_j(\tau)\wedge x_k(\tau)) = \min \{a_i : j\le i <k\}.
$$
\end{thm}
Coalescence levels can be seen on Fig. \ref{fig : smalltreecoal}.

\begin{rem}
\label{rem : Pop}
Taking $\Lambda(dz) = be^{-bz}dz$, one can recover Lemma 3 in \cite{P}. Namely, since an elementary calculation yields $W(x)=1+bx$,  
$$
\PP(C\in d\sigma) =\frac{1+b\tau}{\tau}\frac{1}{(1+b\sigma)^2}\qquad \sigma\le\tau.
$$
\end{rem}
For the proof of Theorem \ref{thm : coalescent PP}, we will need the following claim on chronological trees.

\paragraph{Claim.} For any chronological tree $\TT$, and any three points $x,y,z\in\TT$ such that $x\le y\le z$, then $x\wedge z \prec x\wedge y$ and $x\wedge z \prec y\wedge z$. In addition, 
\begin{equation}
\label{eqn : wedges}
x\wedge z \in\{ x\wedge y, y\wedge z\}.
\end{equation}

\paragraph{Proof.}
Recall from Theorem \ref{thm : proprietes jccp}(ii) that $y\prec x$ iff $\hat{t}\le s\le t$. Now let $x\le y\le z$, and set $a=\varphi(x)$, $b=\varphi(y)$, $c=\varphi(z)$, $s=\varphi(x\wedge y)$, and $t=\varphi(y\wedge z)$. Then since $x\wedge y\prec x,y$, and $y\wedge z \prec y,z$, we get
$$
\hat{s}\le a\le b\le s \quad \mbox{ and }\quad \hat{t}\le b\le c\le t. 
$$
First, we show that $x\wedge z \prec y$. Indeed, writing $r=\varphi(x\wedge z)$, we have $\hat{r}\le a\le c\le r$, but since $a\le b\le c$, we get $\hat{r}\le b\le r$, which is exactly $x\wedge z \prec y$. Then notice that since $x\wedge z \prec y$ and $x\wedge z \prec x, z$, we always have $x\wedge z \prec y\wedge z$ and $x\wedge z \prec x\wedge y$. This is the first assertion of the Claim.

Next, we show that \eqref{eqn : wedges}$\Leftrightarrow$ \eqref{eqn : equiv wedges}, where
\begin{equation}
\label{eqn : equiv wedges}
x\wedge y\prec z\quad \mbox{ or }\quad y\wedge z\prec x.
\end{equation}
Assume \eqref{eqn : equiv wedges}. If $x\wedge y\prec z$, then since $x\wedge y\prec x$, we get $x\wedge y\prec x\wedge z$; on the other hand, thanks to the first assertion of the Claim, $x\wedge z \prec x\wedge y$, so that $x\wedge z = x\wedge y$. The same reasoning shows that if $y\wedge z\prec x$, then $x\wedge z = y\wedge z$. Thus we have shown \eqref{eqn : equiv wedges} $\Rightarrow$ \eqref{eqn : wedges}, and the converse implication is straightforward.

Now observe that \eqref{eqn : equiv wedges} is equivalent to 
$$
\hat{s}\le c\le s \quad \mbox{ or }\quad \hat{t}\le a\le t. 
$$
To prove \eqref{eqn : wedges}, we assume that the last assertion does not hold, and show that this will bring a contradiction. Indeed, if the last display does not hold, then $c>s$ and $a>t$, so that
$$
\hat{s}\le a<\hat{t}\le b\le s<c\le t.
$$
In particular, we can extract from the last display the following two double inequalities $\hat{s}\le \hat{t}\le s$ as well as $\hat{t}\le s\le t$, which reads also
$$
\varphi^{-1}(t)\prec \varphi^{-1}(s) \prec \varphi^{-1}(\hat{t}).
$$
Writing $\varphi^{-1}(s) = (u,\sigma)$ and $\varphi^{-1}(t) = (v,\tau)$, and recalling from Theorem \ref{thm : proprietes jccp}(i) that $\varphi^{-1}(\hat{s}) = (u,\omega(u))$ and $\varphi^{-1}(\hat{t}) = (v,\omega(v))$, we see that $(u,\sigma)$ is in the segment $[(v,\tau) ,(v,\omega(v))]$, from which we conclude that $u=v$. Therefore, $\hat{s}=\hat{t}$, which contradicts the fact that $\hat{s}\le a<\hat{t}$. \hfill$\Box$

\paragraph{Proof of Theorem \ref{thm : coalescent PP} .}
Let us start with the last assertion of the theorem.
A straightforward consequence of the Claim is that for any $x\le y\le z$, we have $p_2(x\wedge z)=\min(p_2(x\wedge y), p_2(y\wedge z))$. Indeed, from $x\wedge z \prec x\wedge y$ and $x\wedge z \prec y\wedge z$, we get $p_2(x\wedge z)\le\min(p_2(x\wedge y), p_2(y\wedge z))$, and the converse inequality stems from \eqref{eqn : wedges}. The last assertion of the theorem then follows from a recursive application of this property. 

Next, recall from Theorem \ref{thm : jccp} the following recursive definition: $t_0=0$ and $t_{i+1}=\inf\{t> t_i : \Xtaut \in\{0,\tau\}\}$. Then abusing notation (confounding those for $\TT$ and $C_\tau(\TT)$), we have $t_i=\varphi (x_i)$, and by Theorem \ref{thm : proprietes jccp}(iii), we get that for any $1\le i \le \Xi_\tau$, $a_i$ is given by
$$
a_i=\inf_{t_i\le t\le t_{i+1}} \Xtaut.
$$ 
But thanks to Theorem \ref{thm : jccp}, under $\PP_\chi (\cdot \mid \Xi_\tau\not=0)$, the killed paths $e_i:= (\Xtauf{t_i+t}, 0\le t< t_{i+1}-t_i)$, $i\ge 1$, form a sequence of i.i.d. excursions, distributed as the Lévy process $Y$ started at $\tau$, killed upon exiting $(0,\tau]$, ending at the first excursion hitting 0 before $(\tau,+\infty)$. It is then straightforward that $(a_i; 1\le i \le \Xi_\tau)$ is a sequence of i.i.d. r.v., stopped at its first 0, and distributed as $A$, where, thanks to the Preliminaries on Lévy processes,
$$
\PP(A\le \sigma)=P_\tau (Y \mbox{ exits } (\sigma,\tau] \mbox{ at the bottom}) = W(0)/W(\tau-\sigma)= 1/W(\tau-\sigma)\qquad \sigma \in[0,\tau].
$$
The fact that $W(0)=1$ can be deduced from a Tauberian theorem, checking that $\psi(\lbd)/\lbd$ converges to 1 as $\lbd\tendinfty$.
In particular, the duration $C=\tau-A$ elapsed since coalescence between two consecutive individuals has
$$
\PP(C\le \sigma)=%\PP(A>\tau-\sigma\mid A \not=0) =   
\frac{\PP(A>\tau-\sigma)}{\PP( A \not=0)}=\frac{1-1/W(\sigma)}{1-1/W(\tau)} ,
$$
which completes the proof.
\hfill$\Box$

\subsection{Ages and residual lifetimes}

We hope that the last subsection has convinced the reader that multiple uses of Theorem \ref{thm : jccp} can be made by applying  standard results on Lévy processes to one's favourite question of population biology. Hereafter, we give a simple example, in the form of a statement on ages and residual lifetimes of living individuals. The proof of this statement is left to the reader, since it relies on the same analysis as done in that of Theorem \ref{thm : coalescent PP}.\\
\\
Fix $\tau>0$. As in the previous subsection, we let $(x_i(\tau);1\le i \le \Xi_\tau)$ denote the ranked points $x_1\le x_2 \le \cdots$ of $\TT$ such that $p_2(x_i)=\tau$. Set also $u_i:=p_1(x_i)$ the $i$-th individual of $\cal T$ who is alive at level $\tau$. Then for $1\le i \le \Xi_\tau$,
$$
A_i:=\tau-\alpha(u_i)\quad\mbox{ and }\quad R_i:=\omega(u_i) -\tau
$$
define respectively the \emph{age} and \emph{residual lifetime} of $u_i$ at `time' $\tau$.

\begin{prop}
\label{prop : ages}
Under $\PP_\chi$, conditional on $\Xi_\tau=n\ge 2$, the individuals $(u_i;2\le i\le n)$ alive at time $\tau$ except $u_1$, have i.i.d. ages and residual lifetimes, whose common distribution is independent of $\chi$ and equal to that of a pair $(A_\tau, R_\tau)$ such that
$$
\PP(A_\tau\in dx, R_\tau\in dy) = P_0 (-Y_{T_{(0,\infty]}-}\in dx, Y_{T_{(0,\infty]}} \in dy\mid T_{(0,\infty]}<T_{-\tau}) \qquad 0<x<\tau, y>0,
$$
where $Y$ is the spectrally positive Lévy process with Laplace exponent $\psi$.
\end{prop}

\subsection{Properties of the Crump--Mode--Jagers process}

\begin{figure}[ht]

\unitlength 1.8mm % = 5.69pt
\linethickness{0.4pt}
\ifx\plotpoint\undefined\newsavebox{\plotpoint}\fi % GNUPLOT compatibility
\begin{picture}(52,40)(-15,0)
\put(3,3){\vector(0,1){37}}
\put(2,14){\vector(1,0){50}}
\put(6.63,39.13){\makebox(0,0)[cc]{$\psi(\lbd)$}}
\put(47.25,16.75){\makebox(0,0)[cc]{$\lbd$}}
\put(4.88,15.63){\makebox(0,0)[cc]{$0$}}
\put(10.25,29.63){\makebox(0,0)[cc]{$\psi'(0^+)<0$}}
\put(4.75,1.75){\makebox(0,0)[cc]{$-b$}}
%\emline(3,3.5)(43,38)
\multiput(3,3.5)(.024434942,.0210751374){1637}{\line(1,0){.024434942}}
%\end
\qbezier(3,14)(3.75,4.38)(40,35.75)
\put(16.63,11.5){\makebox(0,0)[cc]{$b$}}
\put(10.5,15.75){\makebox(0,0)[cc]{$\eta$}}
\end{picture}

\caption{Graph of the Laplace exponent $\psi$ of the Lévy process $Y$ (in the finite case). On this example is shown a supercritical exponent, with positive largest root $\eta$, such that $\psi(0)=0$ (i.e. $\Lambda(\{+\infty\})=0$). The root $\eta$ is the  Malthusian parameter of the CMJ process.
In the finite case, note that because $\psi$ is convex, and $\psi(\lbd)\ge \lbd -b$, one has $\eta <b$.}
\label{fig : exposantLaplace}
\end{figure}
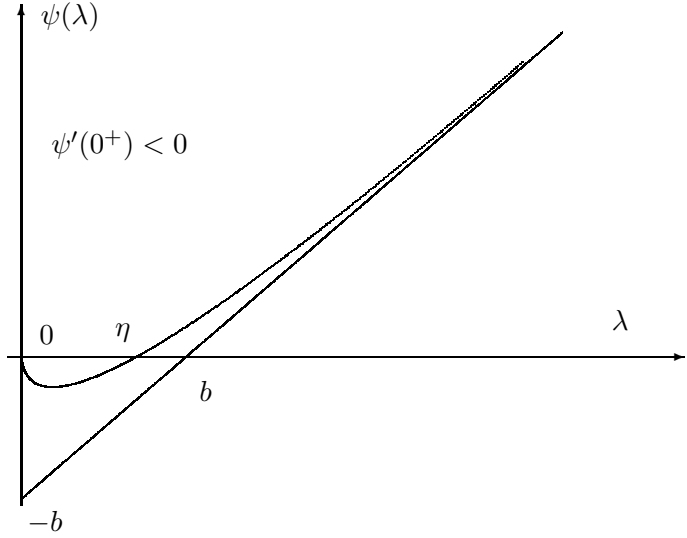

\begin{prop}
\label{prop : marginal}
The probability of extinction is $\PP_\chi(\mathrm{Ext})= e^{-\eta\chi}$.
In addition, the one-dimensional marginal of $\Xi$ is given by
$$
\PP_\chi (\Xi_\tau=0)=W(\tau-\chi)/W(\tau),
$$
and conditional on being nonzero, $\Xi_\tau$ has a \emph{geometric distribution} with success probability  $1/W(\tau)$.
In particular, $\EE_\chi(\Xi_\tau\mid\Xi_\tau \not=0 )= W(\tau)$.
\end{prop}

\paragraph{Proof.}
The value of the extinction probability stems straighforwardly from Proposition \ref{prop : properties Jirina} and Theorem \ref{thm : Jirina}.

Since $\Xi_\tau=0$ iff the first excursion of $\Xtau$ hits 0 before $(\tau,\infty)$, Theorem \ref{thm : jccp}(i) entails
$$
\PP_\chi (\Xi_\tau=0) = P_\chi(T_0< T_{(\tau,+\infty)})=W(\tau-\chi)/W(\tau) .
$$
Now thanks to Theorem \ref{thm : jccp}(i) again, conditional on being nonzero, $\Xi_\tau$ is the length of a sequence of i.i.d. excursions of $Y$ starting from $\tau$ and killed upon exiting $(0,\tau]$, stopped at the first one that exits it from the bottom. As a consequence, the conditional distribution of $\Xi_\tau$ is geometric with success probability 
$$
P_\tau(T_0< T_{(\tau,+\infty)})=1/W(\tau).
$$
The computation of $\EE_\chi(\Xi_\tau\mid\Xi_\tau \not=0)$ is then straightforward.
\hfill$\Box$

\begin{prop}
\label{prop : cond  CMJ ext}
Set $\PP^\natural:=\PP(\cdot\mid \mathrm{Ext})$. Under $\PP^\natural$, the (supercritical) splitting tree has the same law as the subcritical splitting tree with lifespan measure $e^{-\eta r} \Lambda(dr)$. In particular, in the finite case, its birth rate equals $b-\eta$.
\end{prop}

\paragraph{Proof.}
Thanks to Theorem \ref{thm : jccp}(ii), we know that under $\PP^\natural$, the JCCP $X$ is a spectrally positive Lévy process with Laplace exponent $\psi$ conditioned on hitting 0. Now it is known \cite{B} that this conditioned Lévy process is a spectrally positive Lévy process with Laplace exponent $\psi^\natural$, where $\psi^\natural (\lbd) =\psi(\lbd+\eta)$. It is straightforward to obtain the following equality
$$
 \psi^\natural (\lbd)= \lbd- \intgen e^{-\eta r}(1-e^{-\lbd r})\Lambda(dr)\qquad\lbd \ge 0,
$$
which proves that under $\PP^\natural$, $X$ is the JCCP of a splitting tree with lifespan measure $e^{-\eta r} \Lambda(dr)$. It is indeed a subcritical splitting tree, since $\psi^\natural(0)=\psi(\eta)=0$ and $\psi^{\natural\prime}(0^+)=\psi'(\eta)>0$.
Another elementary calculation shows that in the finite case, the birth rate $b^\natural=\intgen e^{-\eta r} \Lambda(dr)$ of this splitting tree equals $b-\eta$.
%, so that
%$$
%b^\star = b-\intgen (1-e^{-\eta r}) \Lambda(dr)=b - (\eta-\psi(\eta)) = b-\eta,
%$$ 
%since $\psi(\eta) =0$, which ends the proof.
\hfill$\Box$

\begin{prop} 
\label{prop : asymptotic CMJ}
The asymptotic behaviour of $\Xi$ is as follows.\\
\indent
(i)
(Yaglom's distribution)
In the subcritical case, 
$$
\lim_{\tau\tendinfty}\ \PP(\Xi_\tau = n \mid \Xi_\tau \not=0) = m^{n-1}(1-m) \qquad n\ge 1. 
$$
\indent
(ii)
In the critical case, provided that $\int^\infty r^2\Lambda(dr) <\infty$,
$$
\lim_{\tau\tendinfty}\ \PP(\Xi_\tau/\tau > x \mid \Xi_\tau \not=0) = \exp(-\psi''(0+)\ x/2)\qquad x\ge0. 
$$
\indent
(iii)
In the supercritical case, conditional on $\mathrm{Ext}^c$,
$$
\lim_{\tau\tendinfty}\ e^{-\eta\tau}\ \Xi_\tau \stackrel{\cal L}{=} \xi,
$$
where $\xi$ is an exponential variable with parameter $\psi'(\eta)$. 
\end{prop}
\begin{rem}
\label{rem : a.s. convergence} \label{61}
Whereas the convergence stated in (iii) certainly holds solely in distribution when $\Lambda$ is infinite (because of Theorem \ref{thm : exceptional points}), there is actually pathwise convergence in the finite case, provided that $\int^\infty r \log(r)\Lambda(dr) <\infty$ \cite{Ne}.
Thanks to this result, $\eta$ is seen to be the so-called \emph{Malthusian parameter}.
\end{rem}

\paragraph{Proof.}
We know from Proposition \ref{prop : marginal} that conditional on being nonzero, $\Xi_\tau$ is geometric  with success probability  $1/W(\tau)$. 

In the subcritical case ($\psi(0)=0$ and $\psi'(0^+)=1-m>0$), $\psi(\lbd)/\lbd$ converges to $1-m$ as $\lbd$ vanishes, so a Tauberian theorem entails that $W(\tau)$ converges to $1/(1-m)$ as $\tau\tendinfty$. This proves (i).

In the critical case ($\psi(0)=0$ and $\psi'(0^+)=0$), $\psi(\lbd)/\lbd^2$ converges to $\psi''(0^+)/2$ as $\lbd$ vanishes, so a Tauberian theorem entails that $W(\tau)/\tau$ converges to $2/\psi''(0^+)$ as $\tau\tendinfty$. This proves (ii).

As for (iii), set $W^\natural(x):=W(x)e^{-\eta x}$. Then the Laplace transform of $W^\natural$ is 
$$
\intgen W^\natural (x) e^{-\lbd x}\, dx = \frac{1}{\psi(\lbd+\eta)}\qquad \lbd\ge 0,  
$$
which is equivalent to $1/\psi'(\eta)\lbd$ as $\lbd$ vanishes. As previously, we deduce that $W^\natural(\tau)$ converges to $1/\psi'(\eta)$ as $\tau\tendinfty$, which reads
$$
\lim_{\tau\tendinfty}W(\tau)e^{-\eta \tau}=\frac{1}{\psi'(\eta)} . 
$$
The convergence in distribution then follows from the asymptotic equality between $\PP(\Xi_\tau e^{-\eta \tau}\in\cdot \mid \Xi_\tau\not=0)$ and $\PP(\Xi_\tau e^{-\eta \tau}\in\cdot \mid \mbox{Ext}^c)$. 
\hfill$\Box$\\
\\
Actually, we can go a little further than the previous proposition, by distinguishing between points whose descendance is either finite or infinite. Recall that finiteness here refers to the first projection of the tree (discrete part), and set \label{62}
$$
\Xi_\tau^\infty:=\mbox{Card} \{x\in\TT: p_2(x)=\tau, \theta(x)\mbox{ is infinite}\} ,
$$
$$
\Xi_\tau^f:=\mbox{Card} \{x\in\TT: p_2(x)=\tau, \theta(x)\mbox{ is finite}\} .
$$
In particular, $\Xi_\tau^\infty+\Xi_\tau^f=\Xi_\tau$. 
\begin{prop}
\label{prop : }
Set $p:=\psi'(\eta)\le 1$. In the supercritical case, conditional on $\mathrm{Ext}^c$,
$$
\lim_{\tau\tendinfty}\ e^{-\eta\tau}\ (\Xi_\tau^\infty,\Xi_\tau^f) \,\ \stackrel{\cal L}{=} \,\ (p\xi, (1-p)\xi),
$$
where $\xi$ is an exponential variable with parameter $p$. In particular, $e^{-\eta\tau} \Xi_\tau^\infty$ converges in distribution to an exponential variable with parameter 1.
\end{prop}
\begin{rem}
In a work in preparation \cite{Lprep2}, we show that actually, $(\Xi_\tau^\infty;\tau\ge 0)$ is Markovian, and, more precisely, that it is a \emph{Yule process} with birth rate $\eta$. In particular, the convergence of $e^{-\eta\tau}\ \Xi_\tau^\infty$ to an exponential variable with parameter 1 is known to be \emph{pathwise}.
\end{rem}

\paragraph{Proof.}
We use the notation of Proposition \ref{prop : ages}.
Each individual $u_i$ alive at level $\tau$ has infinite descendance iff the splitting tree starting from $u_i\times(\tau, \omega(u_i))$ is infinite, which, conditionally on the residual lifetime $R_i:=\omega(u_i)-\tau$, occurs with probability $1-e^{-\eta R_i}$. Now thanks to Proposition \ref{prop : ages}, conditionally on $\Xi_\tau=n\ge 2$, the residual lifetimes of the individuals $(u_i;2\le i\le n)$ are i.i.d., all distributed as $R_\tau$. As a consequence, by the branching property, the conditional distribution of $(\Xi_\tau^\infty,\Xi_\tau^f)$, modulo the first individual $u_1$ (which is negligible as the population size goes to infinity), is that of $(B, n-1-B)$, where $B$ is a binomial variable with parameters $p_\tau$ and $n-1$, and $p_\tau:=E(1-e^{-\eta R_\tau})$. Conditional on non-extinction, the only change is that $B$ is conditioned to be greater than 1. Now since the population size goes to infinity as $\tau\tendinfty$, this last conditioning vanishes, but above all, $(\Xi_\tau^\infty,\Xi_\tau^f)/\Xi_\tau$ converges in probability to $(p_\infty, 1-p_\infty)$, with $p_\infty:=\lim_\tau p_\tau$. Therefore, it only remains to show that $p_\infty=p$, since we know from Proposition \ref{prop : asymptotic CMJ} that $e^{-\eta\tau}\Xi_\tau$ converges in law to $\xi$. Now recall that $p_\tau=E(1-e^{-\eta R_\tau})$, and observe that $R_\tau$ converges in distribution, as $\tau\tendinfty$, to the law of $Y_{T_{(0,\infty]}}$ under $P_0$. Summing up, we have to prove that
$$
p=E_0 \left(1-\exp\big(-\eta Y_{T_{(0,\infty]}}\big)\right).
$$
As a first step, we apply the compensation formula to the Poisson point process $(\Delta_s; s\ge 0)$ of jumps of $Y$. Denoting by $S_t$ the past supremum of $Y$ at time $t$, we get, for any $\lbd\ge 0$, 
\debeq
E_0\left(\exp\big(-\lbd Y_{T_{(0,\infty]}}\big)\right)
		&=&		E_0 \sum_{t: \Delta_t>0}\indic{S_{t-}=0}\exp\big(-\lbd (Y_{t-}+\Delta_t)\big)\indic{Y_{t-}+\Delta_t>0}\\
			&=&		E_0	\intgen dt\intgen \Lambda(dr)\, \indic{S_{t}=0}\exp\big(-\lbd (Y_{t}+r)\big)\indic{Y_{t}+r>0}\\
	%&=&		E_0 \intgen \Lambda(dr)\, e^{-\lbd r}\intgen dt\, \indicbis{S_{t}=0}e^{-\lbd Y_{t}}\indicbis{-Y_{t}<r}\\
		&=&		 \intgen \Lambda(dr) e^{-\lbd r}\int_0^r \mu_0(dx) e^{\lbd x} ,
\fineq
where we have set
$$
\mu_0(dx) := E_0\int_0^{T_{(0,\infty]}} dt\, \indic{-Y_t\in dx}
%= \intgen dt\, P_0 (S_t =0, -Y_t \in dx) 
\qquad x>0.
$$
In order to compute an expression for $\mu_0$, we also define
$$
\mu(dx):= E_0\intgen dt\, \indic{-Y_t\in dx} \qquad x>0.
$$
Because $Y$ has derivative a.s. equal to $-1$  when it hits $-x$, we deduce that $\mu(dx) = E_0(N_x)\, dx$, where $N_x$ is the number of passage times at $-x$. As a consequence, applying the strong Markov property at $T_{-x}$, we get (see Preliminaries)
$$
\mu(dx) = P_0(T_{-x}<\infty)E_{-x}(N_x)\, dx = E_0(N_0) e^{-\eta x} \,dx\qquad x>0.
$$
On the other hand, applying the strong Markov property at the successive hitting times of 0 under $P_0$, say $0=S_1<S_2\cdots<S_{N_0}<S_{N_0+1}=\infty$, we get 
$$
\mu(dx) = E_0\sum_{n= 1}^{N_0}\int_{S_n}^{S_{n+1}} dt\, \indic{-Y_t\in dx} = E_0(N_0) \mu_0(dx),
$$
so that, in conclusion, 
$$
\mu_0(dx) =  e^{-\eta x} \,dx\qquad x>0.
$$
Therefore, for any $\lbd\ge 0$,
$$
E_0\left(\exp\big(-\lbd Y_{T_{(0,\infty]}}\big)\right)
=  \intgen \Lambda(dr) e^{-\lbd r}\int_0^r  e^{(\lbd-\eta) x},
$$
so in particular,
$$
E_0\left(\exp\big(-\eta Y_{T_{(0,\infty]}}\big)\right)
= \intgen \Lambda(dr)\, r \, e^{-\eta r} = 1-\psi'(\eta), 
$$
which is the desired result.\hfill $\Box$

\paragraph{The Markovian (finite) case.}
There are two possibilities for $(\Xi_\tau;\tau\ge 0)$ to be Markovian. 

First, when $\Lambda$ is a Dirac mass at $\{\infty\}$ (wih mass $b$), $(\Xi_\tau;\tau\ge 0)$ is a pure-birth process (with rate $b$), and since $\psi(\lbd)=\lbd -b$, we get $W(x) = \exp(bx)$, so that the size of the population at `time' $\tau$ is (shifted) geometric with success probability $\exp(-b\tau)$.

Second, when $\Lambda$ is exponential, with parameter, say $d$, then $(\Xi_\tau;\tau\ge 0)$ is a birth--death process with birth rate $b$ and death rate $d$ (supercritical iff $b>d$), so that $\Lambda(dr) = bde^{-dr}\,dr$.
%In the previous subsection, Remark \ref{rem : Pop}, we have displayed some calculations in the critical case, when $d=b$. 
Then we have
$$
\psi(\lbd) = \lbd - \frac{b\lbd}{d+\lbd}\qquad \lbd\ge 0.
$$
In particular, $\eta= b-d$ and $\psi'(\eta) = 1-(d/b)$. The scale function (see also Remark \ref{rem : Pop}) is $W(x) = 1+bx$ in the critical case ($d=b$), and in all other cases,
$$
W(x) = \frac{d-be^{(b-d)x}}{d-b}\qquad x\ge 0.
$$
One recovers that conditional on being nonzero, the size of the population at `time' $\tau$ is geometric with success probability $(d-b)/(d-be^{(b-d)\tau})$. In the supercritical case, conditional on non-extinction, the fraction of individuals with \emph{finite} descendance converges in probability to $d/b$, whereas the number of individuals with \emph{infinite} descendance is a Yule process with birth rate $b-d$.

\begin{figure}[ht]

\unitlength 1.2mm % = 2.85pt

\linethickness{0.4pt}
\begin{picture}(34,130)(0,0)
\end{picture}
\begin{picture}(34,200)(0,0)
\thicklines
\put(4,196){\line(1,0){12}}
\thinlines
%\vector{dash}(13,196)(13,193)
\put(13,193){\vector(0,-1){.07}}\put(12.93,195.93){\line(0,-1){.75}}
\put(12.93,194.43){\line(0,-1){.75}}
%\end
\thicklines
\put(13,193){\line(1,0){8}}
\thinlines
%\vector{dash}(19,193)(19,190)
\put(19,190){\vector(0,-1){.07}}\put(18.93,192.93){\line(0,-1){.75}}
\put(18.93,191.43){\line(0,-1){.75}}
%\end
\thicklines
\put(19,190){\line(1,0){5}}
\thinlines
%\vector{dash}(17,193)(17,188)
\put(17,188){\vector(0,-1){.07}}\put(16.93,192.93){\line(0,-1){.833}}
\put(16.93,191.26){\line(0,-1){.833}}
\put(16.93,189.6){\line(0,-1){.833}}
%\end
\thicklines
\put(17,188){\line(1,0){4}}
\put(18,185){\line(1,0){5}}
\put(8,179){\line(1,0){6}}
\thinlines
%\vector{dash}(12,179)(12,176)
\put(12,176){\vector(0,-1){.07}}\put(11.93,178.93){\line(0,-1){.75}}
\put(11.93,177.43){\line(0,-1){.75}}
%\end
\thicklines
\put(12,176){\line(1,0){15}}
\thinlines
%\vector{dash}(25,176)(25,173)
\put(25,173){\vector(0,-1){.07}}\put(24.93,175.93){\line(0,-1){.75}}
\put(24.93,174.43){\line(0,-1){.75}}
%\end
\thicklines
\put(25,173){\line(1,0){3}}
\thinlines
%\vector{dash}(18,188)(18,185)
\put(18,185){\vector(0,-1){.07}}\put(17.93,187.93){\line(0,-1){.75}}
\put(17.93,186.43){\line(0,-1){.75}}
%\end
%\vector{dash}(26,173)(26,170)
\put(26,170){\vector(0,-1){.07}}\put(25.93,172.93){\line(0,-1){.75}}
\put(25.93,171.43){\line(0,-1){.75}}
%\end
\thicklines
\put(26,170){\line(1,0){4}}
\thinlines
%\vector{dash}(29,170)(29,167)
\put(29,167){\vector(0,-1){.07}}\put(28.93,169.93){\line(0,-1){.75}}
\put(28.93,168.43){\line(0,-1){.75}}
%\end
\thicklines
\put(29,167){\line(1,0){3}}
\thinlines
%\vector{dash}(28,170)(28,165)
\put(28,165){\vector(0,-1){.07}}\put(27.93,169.93){\line(0,-1){.833}}
\put(27.93,168.26){\line(0,-1){.833}}
\put(27.93,166.6){\line(0,-1){.833}}
%\end
\thicklines
\put(28,165){\line(1,0){2}}
\thinlines
%\vector{dash}(21,176)(21,168)
\put(21,168){\vector(0,-1){.07}}\put(20.93,175.93){\line(0,-1){.889}}
\put(20.93,174.15){\line(0,-1){.889}}
\put(20.93,172.37){\line(0,-1){.889}}
\put(20.93,170.6){\line(0,-1){.889}}
\put(20.93,168.82){\line(0,-1){.889}}
%\end
\thicklines
\put(21,168){\line(1,0){4}}
\thinlines
%\vector{dash}(15,176)(15,163)
\put(15,163){\vector(0,-1){.07}}\put(14.93,175.93){\line(0,-1){.929}}
\put(14.93,174.07){\line(0,-1){.929}}
\put(14.93,172.22){\line(0,-1){.929}}
\put(14.93,170.36){\line(0,-1){.929}}
\put(14.93,168.5){\line(0,-1){.929}}
\put(14.93,166.64){\line(0,-1){.929}}
\put(14.93,164.79){\line(0,-1){.929}}
%\end
\thicklines
\put(15,163){\line(1,0){8}}
\thinlines
%\vector{dash}(21,163)(21,160)
\put(21,160){\vector(0,-1){.07}}\put(20.93,162.93){\line(0,-1){.75}}
\put(20.93,161.43){\line(0,-1){.75}}
%\end
\thicklines
\put(21,160){\line(1,0){4}}
\thinlines
%\vector{dash}(19,163)(19,158)
\put(19,158){\vector(0,-1){.07}}\put(18.93,162.93){\line(0,-1){.833}}
\put(18.93,161.26){\line(0,-1){.833}}
\put(18.93,159.6){\line(0,-1){.833}}
%\end
\thicklines
\put(19,158){\line(1,0){4}}
\thinlines
%\vector{dash}(22,158)(22,155)
\put(22,155){\vector(0,-1){.07}}\put(21.93,157.93){\line(0,-1){.75}}
\put(21.93,156.43){\line(0,-1){.75}}
%\end
\thicklines
\put(22,155){\line(1,0){7}}
\thinlines
%\vector{dash}(25,155)(25,150)
\put(25,150){\vector(0,-1){.07}}\put(24.93,154.93){\line(0,-1){.833}}
\put(24.93,153.26){\line(0,-1){.833}}
\put(24.93,151.6){\line(0,-1){.833}}
%\end
%\vector{dash}(23,155)(23,148)
\put(23,148){\vector(0,-1){.07}}\put(22.93,154.93){\line(0,-1){.875}}
\put(22.93,153.18){\line(0,-1){.875}}
\put(22.93,151.43){\line(0,-1){.875}}
\put(22.93,149.68){\line(0,-1){.875}}
%\end
\thicklines
\put(23,148){\line(1,0){3}}
\thinlines
%\vector{dash}(16,163)(16,146)
\put(16,146){\vector(0,-1){.07}}\put(15.93,162.93){\line(0,-1){.944}}
\put(15.93,161.04){\line(0,-1){.944}}
\put(15.93,159.15){\line(0,-1){.944}}
\put(15.93,157.26){\line(0,-1){.944}}
\put(15.93,155.37){\line(0,-1){.944}}
\put(15.93,153.49){\line(0,-1){.944}}
\put(15.93,151.6){\line(0,-1){.944}}
\put(15.93,149.71){\line(0,-1){.944}}
\put(15.93,147.82){\line(0,-1){.944}}
%\end
\thicklines
\put(16,146){\line(1,0){5}}
\thinlines
%\vector{dash}(19,146)(19,143)
\put(19,143){\vector(0,-1){.07}}\put(18.93,145.93){\line(0,-1){.75}}
\put(18.93,144.43){\line(0,-1){.75}}
%\end
\thicklines
\put(19,143){\line(1,0){3}}
\thinlines
%\vector{dash}(10,179)(10,141)
\put(10,141){\vector(0,-1){.07}}\put(9.93,178.93){\line(0,-1){.974}}
\put(9.93,176.98){\line(0,-1){.974}}
\put(9.93,175.03){\line(0,-1){.974}}
\put(9.93,173.08){\line(0,-1){.974}}
\put(9.93,171.13){\line(0,-1){.974}}
\put(9.93,169.19){\line(0,-1){.974}}
\put(9.93,167.24){\line(0,-1){.974}}
\put(9.93,165.29){\line(0,-1){.974}}
\put(9.93,163.34){\line(0,-1){.974}}
\put(9.93,161.39){\line(0,-1){.974}}
\put(9.93,159.44){\line(0,-1){.974}}
\put(9.93,157.49){\line(0,-1){.974}}
\put(9.93,155.55){\line(0,-1){.974}}
\put(9.93,153.6){\line(0,-1){.974}}
\put(9.93,151.65){\line(0,-1){.974}}
\put(9.93,149.7){\line(0,-1){.974}}
\put(9.93,147.75){\line(0,-1){.974}}
\put(9.93,145.8){\line(0,-1){.974}}
\put(9.93,143.85){\line(0,-1){.974}}
\put(9.93,141.9){\line(0,-1){.974}}
%\end
\thicklines
\put(10,141){\line(1,0){3}}
\thinlines
%\vector{dash}(8,196)(8,179)
\put(8,179){\vector(0,-1){.07}}\put(7.93,195.93){\line(0,-1){.944}}
\put(7.93,194.04){\line(0,-1){.944}}
\put(7.93,192.15){\line(0,-1){.944}}
\put(7.93,190.26){\line(0,-1){.944}}
\put(7.93,188.37){\line(0,-1){.944}}
\put(7.93,186.49){\line(0,-1){.944}}
\put(7.93,184.6){\line(0,-1){.944}}
\put(7.93,182.71){\line(0,-1){.944}}
\put(7.93,180.82){\line(0,-1){.944}}
%\end
%\vector{dash}(14,193)(14,183)
\put(14,183){\vector(0,-1){.07}}\put(13.93,192.93){\line(0,-1){.909}}
\put(13.93,191.11){\line(0,-1){.909}}
\put(13.93,189.29){\line(0,-1){.909}}
\put(13.93,187.48){\line(0,-1){.909}}
\put(13.93,185.66){\line(0,-1){.909}}
\put(13.93,183.84){\line(0,-1){.909}}
%\end
\thicklines
\put(14,183){\line(1,0){3}}
\thinlines
%\vector{dash}(16,183)(16,181)
\put(16,181){\vector(0,-1){.07}}\put(15.93,182.93){\line(0,-1){.667}}
\put(15.93,181.6){\line(0,-1){.667}}
%\end
\thicklines
\put(16,181){\line(1,0){2}}
\thinlines
\put(4,199){\vector(1,0){30}}
%\dashline(4,200)(4,135.87)
\put(3.93,199.93){\line(0,-1){.987}}
\put(3.93,197.96){\line(0,-1){.987}}
\put(3.93,195.98){\line(0,-1){.987}}
\put(3.93,194.01){\line(0,-1){.987}}
\put(3.93,192.04){\line(0,-1){.987}}
\put(3.93,190.06){\line(0,-1){.987}}
\put(3.93,188.09){\line(0,-1){.987}}
\put(3.93,186.12){\line(0,-1){.987}}
\put(3.93,184.14){\line(0,-1){.987}}
\put(3.93,182.17){\line(0,-1){.987}}
\put(3.93,180.2){\line(0,-1){.987}}
\put(3.93,178.22){\line(0,-1){.987}}
\put(3.93,176.25){\line(0,-1){.987}}
\put(3.93,174.28){\line(0,-1){.987}}
\put(3.93,172.3){\line(0,-1){.987}}
\put(3.93,170.33){\line(0,-1){.987}}
\put(3.93,168.36){\line(0,-1){.987}}
\put(3.93,166.38){\line(0,-1){.987}}
\put(3.93,164.41){\line(0,-1){.987}}
\put(3.93,162.44){\line(0,-1){.987}}
\put(3.93,160.47){\line(0,-1){.987}}
\put(3.93,158.49){\line(0,-1){.987}}
\put(3.93,156.52){\line(0,-1){.987}}
\put(3.93,154.55){\line(0,-1){.987}}
\put(3.93,152.57){\line(0,-1){.987}}
\put(3.93,150.6){\line(0,-1){.987}}
\put(3.93,148.63){\line(0,-1){.987}}
\put(3.93,146.65){\line(0,-1){.987}}
\put(3.93,144.68){\line(0,-1){.987}}
\put(3.93,142.71){\line(0,-1){.987}}
\put(3.93,140.73){\line(0,-1){.987}}
\put(3.93,138.76){\line(0,-1){.987}}
\put(3.93,136.79){\line(0,-1){.987}}
%\end
\thicklines
\put(25,150){\line(1,0){2}}
\thinlines
%\vector{dash}(27,155)(27,153)
\put(27,153){\vector(0,-1){.07}}\put(26.93,154.93){\line(0,-1){.667}}
\put(26.93,153.6){\line(0,-1){.667}}
%\end
\thicklines
\put(27,153){\line(1,0){3}}
\thinlines
%\vector{dash}(28,153)(28,151)
\put(28,151){\vector(0,-1){.07}}\put(27.93,152.93){\line(0,-1){.667}}
\put(27.93,151.6){\line(0,-1){.667}}
%\end
\thicklines
\put(28,151){\line(1,0){3}}

\put(4,129){\vector(1,0){30}}
%\emline(16,129)(13,126)
\multiput(16,129)(-.0337079,-.0337079){89}{\line(0,-1){.0337079}}
%\end
%\emline(21,126)(19,124)
\multiput(21,126)(-.0333333,-.0333333){60}{\line(0,-1){.0333333}}
%\end
\put(24,124){\line(-1,-1){7}}
%\emline(21,117)(18,114)
\multiput(21,117)(-.0337079,-.0337079){89}{\line(0,-1){.0337079}}
%\end
\put(23,114){\line(-1,-1){9}}
%\emline(17,105)(16,104)
\multiput(17,105)(-.033333,-.033333){30}{\line(0,-1){.033333}}
%\end
\put(18,104){\line(-1,-1){10}}
%\emline(14,94)(12,92)
\multiput(14,94)(-.0333333,-.0333333){60}{\line(0,-1){.0333333}}
%\end
%\emline(27,92)(25,90)
\multiput(27,92)(-.0333333,-.0333333){60}{\line(0,-1){.0333333}}
%\end
%\emline(28,90)(26,88)
\multiput(28,90)(-.0333333,-.0333333){60}{\line(0,-1){.0333333}}
%\end
\put(26,88){\line(1,0){4}}
%\emline(30,88)(29,87)
\multiput(30,88)(-.033333,-.033333){30}{\line(0,-1){.033333}}
%\end
\put(29,87){\line(1,0){3}}
\put(32,87){\line(-1,-1){4}}
\put(28,83){\line(1,0){2}}
\put(30,83){\line(-1,-1){9}}
\put(21,74){\line(1,0){4}}
\put(25,74){\line(-1,-1){10}}
\put(15,64){\line(1,0){8}}
%\emline(23,64)(21,62)
\multiput(23,64)(-.0333333,-.0333333){60}{\line(0,-1){.0333333}}
%\end
\put(21,62){\line(1,0){4}}
\put(25,62){\line(-1,-1){6}}
\put(19,56){\line(1,0){4}}
%\emline(23,56)(22,55)
\multiput(23,56)(-.033333,-.033333){30}{\line(0,-1){.033333}}
%\end
\put(22,55){\line(1,0){7}}
%\emline(29,55)(27,53)
\multiput(29,55)(-.0333333,-.0333333){60}{\line(0,-1){.0333333}}
%\end
\put(27,53){\line(1,0){3}}
%\emline(30,53)(28,51)
\multiput(30,53)(-.0333333,-.0333333){60}{\line(0,-1){.0333333}}
%\end
\put(28,51){\line(1,0){3}}
\put(31,51){\line(-1,-1){6}}
\put(25,45){\line(1,0){2}}
\put(27,45){\line(-1,-1){4}}
\put(23,41){\line(1,0){3}}
\put(26,41){\line(-1,-1){10}}
\put(16,31){\line(1,0){5}}
%\emline(21,31)(19,29)
\multiput(21,31)(-.0333333,-.0333333){60}{\line(0,-1){.0333333}}
%\end
\put(19,29){\line(1,0){3}}
\put(22,29){\line(-1,-1){12}}
\put(10,17){\line(1,0){3}}
\put(13,17){\line(-1,-1){9}}
\put(4,130){\vector(0,-1){126}}
\put(13,126){\line(1,0){8}}
\put(19,124){\line(1,0){5}}
\put(17,117){\line(1,0){4}}
\put(18,114){\line(1,0){5}}
\put(14,105){\line(1,0){3}}
\put(16,104){\line(1,0){2}}
\put(8,94){\line(1,0){6}}
\put(12,92){\line(1,0){15}}
\put(25,90){\line(1,0){3}}
\end{picture}

\caption{ A chronological tree and the  associated jumping chronological contour process (JCCP), with jumps in solid line.}\label{fig : JCCP}
\end{figure}

\paragraph{Acknowledgements.}
I would like to thank Jean Bertoin, Thomas Duquesne, Grégory Miermont and Lea Popovic, for some nice and interactive discussions. I also want to acknowledge the clear-sightedness of an anonymous referee who pointed several awkward passages, and thereby allowed this work to attain its present form.

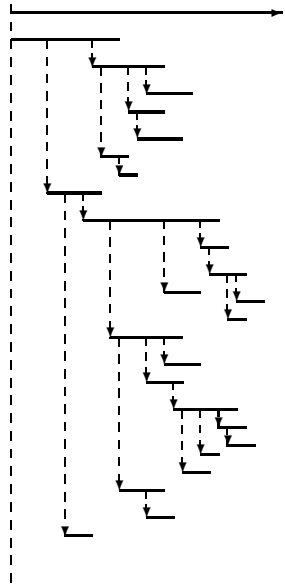
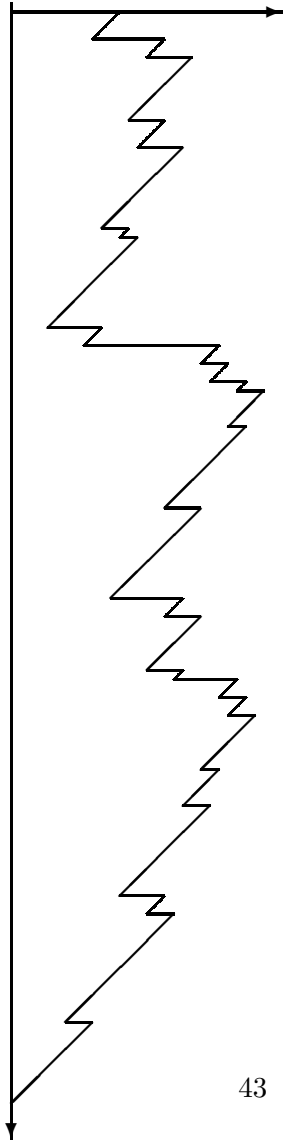
\end{document}